\documentclass[UTF-8,reqno]{amsart}
\usepackage{enumerate}
\usepackage{mhequ}
\usepackage[margin=1.2in]{geometry}
\usepackage{amssymb,url,color, booktabs,nccmath}
\usepackage{mathrsfs}
\usepackage{enumitem}
\usepackage{graphicx}
\usepackage{tikz}
\usepackage{amsthm}
\usepackage{color}
\usepackage[colorlinks=true]{hyperref}
\hypersetup{
    linkcolor=blue,          
    citecolor=red,        
    filecolor=blue,      
    urlcolor=cyan
}

\definecolor{darkergreen}{rgb}{0.0, 0.5, 0.0}

\numberwithin{equation}{section}

\newcommand{\be}{\begin{eqnarray}}
\newcommand{\ee}{\end{eqnarray}}
\newcommand{\ce}{\begin{eqnarray*}}
\newcommand{\de}{\end{eqnarray*}}
\newtheorem{theorem}{Theorem}[section]
\newtheorem{lemma}[theorem]{Lemma}
\newtheorem{remark}[theorem]{Remark}
\newtheorem{definition}[theorem]{Definition}
\newtheorem{proposition}[theorem]{Proposition}
\newtheorem{Examples}[theorem]{Example}
\newtheorem{corollary}[theorem]{Corollary}

\newtheorem*{theorem*}{Theorem}
\newtheorem*{remark*}{Remark}

\newcommand{\assign}{:=}
\newcommand{\cdummy}{\cdot}
\newcommand{\comma}{{,}}
\newcommand{\mathd}{\mathrm{d}}
\newcommand{\nobracket}{}
\newcommand{\nocomma}{}

\newcommand{\tmop}[1]{\ensuremath{\operatorname{#1}}}

\newcommand{\udiv}{\, \mathrm{div}}
\newcommand{\ud}{\,\mathrm{d}}

\def\eps{\varepsilon}

\def\p{\partial}

\def\[{{\Big[}}
\def\]{{\Big]}}
\def\<{{\langle}}
\def\>{{\rangle}}
\def\({{\Big(}}
\def\){{\Big)}}

\def\bx{{\mathbf{x}}}

\def\dif{{\mathord{{\rm d}}}}

\def\={&\!\!=\!\!&}

\def\cF{{\mathcal F}}

\def\cS{{\mathcal S}}

\def\mN{{\mathbb N}}

\def\mP{{\mathbb P}}

\def\mR{{\mathbb R}}

\def\mT{{\mathbb T}}

\def\1{{\mathbf{1}}}

\def\geq{\geqslant}
\def\leq{\leqslant}

\def\div{\mathord{{\rm div}}}

\def\eps{\varepsilon}

\def\p{\partial}

\def\[{{\Big[}}
\def\]{{\Big]}}
\def\<{{\langle}}
\def\>{{\rangle}}
\def\({{\Big(}}
\def\){{\Big)}}

\def\bx{{\mathbf{x}}}

\def\dif{{\mathord{{\rm d}}}}

\def\={&\!\!=\!\!&}
\def\bt{\begin{theorem}}
\def\et{\end{theorem}}
\def\bl{\begin{lemma}}
\def\el{\end{lemma}}
\def\br{\begin{remark}}
\def\er{\end{remark}}
\def\bx{\begin{Examples}}
\def\ex{\end{Examples}}
\def\bd{\begin{definition}}
\def\ed{\end{definition}}
\def\bp{\begin{proposition}}
\def\ep{\end{proposition}}
\def\bc{\begin{corollary}}
\def\ec{\end{corollary}}

\def\geq{\geqslant}
\def\leq{\leqslant}

\def\div{\mathord{{\rm div}}}

\def\<{\langle} \def\>{\rangle}

\allowdisplaybreaks

\begin{document}

\title{Gaussian fluctuations for interacting  particle systems with singular kernels }

\author{Zhenfu Wang}
\address[Z. Wang]{Beijing International Center for Mathematical Research, Peking University, Beijing 100871, China}
\email{zwang@bicmr.pku.edu.cn}

\author{Xianliang Zhao}
\address[X. Zhao]{ Academy of Mathematics and Systems Science,
	Chinese Academy of Sciences, Beijing 100190, China; Fakult\"at f\"ur Mathematik, Universit\"at Bielefeld, D-33501 Bielefeld, Germany}
\email{xzhao@math.uni-bielefeld.de}

\author{Rongchan Zhu}
\address[R. Zhu]{Department of Mathematics, Beijing Institute of Technology, Beijing 100081, China; Key Laboratory on MCAACI, Beijing, China
}
\email{zhurongchan@126.com}

\begin{abstract}
	We consider  the asymptotic behaviour of  the fluctuations for the empirical measures of interacting particle systems with singular kernels. We prove that the sequence of fluctuation processes  converges in distribution to a generalized Ornstein-Uhlenbeck process.
	
	Our result considerably extends classical results  to singular kernels, including the Biot-Savart law.  The result applies to the point vortex model approximating  the 2D incompressible Navier–Stokes equation and the 2D Euler equation.  We also obtain Gaussianity and optimal regularity of the limiting Ornstein-Uhlenbeck process.
	The method relies on  the martingale approach  and  the  Donsker-Varadhan variational formula, which transfers the uniform estimate to some exponential integrals. Estimation of those exponential integrals follows by  cancellations and combinatorics techniques and is of the type  of large deviation principle.
\end{abstract}


\date{\today}

\maketitle

\setcounter{tocdepth}{1}

\tableofcontents

\section{Introduction}

In this article, we consider interacting particle systems characterized by the following system of SDEs on the torus $\mathbb{T}^d$, $d\geq 2$,
\begin{equation}
	\mathd \nocomma X_i = \frac{1}{N}  \sum_{j \ne i} K (X_i - X_j) \mathd
	\nocomma t + F (X_i) \mathd t + \sqrt{2 \sigma_N} \mathd \nocomma B^i_t,
	\quad i = 1, ..., N, \label{equation pa}
\end{equation}
with random initial data $\{X_i(0)\}_{i=1}^N$. The collection $\{B^i_\cdot \}_{i=1}^N$ consists of $N$ independent $d$ dimensional Brownian motions  on a stochastic basis, i.e. $(\Omega,\mathcal{F},\mathbb{P})$ with a normal filtration $(\mathcal{F}_t)$, induced by the Laplacian operator on the torus,
 independent of $\{X_i(0)\}_{i=1}^N$.  The coefficient $\sigma_N \geq 0$ is a non-negative scalar for simplicity.  In this model, $X^N(t)\assign(X_1(t),...,X_N(t))\in (\mathbb{T}^d)^N $ represents the positions of particles, which are interacting through the kernel  $K$ and confined by the exterior force  $F$.

Many particle systems written in the canonical form \eqref{equation pa} or its  variant are now quite ubiquitous. Such systems are usually formulated by first-principle agent based models which are conceptually simple. For instance, in physics those particles $X_i$ can represent ions and electrons in plasma physics \cite{dobrushin1979vlasov}, or molecules in a fluid \cite{jabin2004identification} or even large scale galaxies \cite{Jean} in some cosmological models; in biological sciences, they typically model the collective behavior of animals or micro-organisms (for instance flocking, swarming and chemotaxis and other aggregation phenomena \cite{carrillo2014derivation}); in economics or social sciences particles are usually individual ``agents” or ``players” for instance in opinion dynamics \cite{friedkin1990social} or in the study of mean-field games \cite{lasry2007mean,huang2006large}. Motivation even extends to the analysis of large biological \cite{bossy2015clarification} or artificial \cite{mei2018mean} neural networks in neuroscience or in machine learning.

Under mild assumptions, it is well-known  
(see for instance \cite{mckeanpropagation,braun1977vlasov,dobrushin1979vlasov,osada1986propagation,sznitman1991topics,fournier2014propagation,jabin2018quantitative,serfaty2020mean,jabin2014review,bresch2020mean} and Section \ref{sec:i-m} for more details)  that  the empirical measure  $\mu_N(t):=\frac{1}{N}\sum^N_{i=1}\delta_{X_i(t)}$ of  the particle system \eqref{equation pa} converges to  the solution $\bar \rho (t) $ of the following deterministic mean-field  PDE
\begin{equation}
	\partial_t  \bar{\rho} = \sigma \Delta \bar{\rho}- \tmop{div} \nocomma ([F + K \ast \bar{\rho}]
	\bar{\rho}), \label{equation me}
\end{equation}
as $N \to \infty$, where $\sigma = \lim_{N \to \infty} \sigma_N$.  This is equivalent to the {\em propagation of chaos}, i.e.  the $k-$th marginal $\rho_{N, k}$ of the particle system \eqref{equation pa} will converge to the tensor product of the limit law $\bar \rho^{\otimes k}$  as $N$ goes to infinity, given for instance the i.i.d. initial data.   This law of large numbers type result  implies that the continuum model \eqref{equation me} is a suitable approximation  to the particle system \eqref{equation pa} when  $N$ is large, i.e. $\mu_N\approx \bar{\rho}+o(1).$

Inspired in particular by quantitative estimates of  propagation of chaos by Jabin and Wang \cite{jabin2018quantitative}, which is $\|\rho_{N, k} - \bar \rho^{\otimes k} \|_{L^\infty ([0, T], L^1)  } \leq C_T / \sqrt{N}$, we aim to  study the  central limit theorem of \eqref{equation pa}, which  provides a  better continuum approximation to \eqref{equation pa}.  More precisely,  we study
 the  limit of fluctuation measures around the mean-field law, which are  defined by
\begin{align}\label{def:eta}
\eta^N:=\sqrt{N}(\mu_N-\bar{\rho})=\frac{1}{\sqrt{N}}\sum_{i = 1}^N\left( \delta_{X_i}-\bar{\rho}\right) .\end{align}
In this article, we establish   that the fluctuation measure  $\eta^N$ converges in distribution as $N\rightarrow \infty$  to an infinite-dimensional continuous Gaussian process $\eta$ for a large class of particle systems \eqref{equation pa}. This implies that there exists a  continuum model $\eta$ such that
$$\mu_N\overset{\text{d}}{\approx} \bar{\rho}+\frac{1}{\sqrt{N}}\eta+o(\frac{1}{\sqrt{N}}),$$
where  $\overset{\text{d}}{\approx}$ means that the approximation holds in distribution.

\subsection{Assumptions}\label{sec:assumption}
To state our main results we first give the  framework in this article. Recall that the relative entropy $H(\mu \vert \nu)$ between probability measures $\mu$ and $\nu$ on a Polish space $E$  is defined by
\begin{equation*}
H(\mu \vert \nu):=\begin{cases} \int_E  \frac{\mathd \mu}{\mathd \nu}	\log\frac{\mathd \mu}{\mathd \nu} \mathd \nu & \text{ if } \mu\ll \nu, \\
\infty	& \text{ otherwise,}
\end{cases}
\end{equation*}
where $\frac{\mathd \mu}{\mathd \nu}$ is the Radon-Nikodym derivative of  $\mu$ with respect to $\nu$. Note that throughout this article, all the relative entropy is of the classical form. We will not normalize it as what have been done  for instance in \cite{jabin2018quantitative}.

Our assumptions are listed as follows.
\paragraph*{\bf{(A1)}-CLT for initial values.}  
There exists $\eta_0$, which  belongs to the space of tempered distributions $\cS'(\mT^d)$,  
such that the  sequence $\{\eta^N_0\}_{N\geq 1}$ converges in distribution to $\eta_0$ in $\cS'(\mT^d)$.  Here $\eta_{0}$ will be the initial data for our expected limit SPDE \eqref{equation
	li} below.
$\\$
{\bf{(A2)}-Regularity of the kernel. } The kernel $K : \mathbb{T}^d \rightarrow \mathbb{R}^d$, $d\geq 2$, satisfies one of  the following conditions
\begin{enumerate}
	\item $K$ is bounded.
	\item For each $x \in \mathbb{T}^d$, $K (x) = - K (- x)$ and $|x | \nocomma
	K (x) \in L^{\infty}$.
	\end{enumerate}
\paragraph*{\bf{(A3)-Uniform  relative entropy bound.}} Let $X^N(t)=(X_1(t),...,X_N(t))$ be a solution to the particle system \eqref{equation pa}, and let  $\rho_N(t)$ represent  the joint distribution of $X^N(t)$.  It holds that
\begin{equation*}
\sup_{ t \in [0, T]}\sup_NH(\rho_N(t)|\bar{\rho}_N(t))<\infty,
\end{equation*}
where $\bar{\rho}_N(t)$ denotes the tensor product $\bar{\rho}_t^{\otimes N}$. We may use $H_t(\rho_N|\bar{\rho}_N)$ to represent $H(\rho_N(t)|\bar{\rho}_N(t))$ for simplicity.

The global well-posedness of the limit equation \eqref{equation li} will be obtained by two different approaches, depending on the diffusion coefficient $\sigma$ is positive  or zero. Hence we distinguish the  extra assumptions  into the following two cases.

For the case when $\sigma >0$, in addition to the assumptions {\bf{(A1)-(A3)}}, we need the following extra assumption:
\paragraph*{\bf{(A4)-The case with non-vanishing diffusion}}
\begin{enumerate}
	\item$\sigma > 0$ and $| \sigma_N - \sigma | =\mathcal{O} \left( \frac{1}{N}
	\right)$.
	\item There exists some $\beta > d / 2$  such that $\bar{\rho} \in C ([0, T],
	C^{\beta}(\mathbb{T}^d))$ and $F \in C^{\beta}(\mathbb{T}^d)$, where $\bar{\rho}$ solves  equation \eqref{equation me} in the weak sense.
\end{enumerate}

{On the other hand}, for the case with vanishing diffusion,   besides {\bf{(A1)-(A3)}},  we require that
\paragraph*{\bf{(A5)-The case with vanishing diffusion}} The diffusion coefficients  and  the mean-field equation \eqref{equation me} satisfy that,
\begin{enumerate}
	\item$\sigma = 0$ and $| \sigma_N - \sigma | =\mathcal{O} \left( \frac{1}{N}
	\right)$.
	\item $\bar{\rho} \in C^1 ([0, T], C^{\beta+2}(\mT^d))$ and $F \in C^{\beta+1}(\mT^d)$ with $\beta>d/2$ where $\bar{\rho}$  solves \eqref{equation me}.
	\item $\div K\in L^1$.
\end{enumerate}


We make several remarks on our assumptions. Firstly, when $\{X_i(0)\}_{i\in \mathbb{N}}$ are i.i.d. with a common probability density function $\mu$,  which is the usual setting to study the fluctuations,  
	one can easily check that {\bf{(A1)}}  holds true. Indeed,  for each $\varphi\in C^{\infty}(\mT^d)$, we have
	\begin{align*}
		\left\langle \eta^N_0,\varphi\right\rangle =\frac{1}{\sqrt{N}}\sum_{i=1}^N\[\varphi(X_i(0))-\left\langle \varphi,\mu \right\rangle \]\xrightarrow{N\rightarrow \infty} \mathcal{N}\(0,\left\langle \varphi^2,\mu\right\rangle-\left\langle \varphi,\mu\right\rangle^2  \),
	\end{align*}
where $\mathcal{N}(0,a)$ denotes the centered Gaussian distribution on $\mathbb{R}$ with variance $a$. Hereafter we use the bracket $\langle \cdot  , \cdot  \rangle$ as a shorthand notation for integration. We also state a central limit theorem under an assumption on $H(\rho_{ N}(0)|\bar{\rho}_N(0))$ in Section \ref{sec:5.2}, where $\{X_i(0)\}_{i\in \mathbb{N}}$ can be neither independent nor identically  distributed.

	Assumption {\bf{(A2)}} on interaction kernels allows our framework to cover smooth kernels and some singular kernels, in particular the Biot-Savart kernel related to the vorticity formulation  of 2D Navier-Stokes/Euler equation on the torus. See Theorem \ref{th:main} and Section \ref{sec:exm} for more details.

 Assumption {\bf{(A3)}} seems quite nontrivial and demanding, but fortunately it has been established by Jabin and Wang  in \cite{jabin2018quantitative} for a quite large family of interacting kernels, including all the kernels satisfying {\bf{(A2)}}. Indeed, once we have that  the relative entropy between the joint distribution $\rho_N$ of the  interacting  particle system \eqref{equation pa} and the tensorized law $\bar{\rho}^{\otimes N}$ of the  mean-field  PDE \eqref{equation me} is  uniformly bounded with respect to $N$, then easily  the particle  system \eqref{equation pa} converges to the mean-field equation \eqref{equation me} with a rate  $\frac{C_T}{\sqrt{N}}$,  in the total variation norm or the Wasserstein metric.  More precisely, since all particles in \eqref{equation pa} are indistinguishable, the joint distribution $\rho_N $ is thus assumed to be  symmetric/exchangeable,   so is any $k$-marginal  distribution $\rho_{N, k}$  of $\rho_N$,  which is defined as
	\begin{equation*}
		\rho_{N, k}(t,x_1,...,x_k):=\int_{\mathbb{T}^{d(N-k)}}\rho_N(t,x_1,...,x_N)\mathd x_{k+1}...\mathd x_N.
	\end{equation*}
	Then  by the  sub-additivity  of relative entropy, in particular  $H(\rho_{N, k}|\bar{\rho}^{\otimes k})\leq \frac{k}{N} H(\rho_N \vert \bar{\rho}^{\otimes N})$  and  the classical  Csisz{\'a}r--Kullback--Pinsker inequality \cite[(22.25)]{villani2008optimalbook},  it follows that for fixed $k\in \mathbb{N}$,
	\begin{equation}\label{CKP}
		W_1\(\rho_{N, k}(t),\bar{\rho}^{\otimes k}(t)\)\lesssim	\|\rho_{N, k}(t)-\bar{\rho}^{\otimes k}(t)\|_{TV}  \leq  \sqrt{2H_t(\rho_{N, k}|\bar{\rho}^{\otimes k})}\lesssim \sqrt{\frac{k}{N}}\longrightarrow0,
	\end{equation}
	where $W_1(\cdot,\cdot)$ denotes the  $1$-Wasserstein distance, $\|\cdot\|_{TV}$ denotes the  total variation norm and the first inequality is guaranteed by \cite[Theorem 6.15]{villani2008optimalbook} since now $\mathbb{T}^d$ is compact.
	
\subsection{Main results}
Under the assumptions  {\bf{(A1)-(A3)}} and either {\bf{(A4)}} or  {\bf{(A5)}}, depending on $\sigma > 0$ or $\sigma =0$,   we establish that as $N \to \infty$,  the sequence of the fluctuation measures $(\eta_\cdot^N )$ converges  in  distribution to  the centered Gaussian  process $\eta$
solving the following stochastic PDE (SPDE)
\begin{equation}
	\partial_t \eta = \sigma \Delta \eta - \nabla \cdot (\bar{\rho} K \ast \eta)
	- \nabla \cdot (\eta K \ast \bar{\rho}) - \nabla \cdot (F \eta) - \sqrt{2
		\sigma} \nabla \cdot \left( \sqrt{\bar{\rho}} \xi \right) ,\qquad \eta(0)=\eta_0,\label{equation
		li}
\end{equation}
where $\eta_0$ is given in Assumption {\bf (A1)} and $\xi$ is vector-valued space-time white noise on $\mR^+\times\mT^d$,  i.e. a family of centered Gaussian random variables $\{\xi(h):h\in L^2(\mR^+\times \mT^d;\mR^d)\}$ such that $\mathbb{E}[|\xi(h)|^2]=\|h\|^2_{L^2(\mR^+\times \mT^d;\mR^d)}$, and $\bar{\rho}$ solves the mean-field equation \eqref{equation me}.  But when $\sigma =0$, the SPDE \eqref{equation li} becomes a deterministic PDE.

As a first step, we need a proper notion of solutions to the SPDE \eqref{equation li}.  When $\sigma>0$, it turns out to be  the martingale solutions defined as below.
\begin{definition}
	\label{def:mar}
	We call $\eta$ a martingale
	solution to the SPDE {\eqref{equation li}} on some stochastic basis $(\Omega, \mathcal{F},
	\mathcal{F}_t, \mathbb{P})$  if
	\begin{enumerate}
		\item$\eta$  is a continuous $(\mathcal{F}_t)$-adapted process with values
		in $H^{- \alpha - 2}$  and  $\eta\in L^2 ([0, T], H^{-
			\alpha}) $  for every $\alpha > d / 2$,  $\mathbb{P}$-almost surely.
		\item For each $\varphi \in C^{\infty} (\mathbb{T}^d)$ and $t \in [0, T]$, it
		holds that
		
		\begin{align}
			\langle \eta_t, \varphi \rangle - \langle \eta_0, \varphi \rangle =
			\int^t_0 \langle \sigma \Delta \varphi, \eta \rangle \mathd \nocomma s +
			\int^t_0 \langle \nabla \varphi, \bar{\rho} K \ast \eta +\eta K \ast \bar{\rho}+F \eta\rangle \mathd
			\nocomma s + \mathcal{M}_t (\varphi)
			,  \nonumber
		\end{align}
		where $\mathcal{M}$  is  a continuous $(\mathcal{F}_t)$-adapted centered Gaussian process
		with values in $H^{- \alpha - 1}$ for every $\alpha > d / 2$ and its
		covariance given by
		\begin{equation*}
			\mathbb{E} [\mathcal{M}_t (\varphi_1) \mathcal{M}_s (\varphi_2)] = 2
			\sigma \int^{s \wedge t}_0 \langle \nabla \varphi_1 \cdot\nabla \varphi_2,
			\bar{\rho}_r \rangle \mathd r,
		\end{equation*}
		for each $\varphi_1,\varphi_2\in C^{\infty}(\mT^d)$ and $s,t\in [0,T]$.
	\end{enumerate}
\end{definition}
\begin{remark}
	\begin{enumerate}
	\item The stochastic basis in Definition \ref{def:mar} might be different from the stochastic basis where the particle system \eqref{equation pa} lives.
	\item By Lemma \ref{lemma triebel} and Lemma \ref{lemma convolution} given in Appendix \ref{sec:appa},  $\bar{\rho} K \ast \eta$,  $\eta K \ast \bar{\rho}$, and $F \eta$ are all well-defined under Assumption {\bf{(A4)}}.
	\item The noise $\mathcal{M}$ is equivalent to be  characterized as:  for each $\varphi\in C^{\infty}$, $\mathcal{M}(\varphi)$ is a continuous $(\mathcal{F}_t)$-adapted  martingale with quadratic variation given by
	\begin{equation*}
		\mathbb{E} \[|\mathcal{M}_t (\varphi)|^2\] = 2
		\sigma \int^{ t}_0 \langle |\nabla \varphi|^2,
		\bar{\rho}_r \rangle \mathd r.
	\end{equation*}
		\end{enumerate}
\end{remark}

Similarly, when $\sigma=0$, the equation \eqref{equation li} actually becomes a deterministic  PDE.  We  define solutions to this first order PDE as follows.
\begin{definition} Given that $\sigma =0$, 	we call $\eta$ a solution to the  PDE \eqref{equation li}  with random initial data $\eta_{0}$,  if
\begin{enumerate}
	\item $ \eta \in L^2 ([0, T], H^{-
		\alpha}) \cap C ([0, T], H^{- \alpha - 2})$ for every
	$\alpha > d / 2$ almost surely.
	
	\item For each $\varphi \in C^{\infty} (\mathbb{T}^d)$ and $t \in [0, T]$, it
	holds that 	
	\begin{align}
		\langle \eta_t, \varphi \rangle = \langle \eta_0, \varphi \rangle +
		\int^t_0 \langle \nabla \varphi, \bar{\rho} K \ast \eta +\eta K \ast \bar{\rho}+F \eta\rangle \mathd
		\nocomma s . \nonumber
	\end{align}
\end{enumerate}
\end{definition}


Our first main result gives   the convergence of fluctuation measures when the diffusion coefficient $\sigma $ is positive (which  may be generalized to the case with a non-degenerate coefficient matrix though).
\begin{theorem}\label{thm:1}
	Under the assumptions {\bf{(A1)-(A4)}}, the sequence $\eta^N$ defined in \eqref{def:eta} converges in distribution to $\eta$ in the space $L^2 ([0, T], H^{-
		\alpha}) \cap C ([0, T], H^{- \alpha - 2})$  for every $\alpha>d/2$, where $\eta$ is the unique  martingale solution to the SPDE \eqref{equation li}.
\end{theorem}
The proof of Theorem \ref{thm:1} will be  given in Section \ref{sec:chara}.

It is  worth emphasizing that the condition $\alpha>d/2$ is  optimal due to the optimal regularity of $\eta$ established in Proposition \ref{pro:op}. Since the driven noise of equation \eqref{equation li} is very rough, so are the solutions. In Section \ref{sec:op}, we rewrite the equation  \eqref{equation li} in the  mild form and study the   regularity of  the stochastic part by Kolmogorov's theorem.  Using the Schauder estimate,   we obtain in Proposition \ref{pro:op}  the optimal regularity  of $\eta$ given by $C([0,T],C^{-\alpha})$  $\mP$-a.s. for every $\alpha>d/2$.

Comparing to the previous result by Fernandez and M\'el\'eard \cite{fernandez1997hilbertian}, Theorem \ref{thm:1} requires less regularity of  the kernel but  more regularity of the solution to the mean-field equation, which eventually would lead to a more restrictive condition on the initial value $\bar{\rho}(0)$. The extra assumption on the mean-field equation  is indeed fairly reasonable, moreover,  by regularity analysis,  the condition $\beta>d/2$  in {\bf{(A4)}} is  optimal on the scale of  H\"older spaces.

 Furthermore, $\eta$ is characterized as a  solution to the following  generalized Ornstein-Uhlenbeck process  in the weak formulation:
\begin{equation}
\left\langle \eta_t,\varphi \right\rangle =\left\langle \eta_{0},Q_{0,t}\varphi\right\rangle +\sqrt{2\sigma}\int_0^t\int_{\mathbb{T}^d}(\nabla Q_{s,t}\varphi)\sqrt{\bar{\rho}_s} \, \xi(\mathd x,\mathd s), \label{eqt:ou}
\end{equation}
for each $\varphi\in C^{\infty}$.  Here the time evolution operators $\{Q_{s,t}\}_{0\leq s\leq t\leq T}$ is defined  for each $t\in [0,T] $ and $\varphi\in C^{\infty}$,
\begin{align}\label{eq:Q}
Q_{\cdot,t}\varphi:=f(\cdot),
\end{align}
with
\begin{enumerate}
	\item $f\in L^2([0,t], H^{\beta+2})\cap C([0,t],H^{\beta+1})$ with  $\p_t f\in L^2([0,t],H^{\beta})$ for $\beta >d/2$.
	
	\item $f$ is the unique solution with terminal value $\varphi$ to the following backward equation
	\begin{equation}
	f_s =\varphi +\sigma \int_s^t\Delta f_r\mathd r+\int_s^t\[K*\bar{\rho}_r\cdot\nabla f_r+K(-\cdot)*(\nabla f_r\bar{\rho}_r) +F\nabla f_r \] \mathd r,\quad s\in [0,t], \nonumber
	\end{equation}
\end{enumerate}
where $K(-\cdot) * g (x) := \int K(y-x) g(y) \mathd y$ and we use this convention throughout the article. For the definition of $\{Q_{s,t}\}$ we refer to Section \ref{sec:gauss} for more details.
The formulation \eqref{eqt:ou} gives  rise to the Gaussianity of the  limit of fluctuation measures. We state the result as follows and we give the proof in Section \ref{sec:gauss}.
\begin{proposition}\label{prop:gauss}
	Under the assumptions {\bf{(A1)-(A4)}}, for the $\eta$ obtained in Theorem \ref{thm:1}, assume in addition that $\bar{\rho}\in C([0,T],C^{\beta+1}(\mathbb{T}^d))$, $F\in C^{\beta+1} (\mathbb{T}^d)$ with $\beta>d/2$, and $\eta_{0}$ in {\bf{(A1)}}  is characterized by
	\begin{equation*}
		\left\langle \eta_0,\varphi\right\rangle \sim \mathcal{N}(0,\left\langle \varphi^2,\bar{\rho}_0 \right\rangle-\left\langle \varphi,\bar{\rho}_0\right\rangle^2  ), \quad \varphi\in C^{\infty}(\mathbb{T}^d).
	\end{equation*}
Then it holds for each test function  $\varphi\in C^\infty$ and $t\in [0,T]$ that
	\begin{equation}
		\left\langle 	\eta_t,\varphi \right\rangle \sim\mathcal{N}\(0,\left\langle |Q_{0,t}\varphi|^2,\bar{\rho}_0\right\rangle-\left\langle Q_{0,t}\varphi,\bar{\rho}_0\right\rangle^2 +2\sigma\int_0^t \left\langle |\nabla Q_{s,t}\varphi|^2,\bar{\rho}_s \right\rangle \mathd s  \).\nonumber
	\end{equation}
\end{proposition}

\medskip
We now focus on the fluctuation problem for the case with vanishing diffusion. In contrast to the non-degenerate case with $\sigma >0$,  due to the vanishing diffusion, the limit equation becomes a deterministic  PDE but  with random initial data. We then analyze the limit equation with the method of characteristics, and obtain the following result in Section \ref{sec:vani}.
\begin{theorem}\label{thm:2}
		Under the assumptions {\bf{(A1)-(A3)}} and  {\bf{(A5)}}, assume further  that $\eta_{0}$ in {\bf{(A1)}}  is characterized by
		\begin{equation*}
			\left\langle \eta_0,\varphi\right\rangle \sim \mathcal{N}(0,\left\langle \varphi^2,\bar{\rho}_0 \right\rangle-\left\langle \varphi,\bar{\rho}_0\right\rangle^2  ), \quad \varphi\in C^{\infty}(\mathbb{T}^d).
		\end{equation*}
	Let $\eta$ be the unique solution to \eqref{equation li} with $\sigma=0$ on the same  stochastic basis with the particle system \eqref{equation pa}. Then the sequence $\eta^N$ defined in \eqref{def:eta} converges in probability to $\eta$ in the space $L^2 ([0, T], H^{-
		\alpha}) \cap C ([0, T], H^{- \alpha - 2})$  for every $\alpha>d/2$. Furthermore, $\eta$ satisfies
	\begin{equation*}
			\left\langle 	\eta_t,\varphi \right\rangle=\left\langle \eta_{0},Q_{0,t}\varphi \right\rangle  \sim\mathcal{N}\(0,\left\langle |Q_{0,t}\varphi|^2,\bar{\rho}_0\right\rangle-\left\langle Q_{0,t}\varphi,\bar{\rho}_0\right\rangle^2\),
	\end{equation*}
 for each test function  $\varphi$ and $t\in [0,T]$. Here the time evolution operator $\{Q_{0,t}\}_{0\leq t\leq T}$ is given by \eqref{eq:Q} with $\sigma=0$.
\end{theorem}

Our main results  validate  that  the relative entropy bound  $\sup_{t \in [0, T]} \sup_N H(\rho_N\vert \bar \rho_N) \lesssim 1 $ which has been established by Jabin and Wang in \cite{jabin2018quantitative} is actually optimal. But the convergence rate  for the marginal distributions $\|\rho_{N, k} - \bar \rho^{\otimes k} \|_{L^\infty( [0, T], L^1)} \lesssim C_T/\sqrt{N} $ is less optimal possibly due to the naive  application of  the CKP inequality as in \eqref{CKP}.  Notice that very recently  Lacker  \cite{lacker2021hierarchies} obtains a sharp estimate for marginal distributions, which is $\|\rho_{N, k} - \bar \rho^{\otimes k} \|_{L^\infty( [0, T], L^1 )} \lesssim C_T/N $ by  local  relative entropy analysis of the BBGKY hierarchy, but under stronger assumptions $H(\rho_{N, k} (0)\vert \bar \rho_0^{\otimes k} ) \lesssim k^2 /N^2$ as well. Even though it is well-known that the convergence of empirical measures and the k-marginal distributions are equivalent in the qualitative sense for instance in \cite{sznitman1991topics}, their quantitative behaviors can be quite complicated when it comes to the order of $N$. See some related discussions in \cite{lacker2021hierarchies,mouhot2013kac,hauray2014kac, mischler2015new}.

\medskip
As an guiding example for our main results in  Theorem \ref{thm:1} and Theorem \ref{thm:2}, we consider the famous vortex model for approximating the  2D Navier-Stokes equation in the vorticity formulation  when $\sigma>0$  and also the 2D Euler equation when $\sigma=0$.  More precisely, given a sequence of i.i.d. initial random variables $\{X_i(0)\}_{i\in \mathbb{N}}$ with a common  probability density function  $\bar{\rho}_0$ on $\mathbb{T}^2$, and consider  the  particle system
\begin{equation}\label{eq:in}
\mathd \nocomma X_i = \frac{1}{N}  \sum_{j \ne i} K (X_i - X_j) \mathd
\nocomma t + \sqrt{2 \sigma} \mathd \nocomma B^i_t, \quad i=1, 2, \cdots, N,
\end{equation}
with  the Biot–Savart law $K:\mathbb{T}^2\rightarrow \mathbb{R}^2$  defined by
\begin{align}\label{eq:K}
K= \nabla^{\perp}G=(-\partial_2G,\partial_1G)
\end{align}
where $G$ is the Green function of the Laplacian on the torus $\mathbb{T}^2$ with mean $0$.  Note in particular that
\begin{align*}
	K(x) = \frac{1}{2\pi} \frac{x^\bot}{|x|^2} + K_0(x),
 \end{align*}
where $x^\bot = (x_1, x_2)^\bot = (-x_2, x_1) \in \mathbb{R}^2$ and $K_0$ is a smooth correction to periodize $K$ on the torus $\mathbb{T}^2$. Obviously the Biot-Savart kernel $K$ satisfies our assumption {\bf (A2)}.

One major corollary of our main results Theorem \ref{thm:1} and Theorem \ref{thm:2} is the following result.
\begin{theorem}\label{th:main}
		If   $\bar{\rho}_0\in C^3(\mathbb{T}^2)$ when $\sigma>0$ and $\bar{\rho}_0\in C^4(\mathbb{T}^2)$  when $\sigma=0$, and $\inf \bar{\rho}_0>0$ for both cases,  then the sequence of  fluctuation measures $\{\eta^N\}_{N\in \mathbb{N}}$ associated with \eqref{eq:in}  converges in distribution to $\eta$ in the space $L^2 ([0, T], H^{-
		\alpha}) \cap C ([0, T], H^{- \alpha - 2})$  for every $\alpha>1$.  Here $\eta$ is a generalized Ornstein-Uhlenbeck process solving  the equation \eqref{equation li} with $K$ given by \eqref{eq:K} and $F=0$. Moreover, $\left\langle \eta,\varphi\right\rangle $ is a centered continuous Gaussian process with  covariance
	\begin{equation}
	\left\langle |Q_{0,t}\varphi|^2,\bar{\rho}_0\right\rangle-\left\langle Q_{0,t}\varphi,\bar{\rho}_0\right\rangle^2 +2\sigma\int_0^t \left\langle |\nabla Q_{s,t}\varphi|^2,\bar{\rho}_s \right\rangle \mathd s  ,\nonumber
	\end{equation}
	where $\{Q_{s,t}\}$ is  introduced in \eqref{eq:Q} with $F=0$ and $\bar{\rho}$ is the solution to the vorticity formulation of 2D imcompressible  Navier-Stokes equation when $\sigma>0$  and 2D Eurler equation when $\sigma=0$.
\end{theorem}

The point vortex approximation towards  2D Navier-Stokes/Euler equation arouses lots of interests since 1980s. The well-posednesss of the point vortex model \eqref{eq:in} was established in \cite{osada1985stochastic,marchioro2012mathematical,takanobu1985existence,fontbona2007paths}.  The main part is  to show that $X_i(t)\neq X_j(t)$ for all $t\in [0,T]$ and $i\neq j$ almost surely, thus the singularity of the kernel will not  be visited almost surely.  
The routine method for instance in  \cite{takanobu1985existence} is based on  estimating the  quantity $\sum_{i\neq j}G(|X_i-X_j|)$, where  $G$ is the Green function. Using the fact $\nabla G\cdot \nabla^{\perp}G=0 $ and by regularization in the intermediate step, it can be  shown that $\sum_{i\neq j}G(|X_i-X_j|)$ is finite almost surely for all  $t\in [0,T]$.   In  \cite{marchioro2012mathematical} by Marchioro and Pulvirenti and  \cite{fontbona2007paths} by Fontbona and Martinez,  the well-posedness of point vortex model with more general circulations/intensities  was established by estimating the displacements of particles. Osada in \cite{osada1985stochastic} obtained the same  result by an analytic approach, which depends on  Gaussian upper and lower bounds for the fundamental solution and the result from \cite{kanda1967regular}.

Osada \cite{osada1986propagation} firstly  obtained a  propogation of chaos result  for \eqref{eq:in} with bounded initial distribution and large viscosity.  More recently, Fournier, Hauray, and Mischler \cite{fournier2014propagation} obtained entropic propagation of chaos by the compactness argument and  their result applies to all viscosity, as long as it is positive, and all initial distributions with finite $k$-moment ($k>0$) and finite Boltzmann entropy.  A quantitative estimate of propagation of chaos has been established  by Jabin and Wang in \cite{jabin2018quantitative} by evolving  the relative entropy between  the joint distribution of \eqref{eq:in} and the tensorized law at the limit. Note in particular  that \cite{jabin2018quantitative} provided  the uniform relative entropy bound as in  {\bf{(A3)}}  for all the kernels  satisfying {\bf{(A2)}}, including the Biot-Savart law.

To the authors' knowledge, Theorem \ref{th:main} is the first result on the fluctuation problems for the  2D Navier-Stokes/Euler equation.

\subsection{Related literatures}\label{sec:i-m}
 Mean field limit and propagation of chaos for the 1st order system given in our canonical form  \eqref{equation pa}  have been extensively studied over the last decade. The basic idea of deriving some effective PDE describing the large scale behaviour of interacting particle systems dates back to Maxwell and Boltzmann. But in our setting, the very first mathematical investigation can be traced back to McKean in \cite{mckeanpropagation}. See also the classical mean field limit from Newton dynamics towards Vlasov Kinetic PDEs in \cite{dobrushin1979vlasov,braun1977vlasov,jabin2015particles,lazarovici2016vlasov} and  the review \cite{jabin2014review}. Recently much progress has been made in the mean field limit for systems as \eqref{equation pa} with singular interaction kernels, including those results focusing on the vortex model \cite{osada1986propagation,fournier2014propagation}, Dyson Brownian motions \cite{berman2019propagation,song2020high,li2020law} and very recently quantitative convergence results  on  general singular kernels  for example as in \cite{jabin2018quantitative, bresch2020mean}  and   \cite{serfaty2020mean,duerinckx2016mean,rosenzweig2020mean,nguyen2021mean}.  See also the references therein for more complete development on the mean field limit.

However, the study of central limit theory for the system \eqref{equation pa}, in particular for those with singular interactions, is quite limited, due to the lack of proper mathematical  tools. 
The  fluctuation problem around a limiting PDE was  popularized for the Boltzmann equation in 1970-1980s for instance in  \cite{mckean1975fluctuations,tanaka1982fluctuation,tanaka1983some,uchiyama1983fluctuations}, but those results  focus more  on the jump-type particle systems.  We also refer to \cite{bodineau2020statistical} for the recent breakthrough on the deviation of the hard sphere dynamics from the kinetic Boltzmann equation. For the fluctuations of interacting diffusions, which is the focus of our article,  to the best of the authors' knowledge,   one of the earliest results  is due to It\^o \cite{ito1983distribution} , where he showed that  for  the system of 1D independent and identically distributed Brownian motions,  the limit of the corresponding fluctuations is a Gaussian process. In the literature, there are mainly  two type results for the fluctuations of interacting processes, either  in the path space or  in the time marginals. This article focuses  on the later one for bounded kernels and some singular ones.  When studying fluctuations in the path space, one treats  processes $\{X_i\}$ as   random variables valued in some functional space,  for instance  the fluctuation measures may be  defined as $\sqrt{N}\(\frac{1}{N}\sum_{i = 1}^N\delta_{X_i}-\mathcal{L}(X)\)$, with the process $X\in C([0,T],\mathbb{R}^d)$ solves the nonlinear stochastic differential equation
\begin{align*}
	X(t)=X(0)+\int_0^t\int_{ \mathbb{R}^d}K(X(s)-x)\mathd \mu_s(x)+\sqrt{2\sigma}B_t,\quad  \mbox{with }\, \mu_s=\mathcal{L}(X(s)).
\end{align*}
To the best of our knowledge, the fluctuation in path space type result was firstly obtained by  Tanaka and Hitsuda \cite{tanaka1981central} and by Tanaka \cite{tanaka1984limit}  for interacting diffusions.  They proved that the fluctuation measures on the path space  converges to a Gaussian random field when  the interacting kernels are bounded and Lipschitz continuous on $\mathbb{R}^1$, with respect to differentiable test functions on the path space.   Later Sznitmann \cite{sznitman1984nonlinear} removed  the differentiability condition  on test functions and generalized the result to $\mathbb{R}^d$, using   Girsanov’s formula and the method of $U$-statistics.

 The article by Fernandez and M\'el\'eard \cite{ fernandez1997hilbertian}  is probably the one closest  to our article when it comes to the basic setting, where they    studied interacting diffusions with regular enough coefficients, using the so-called Hilbertian approach.  Their result cannot cover kernels which are  only bounded or even singular. The systems they consider  are on the whole space and allow multiplicative independent noises. It is worth emphasizing that the  Hilbertian approach introduced in \cite{ fernandez1997hilbertian} has been amplified to study  various interacting models, see \cite{jourdain1998propagation,chevallier2017fluctuations,chen2016fluctuation,lucon2016transition} etc. The Hilbertian approach is based on the  martingale method (as used in this article and many other stochastic problems), coupling method, and analysis in negative weighted Sobolev spaces. The coupling method, which is based on directly  comparing the $N-$ particle system \eqref{equation pa} and $N-$copies of the limit McKean-Vlasov equation,  is also widely used in classical propagation of chaos result \cite{sznitman1991topics},   but usually requires strong assumptions on the interacting kernels and diffusion coefficients. In contrast, our new method enables us to obtain uniform estimates and hence convergence results through directly comparing the Liouville equation and the limit mean-field equation.


We also quickly review some related central limit theory  result for general  interacting particle systems. In  a classical work  \cite{braun1977vlasov} by  Braun and Hepp,  the authors  established  the stability of characteristic flow  in the phase space $\mathbb{R}^{2d}$ with respect  to the initial  measure and thus established the mean field limit for  Newton dynamics with regular interactions towards the Vlasov kinetic PDE. Furthermore, the authors proved that  the limiting behavior of normalized fluctuations around the mean-field characteristics  is a Gaussian process and a precise SDE governing this limit  was also presented. See  its recent generalization by Lancellotti \cite{lancellotti2009fluctuations}.  Budhiraja and Wu \cite{budhiraja2016some} studied some general interacting systems with possible common factors, which do not necessarily  have the exchangability property as usual. Their result is in the flavor of fluctuation in the path space and its  proof follows the strategy by Sznitmann \cite{sznitman1984nonlinear}, i.e. using Girsanov transform and U-statistics.  Furthermore, Kurtz and Xiong \cite{kurtz2004stochastic} studied some interacting particle system from filtering problems. Those SDE's are driven by common noise. The fluctuation result  is similar to the one  driven by independent noises, but the limiting fluctuation is not  Gaussian in general.


 For particle systems in the stationary state with possible singular interaction kernels, there are also many results in the flavor of central limit theory. We only refer to a  few results and readers can find more in reference therein.  Fluctuations for point vortices charged by canonical Gibbs ensembles with the limits given by the so-called energy-enstrophy Gaussian random distributions has been studied by Bodineau and Guionnet  \cite{bodineau1999stationary}  and recently by Grotto and Romito \cite{grotto2020central}. Those results can be regarded as  stationary counterparts of our main theorem in the 2D Euler setting. See also a recent generalization  \cite{geldhauser2021limit} for more singular point vortex model leading to generalized 2D Euler equation but also in the stationary setting.    Moreover,  Lebl\'e and Serfaty \cite{leble2018fluctuations} and Serfaty \cite{serfaty2020gaussian} considered the fluctuation of Coulomb gas on dimension 2 and 3, where the joint distribution of $N$-particle is given by the following Gibbs measure
 \begin{align*}
 	\mathd \mathbb{P}_{N,\beta}=\frac{1}{Z_{N,\beta}}e^{-\frac{\beta}{2}H_N(X_N)},
 \end{align*}
 where $Z_{N, \beta}$ is the partition function, $\beta$ is the temperature, and $H_N$ is the energy including interacting and confining potentials. Now the fluctuation measure, defined as $\sum_{i = 1}^N\delta_{x_i}-N\mu_0$,  where $\mu_0$ is the equilibrium measure, is shown to  converge to a Gaussian free field by using the  Laplace transform and many delicate analysis. In a  similar context as above,  the large deviation principle  for the empirical measure charged by a Gibbs distribution with possible singular Hamiltonian has been obtained by Liu and Wu in  \cite{liu2020large}, even though its dynamical counterpart is still missing and believed to be challenging.  Moreover,  see for instance \cite{perez2007functional}  on  some study  the fluctuations of eigenvalues of random matrices  and a particular case when the eigenvalues are given by Dyson Brownian motions investigated in  Theorem 4.3.20 \cite{anderson2010introduction} where  the fluctuations of moments, when properly normalized, converge to Gaussian processes.


\subsection{Methodology and difficulties}

 The main result (Theorem \ref{thm:1}) follows by the  martingale approach, which has also been used to study the  fluctuation problem of  interacting diffusions with regular kernels as in \cite{meleard1996asymptotic,fernandez1997hilbertian}. The proof consists of three steps:  tightness, identifying the limits of converging subsequences, and well-posedness of the SPDE \eqref{equation li}. By   It\^{o}'s formula, we have
\begin{align}
	\mathd \langle \eta^N_t, \varphi \rangle =& \langle \sigma \Delta \varphi,
	\eta^N_t \rangle \mathd \nocomma t +	\mathcal{K}_t^N (\varphi)  \mathd t + \langle \nabla \varphi, F \eta_t^N \rangle \mathd t
	 + {\frac{\sqrt{2 \sigma_N}}{\sqrt{N}}}
	\sum_{i = 1}^N \nabla \varphi (X_i) \mathd \nocomma B_t^i \nonumber
	\\&+ \sqrt{N}  (\sigma_N - \sigma)  \langle \Delta \varphi,
	\mu_N (t) \rangle\mathd t,\label{equation
		fluctuation}
\end{align}
$\mathbb{P}$-a.s. for each $  \varphi\in C^{\infty}(\mathbb{T}^d)$.
Here the interacting term $	\mathcal{K}_t^N :C^{\infty}(\mathbb{T}^d)\rightarrow \mathbb{R}$ is defined by
\begin{eqnarray}
	\mathcal{K}_t^N (\varphi)  = & \sqrt{N}  \langle \nabla \varphi, K \ast
	\mu_N (t) \mu_N (t) \rangle - \sqrt{N}  \langle \nabla \varphi, \bar{\rho} K
	\ast \bar{\rho} \rangle .  \label{define k}
\end{eqnarray}
To show the  tightness of $\eta^N$, we need to   derive some  uniform estimates for $\eta^N$ in \eqref{equation fluctuation}.  However, due to the singularity of  kernels $K$ in Assumption {\bf{(A2)}}, it seems challenging  to directly obtain  uniform estimates for terms involving $\eta^N$  in the  negative Sobolev spaces. In fact, the optimal regularity  for the limit $\eta$ obtained in Section \ref{sec:op} is in $C_TC^{-\alpha}$ with $\alpha>d/2$. It is natural to consider the energy estimate for $\eta^N$  in $H^{-\alpha}$ using \eqref{equation fluctuation}.  For the purpose of  illustration, let us assume that $\sigma_N = \sigma \equiv 0$, and the exterior force $F=0$ as well, so we can rewrite  \eqref{equation
	fluctuation} as the following  form
\begin{equation*}
\partial_t \eta^N  + \div ( \mu_N K * \eta^N ) + \udiv (\eta^N   K * \bar \rho ) = 0.
\end{equation*}
To control nonlinear terms appearing in the time evolution  $\frac{\ud }{\ud t } \< \eta^N,\eta^N\>_{H^{-\alpha}} $, such as $\<\nabla \eta^N,K*\mu_N\eta^N\>_{H^{-\alpha}}$, we  need  $K\in C^{\beta}$ with $\beta>d/2$  by multiplicative inequality in Appendix \ref{sec:appa}, 
   which is much more demanding than the  assumptions we made on our kernels $K$. 

We overcome this difficulty caused by the singularity of interaction kernels  by using the   Donsker-Varadhan variational formula   \cite[Proposition 4.5.1]{dupuis2011ldp} (see \eqref{eq:var} below) and two large deviation type estimates, one is from \cite[Theorem 4]{jabin2018quantitative} and the other is our contribution (see Lemma \ref{lemma us}). More precisely, now  the uniform estimate of fluctuation measures can be controlled by two terms, one is the relative entropy $H(\rho_N \vert \bar \rho_N)$ and the other is some exponential integrals  with a tensorized  reference measure $\bar{\rho}_N = \bar \rho^{\otimes N}$ (see \eqref{trans}).
 On  one hand,  the uniform bound on  $H(\rho_N \vert \bar{\rho}_N)$,  as summarized in Assumption {\bf{(A3)}},   has already  been established by Jabin and Wang in \cite{jabin2018quantitative}  for a large family of interaction kernels, in particular including those  specified in  {\bf{(A2)}}.   On the other hand, exploiting  cancellation properties from  the interaction terms for instance $\mathcal{K}_t^N$ would enable us to obtain  a uniform bound of  the exponential integrals (see Lemma \ref{lemma jw}  and Lemma \ref{lemma us} for details).   This large deviation type estimate enables the authors of  \cite{jabin2018quantitative} to conclude quantitative estimates of propagation of chaos and also serves a major technical  contribution in our proof.    See Remark \ref{re:uu} for further comments about  the exponential integrals in terms of  $U$-statistics.

Recall the decomposition \eqref{equation fluctuation},  we also  need  to estimate the martingale part  and  show its convergence  as well.  We shall find a pathwise realization $\mathcal{M}^N$ of the martingale part (see Appendix \ref{sec:appb}) and then establish its  tightness.  The tightness of laws of  fluctuation measures then  follows by applying Arzela-Ascoli theorem.

 When characterizing the limit of a converging subsequence, the difficulty still comes from the singularity of  kernels. For the illustration example ($\sigma_N \equiv 0$ and $F=0$), we  notice that it has the following representation
 \begin{equation}
 	\partial _t \eta^N + \udiv (\eta^N K * \bar \rho ) + \udiv ( \bar \rho K * \bar \eta^N ) + \frac{1}{\sqrt{N}} \udiv (\eta^N K * \eta^N) =0,
 \end{equation} More precisely, the  convergence of the interaction term $\mathcal{K}_t^N(\varphi)$  cannot be directly deduced from  convergence of $\mu_N$ and $\eta^N$. We notice that  the interaction term can be splitted into two terms. One term is a  continuous  function of $\eta^N$, or more precisely $ \bar \rho K * \eta^N + \eta^N K * \bar \rho$,  which definitely  converges as $N$ goes  to infinity. The other term is of the form   $\frac{1}{\sqrt{N}} \eta^N K * \eta^N $,  which is not easy to handle directly  since the  formal limit $K*\eta\eta$ is not well-defined (see Lemma \ref{lemma triebel}) in the classical sense due to the singularity of $K$.  Instead, we obtain a uniform bound of  this singular term  by using the variational formula trick again (see Lemma \ref{lemma uni product} below).
The remaining  part for  identifying limits  is  classical.

The last step to Theorem \ref{thm:1} is the uniqueness of martingale solutions to the SPDE \eqref{equation li}, which follows by pathwise uniqueness (see Lemma \ref{lemma unique pde}) and  Yamada-Watanabe theorem.
 Proposition  \ref{prop:gauss} is obtained by solving the dual backward equation of \eqref{equation li} without noises,  which gives the Gaussianity of the limit process of fluctuation measures. 

For the case with vanishing diffusion (which includes the purely deterministic dynamics with $\sigma_N \equiv 0$),  the only difference is on the well-posedness of the limit equation \eqref{equation li}, which is a first order PDE. The well-posedness follows from  the method of characteristics. Since now the  limit equation is deterministic, by a useful lemma in  \cite{gyongy1996existence} by Gy\"ongy and Krylov (see Lemma  \ref{lemma gk} below) we obatin the convergence in probability of the  fluctuation measures.
\subsection{Notations}Throughout the paper, we use the notation $a\lesssim b$ if there exists a universal constant $C > 0$ such that $a\leq Cb$.
 We shall use $\{e_k \}_{k \in
	\mathbb{Z}^d}$ to represent  the Fourier basis on $\mathbb{T}^d$ or $e_k(x) = e^{\sqrt{-1} k \cdot x }$.  For simplicity, we define $\langle k \rangle: = \sqrt{1 + |k|^2}$.

We will mostly  work on  Sobolev spaces, Besov spaces, and the space of $k$-differentiable functions.  The norm of Sobolev space $H^{\alpha}(\mathbb{T}^d)$, $\alpha\in \mathbb{R}$, is defined by
\begin{equation*}
\|f\|_{H^{\alpha}}^2:=\sum_{k\in\mathbb{Z}^d}\langle k \rangle^{2\alpha}|\langle f,e_k \rangle|^2,
\end{equation*}
with the inner product $\<\cdot,\cdot\>_{H^\alpha}$. Moreover, we also use the bracket $\langle \cdot  , \cdot  \rangle$ to denote integrals when the space and underlying measure are clear from the context.   The precise definition and some basic properties  of Besov spaces on torus $B^{\alpha}_{p,q}(\mathbb{T}^d)$ with  $\alpha\in \mathbb{R}$ and $1\leq p,q\leq \infty$, will be given in the Appendix \ref{sec:appa} for completeness.  We remark that  $B^{\alpha}_{2,2}(\mathbb{T}^d)$ coincides with Sobolev space $H^{\alpha}(\mathbb{T}^d)$. We say $f\in C^\alpha(\mT^d)$, $\alpha\in \mathbb{N}$, if $f$ is $\alpha$-times differentiable. For $\alpha \in \mathbb{R}\setminus \mathbb{N} $, the  $C^{\alpha}(\mathbb{T}^d)$ is given by $C^{\alpha}(\mathbb{T}^d)=B^{\alpha}_{\infty,\infty}(\mathbb{T}^d)$. We will often write $\|\cdot\|_{C^{\alpha}}$ instead of $\|\cdot\|_{B^{\alpha}_{\infty,\infty}}$. In the case $\alpha\in \mathbb{R}^+\setminus \mathbb{N}$, $C^{\alpha}(\mathbb{T}^d)$ coincides with the usual H\"older space. We use $C^{\infty}(\mathbb{T}^d)$ to denote  the space of infinitely differentiable functions on $\mathbb{T}^d$, $\cS(\mR^d)$ to denote the class of Schwartz functions on $\mR^d$ and  $\mathcal{S'}(\mathbb{T}^d)$ to denote the space of tempered distributions. Given a Banach space $E$ with a norm $\|\cdot\|_E$ and $T>0$, we write $C_TE=C([0,T];E)$ for the space of continuous functions from $[0,T]$ to $E$, equipped with the supremum norm $\|f\|_{C_TE}=\sup_{t\in[0,T]}\|f(t)\|_{E}$.  For $p\in [1,\infty]$ we write $L^p_TE=L^p([0,T];E)$ for the space of $L^p$-integrable functions from $[0,T]$ to $E$, equipped with the usual $L^p$-norm.

For simplicity, we may omit the underlying space  $\mathbb{T}^d$ without causing confusions.

\subsection{Structure of the paper}
This paper is organized as follows. Section \ref{sec:estimates} is devoted to obtaining  three main estimates which  are  based on the  variational formula and  the large deviation type, including uniform estimates on terms related to $\eta^N$, $\mathcal{K}^N$, and a singular term derived from $\mathcal{K}^N(\varphi)$.  The critical part is to establish some uniform in $N$ estimate of some partition functions. The proof of Theorem \ref{thm:1}  and Proposition \ref{prop:gauss} is completed in Section \ref{sec:SPDE}. First, in Section \ref{sec:tight}, we obtain tightness of the laws of $\{\eta^N\}$ in the space $C([0,T], H^{-\alpha})$ for every $\alpha>d/2+2$, meanwhile we  prove tightness of laws for the pathwise realizations  $\{\mathcal{M}^N\}$ of the martingale part in \eqref{equation fluctuation}.  Then we identify the limits of converging (in distribution)  subsequences of $\{\eta^N\}$ as  a martingale solution to the SPDE \eqref{equation li} and  finish the proof of Theorem \ref{thm:1} in Section \ref{sec:chara}.  The optimal regularity of solutions to the SPDE \eqref{equation li} is shown in Section \ref{sec:op}. Lastly, we prove  Proposition \ref{prop:gauss} in Section \ref{sec:gauss}.

Section \ref{sec:vani} is   concerned with the case with vanishing diffusion, where we give the proof of Theorem \ref{thm:2}. Section \ref{sec:exm}   focuses on some examples which fulfill assumptions {\bf{(A1)-(A5)}}, including the point vortex model approximating   the vorticity formulations of  the 2D Navier-Stokes/Euler equation on the torus.

Finally in Appendix \ref{sec:appa}, we collect the notations and  lemmas about Besov spaces used throughout the paper for completeness. In Appendix \ref{sec:appb} we give the proof of Lemma \ref{lemma reali}, which  shows existence of pathwise realizations.

\subsection*{Acknowledgments}
We would like to thank Pierre-Emmanuel Jabin for helpful discussion.
Z.W. is supported by the National Key R\&D Program of China, Project Number 2021YFA1002800, NSFC grant No.12171009, Young Elite Scientist Sponsorship Program by China Association for Science and Technology (CAST) No. YESS20200028 and the start-up fund from Peking University.
R.Z. is grateful to the financial supports by National Key R\&D Program of China (No. 2022YFA1006300), the NSFC (No.  12271030) and BIT Science and Technology Innovation Program Project 2022CX01001.
 The financial support by the Deutsche Forschungsgemeinschaft (DFG, German Research Foundation) – Project-ID 317210226--SFB 1283 are greatly acknowledged. 
\section{Large Deviation Type Estimates}\label{sec:estimates}
This section collects uniform estimates on $\mu_N-\bar{\rho}$, the interaction term $\mathcal{K}^N$, and a singular term derived from $\mathcal{K}^N(\varphi)$, where $\mathcal{K}^N$ is defined in \eqref{define k}.
These estimates shall play crucial roles in  obtaining tightness and
identifying the limit in Section \ref{sec:SPDE}.
 Indeed, proving the uniform  estimates  is the main
difficulty and technical contribution  of  this article.  Surprisingly, this type estimate, which has been shown to be very useful for many
purposes,   can be actually obtained through a  simple unified  idea.  The quantity we want to bound  can be put in the integral form  $\int \Phi \rho_N$, 
where $\Phi$ is a nonnegative function on $\mathbb{T}^{d \nocomma N}$.  Applying the famous  variational formula from {\cite[Proposition
	4.5.1]{dupuis2011ldp}}, that is
\begin{equation}\label{eq:var}
	\log\int_{\mathbb{T}^{d \nocomma N}} \bar{\rho}_N e^{\Phi} \mathd X^N =
	\sup_{\nu \in \mathcal{P} (\mathbb{T}^{d \nocomma N}), H (\nu | \nobracket
		\bar{\rho}_N) < \infty} \left\{ \int_{\mathbb{T}^{dN}} \Phi \mathd \nu - H
	(\nu | \nobracket \bar{\rho}_N) \right\}, \quad \forall \Phi \geqslant 0,
\end{equation}
with $X^N\assign(x_1,..,x_N)$, $\mathcal{P} (\mathbb{T}^{d \nocomma N})$ the probability measures on $\mathbb{T}^{d \nocomma N}$,
 one can easily control  $\int \Phi \rho_N$ as follows   
\begin{equation}
	\int_{\mathbb{T}^{dN}} \Phi \rho_N \mathd X^N \leqslant \frac{1}{\kappa N}
	\left( H (\rho_N | \nobracket \bar{\rho}_N) + \log \int_{\mathbb{T}^{d
			\nocomma N}} \bar{\rho}_N e^{\kappa N \Phi} \mathd X^N \right),
	\label{trans}
\end{equation}
for any $\kappa > 0$, simply noticing  that $\rho_N$ plays the role of $\nu$ and replacing  $\Phi$ with $\kappa N \Phi$.  See also a direct proof of this inequality  \eqref{trans} in {\cite[Lemma 1]{jabin2018quantitative}} .   As we will see in Lemma \ref{lemma tight mu}, the extra factor  $\frac{1}{N}$ is essential to obtain uniform estimate for fluctuations $\eta^N = \sqrt{N}(\mu_N - \bar \rho )$, but it
comes with a cost that we have to bound the exponential integral $\int \bar \rho_N \exp(\kappa N \Phi )$
uniformly in $N$. Controlling such exponential integrals will be achieved in Section \ref{sec:ldp}, then the uniform estimates will be stated and proved in Section \ref{sec:uni}.
\subsection{Large deviation type estimates}\label{sec:ldp}

As we mentioned before, the major difficulty of our main estimates  is to bound some exponential
integrals, which can be understood as some proper partition functions. To prove  Lemma \ref{lemma uni mu} below, the following  result from Jabin and Wang
{\cite{jabin2018quantitative}} is crucial, and we adapt  it  a bit  below for
convenience.

\begin{lemma}[Jabin and Wang {\cite[Theorem 4]{jabin2018quantitative}}]
	\label{lemma jw}For any
	probability measure $\bar{\rho}$ on $\mathbb{T}^d$,
	and any $\phi (x, y) \in
	L^{\infty} (\mathbb{T}^{2 d})$ with
	\begin{equation}
		\gamma \assign C \| \phi \|_{L^{\infty}}^2 < 1, \nonumber
	\end{equation}
	where $C$ is a universal constant. Assume that $\phi$ satisfies the
	following cancellations
	\begin{equation}
		\int_{\mathbb{T}^d} \phi (x, y)  \bar{\rho} (x) \mathd \nocomma x = 0
		\quad \forall y, \quad \int_{\mathbb{T}^d} \phi (x, y)  \bar{\rho} (y)
		\mathd \nocomma y = 0 \quad \forall x.\nonumber
	\end{equation}
	Then
	\begin{equation}
	\sup_{N \geq 2 }\, 	\int_{\mathbb{T}^{d \nocomma N}} \bar{\rho}_N \exp \big(N \langle \phi, \mu_N
		\otimes \mu_N \rangle \big) \mathd X^N \leqslant \frac{2}{1 - \gamma} < \infty,\nonumber
	\end{equation}
	where $\mu_N = \frac{1}{N} \sum^N_{i = 1} \delta_{x_i}$, $X^N\assign(x_1,..,x_N)\in \mathbb{T}^{dN}$, and  $\bar \rho_N = \bar \rho^{\otimes N}$.	
	
\end{lemma}

Here we abuse the notation  $\bar{\rho}$ since applications below are for the solution $\bar{\rho}$ to \eqref{equation me}.
In this section we also abuse the notations $\mu_N$ and $X^N$,  but we shall always point out the dependence on time when we mention the empirical measure and vector associated to the particle system \eqref{equation pa}.

\begin{remark}
	\label{remark jw}The  proof of the above lemma in  \cite{jabin2018quantitative} relies on the  observation  that $e^A \leqslant e^A + e^{- A}$
	and it is then sufficient  to control the
	series	
	\begin{align}
		\sum_{k = 0}^{\infty} \frac{1}{(2 k) !} \int_{\mathbb{T}^{d \nocomma N}}
		\bar\rho_N A^{2 k} \mathd X^N . &  \nonumber
	\end{align}
	Then of course under the same assumptions, we also have
	\begin{equation}
		\int_{\mathbb{T}^{d \nocomma N}} \bar{\rho}_N \exp \big(N | \langle \phi,
		\mu_N\otimes \mu_N \rangle | \big) \mathd X^N \leqslant \frac{2}{1 - \gamma} <
		\infty .\nonumber
	\end{equation}
Here adding $|\cdot|$ in the exponential part  will be  convenient for proving Lemma \ref{lemma uni product}.
\end{remark}


We also need the following  novel large deviation type estimate on the uniform in $N$  control of a partition function, and use it to obtain the uniform estimate of the interaction term. The proof below is inspired  by   {\cite[Theorem 4]{jabin2018quantitative}}, using combinatorics techniques and some cancellation properties of functions.
\begin{lemma}
	\label{lemma us}For any probability measure $\bar{\rho}$ on $\mathbb{T}^d$. 
	Assume further that functions $\phi (x, y)
	\in L^{\infty} (\mathbb{T}^{2 d})$ with $\|
	\phi \|_{L^{\infty}}$ is small enough, and that
	\begin{equation} \label{CanRule}
		\int_{\mathbb{T}^{2 d}} \bar{\rho} (x)  \bar{\rho} (y) \phi (x, y)
		\mathd \nocomma x \nocomma \mathd \nocomma y = 0.
	\end{equation}
	Then
	\begin{equation}
		\int_{\mathbb{T}^{dN}} \bar{\rho}_N \exp \left( N
		|\langle \phi, \mu_N \otimes \mu_N \rangle |^2  \right) \mathd X^N \leqslant  1 + \frac{
			\alpha_0}{1 - \alpha_0} + \frac{\beta_0}{1 - \beta_0},
		\nonumber
	\end{equation}
	where
	\begin{equation}
		\alpha_0 \assign e^{9}  \| \phi
		\|^2_{L^{\infty}}  < 1, \quad \beta_0 \assign 4 e
		\| \phi \|^2_{L^{\infty}} < 1.
		\nonumber
	\end{equation}
\end{lemma}
\begin{proof}
	We start with the Taylor expansion:
	\begin{align*}
		& \int_{\mathbb{T}^{dN}} \bar{\rho}_N \exp \left( N | \langle \phi, \mu_N\otimes \mu_N \rangle|^2   \right) \mathd X^N =  \sum_{m = 0}^{\infty} \frac{1}{m !}
		\int_{\mathbb{T}^{dN}} \bar{\rho}_N  \left( N
	 | \langle \phi, \mu_N\otimes \mu_N \rangle|^2 \right)^{m} \mathd X^N .
	\end{align*}	
	For the $m$-th term, we use $\mu_N=\frac1N\sum_{i=1}^N\delta_{x_i}$ to expand it as
	\begin{equation}\label{le:mid1}
		\begin{split}
		 &\frac{1}{ m !}  \int_{\mathbb{T}^{dN}} \bar{\rho}_N  \Big( N
 | \langle \phi, \mu_N\otimes \mu_N \rangle|^2\Big)^{m} \mathd X^N
 = \frac{1}{ m !} N^{ m}
 \int_{\mathbb{T}^{dN}} \bar{\rho}_N \bigg(\frac{1}{N^2}\sum_{i,j=1}^N \phi(x_i,x_j)\bigg)^{2m}
 \mathd X^N
\\		= & \frac{1}{ m !} N^{- 3 m}  \sum^N_{i_1, \ldots, i_{2m}, j_1, \ldots, j_{2 m} = 1}
		\int_{\mathbb{T}^{dN}} \bar{\rho}_N  \prod_{\nu = 1}^{2 m}  \phi (x_{i_{\nu }}, x_{j_{\nu }}) \mathd X^N .
		\end{split}
	\end{equation}
	We shall divide the rest proof into two different cases: $4 m > N$ and $4 m
	\leqslant N$.
	
	In the case of $4m > N$,   the $m-$th term, i.e. \eqref{le:mid1} can be bounded trivially
	\begin{align*}
		& \frac{1}{ m !} N^{- 3 m}  \sum^N_{i_1, \ldots, i_{2m}, j_1,  \ldots, j_{2 m} = 1}
		\int_{\mathbb{T}^{dN}} \bar{\rho}_N  \prod_{\nu = 1}^{2 m}  \phi (x_{i_{\nu }}, x_{j_{ \nu }}) \mathd X^N \\
		& \leqslant \frac{1}{ m !} N^{- 3 m} N^{4 m}  \| \phi \|^{2m}_{L^{\infty}} \leqslant m^{-
			\frac{1}{2}} 4^{ m} e^{ m}  \| \phi
		\|^{2m}_{L^{\infty}} .
	\end{align*}
	Here we used the following Stirling's formula with $x = m$
	\begin{equation}\label{Stirling}
		x! = c_x  \sqrt{2 \pi x}  \left( \frac{x}{e}\right)^x,
	\end{equation}
	where $1 < c_x < \frac{11}{10}$ and $c_x \to 1$ as $x \rightarrow \infty$.
	
	In the case of $4 \leqslant 4 m \leqslant N$, we estimate the $m$-th term via
	counting how many choices of multi-indices $(i_1, \ldots, i_{2m}, j_1, \ldots, j_{2 m})$ that
	lead to a non-vanishing integral. If there exists a couple $(i_{q},
	j_{q})$ such that	
	\begin{align}
		i_{q } \neq j_{q}  \, \,  \mbox{and} \,\,\,   i_{q}, j_q  \notin \{  i_{\nu}, j_{\nu} \}  \textrm{ for any } \nu \ne q ,
		\nonumber
	\end{align}
	then variables $x_{i_{q}}$ and $x_{j_{q}}$ enter exactly once in the
	integration. For simplicity, let $(x_{i \nocomma_{q }}, x_{j_{q}}) =
	(x_1, x_2)$, then by Fubini and the cancellation rule \eqref{CanRule} of $\phi$, 	
	\begin{align*}
		& 	\int_{\mathbb{T}^{dN}} \bar{\rho}_N  \prod_{\nu = 1}^{2 m}  \phi (x_{i_{\nu }}, x_{j_{\nu }}) \mathd X^N \\
		&= \int_{\mathbb{T}^{d \nocomma (N - 2)}}
		\left( \int_{\mathbb{T}^{2 d}} \bar{\rho}(x_1) \bar{\rho}
		(x_2) \phi (x_1, x_2) \mathd \nocomma x_1 \mathd \nocomma x_2 \right) \cdot \left(\prod_{\nu = 2}^{2 m} \phi
		(x_{i_{ \nu}}, x_{j_\nu }) \right) \left(\prod_{i \neq 1, 2}  \bar{\rho} (x_i)\right)\mathd
		\nocomma x_3 \ldots \mathd \nocomma x_N
		\\& = 0.
	\end{align*}

	In this case, we introduce auxiliary notations:
	
	$\bullet$ $l$ denotes the number of $x_{i_{\nu}}$ or $x_{j_{\nu}}$ which appears exactly once in
	the integral.
	
	$\bullet$ $p$ denotes the number of $x_{i_{\nu}}$ or  $x_{j_{\nu}}$ which appears at least twice
	in the integral. $\\
	$ A crucial observation is that for multi-indices $(i_1, \ldots, i_{2m}, j_1, \ldots,  j_{2 m})$ which lead to a non-vanishing
	integral, these $l$ variables enter in different couples. This gives $0 \leq l \leq 2m$.  We
	summarize the  following relations among $\{l, p, m, N\}$ as
	\begin{equation*}
		4 		\leqslant 4 m \leqslant N;  \quad  0 \leqslant l \leqslant 2 m;  \quad 1 \leqslant p
		\leqslant (4 m - l) / 2
	\end{equation*}
	
	For a fixed $(l, p)$, notice that there are $\binom{N}{l} \binom{N - l}{p}$
	choices of variables. For each choice of variables, there exists
	$\binom{2 m}{l} 2^l$ choices of place to arrange the $l$ unique variables.
	Lastly, for each arrangement, there are  at most $l!p^{4 m - l}$ plans where $l$! is
	for the $l$ unique variables while $p^{4 m - l}$ is for the other $p$
	variables.
	
	In conclusion, we have when $4 \leq 4m \leq N$,
	
	\begin{align}
		& \frac{1}{ m !} N^{- 3 m}  \sum^N_{i_1, \ldots, i_{2 m}, j_1, \ldots, j_{2m} = 1}
		\int_{\mathbb{T}^{dN}} \bar{\rho}_N  \prod_{\nu = 1}^{2 m}  \phi (x_{i_{\nu}}, x_{j_{\nu }}) \mathd X^N \nonumber\\
		& \leqslant \frac{1}{ m !} N^{- 3 m}
		\| \phi \|^{2m}_{L^{\infty}}  \sum^{2 m}_{l = 0} \sum^{2 m -
			l / 2}_{p = 1} \binom{N}{l} \binom{N - l}{p} \binom{2 m}{l} 2^l\,  l!\, p^{4 m -
			l} \nonumber\\
		& =  \| \phi \|^{2m}_{L^{\infty}}
  \sum^{2 m}_{l = 0} \sum^{2 m - l / 2}_{p = 1} \frac{N!N^{-
				3 m}}{(N - p - l) !}  \frac{\binom{2 m}{l} \, 2^l \,  p^{4 m - l}}{ m !p!} .
		\label{ooo}
	\end{align}
	Applying Stirling's formula \eqref{Stirling} with $x=m, p$ gives
	\begin{equation}
		\sum^{2 m}_{l = 0} \sum^{2 m - l / 2}_{p = 1} \frac{N!N^{- 3 m}}{(N - p -
			l) !}  \frac{\binom{2 m}{l} 2^l p^{4 m - l}}{ m !p!} \leqslant \sum^{2
			m}_{l = 0} \sum^{2 m - l / 2}_{p = 1} N^{p + l - 3 m} 2^l e^{p +  m}
		\frac{\binom{2 m}{l} p^{4 m - l - p}}{ m^{ m}} . \label{ldp 2}
	\end{equation}
	Furthermore, observe that
	\begin{align}
		\frac{\binom{2 m}{l} p^m}{m^m} \leqslant \frac{2^{2 m}  \left(
			2 m - \frac{l}{2} \right)^{ m}}{ m^{ m}} \leqslant 2^{3 m}, \nonumber
	\end{align}
		where we used $2^{2 m} = \sum^{2 m}_{l = 0} \binom{2 m}{l}$ and $p\leq 2m-\frac{l}{2}$. Taking this
	estimate into {\eqref{ldp 2}} yields
	\begin{align}
		\sum^{2 m}_{l = 0} \sum^{2 m - l / 2}_{p = 1} N^{p + l - 3 m} 2^l e^{p +		m}  \frac{\binom{2 m}{l} p^{4m - l - p}}{ m^{ m}} \leqslant \sum^{
		2	m}_{l = 0} \sum^{2m - l / 2}_{p = 1} N^{p + l - 3 m} 2^{l + 3 m} e^{p + 	m} p^{3m - p - l} . \nonumber
	\end{align}
	Since $p + l - 3 m \leqslant 0$ and $p< N$, the above inequality is bounded by $e^{9
		m} $.
		Combining this with {\eqref{ooo}}, we find for every $m \in \left[
	1, \left\lfloor \frac{N}{4} \right\rfloor \right]$,
	
	\begin{align*}
		& \frac{1}{ m!} N^{- 3 m}  \sum^N_{i_1, \ldots, i_{4 m} = 1}
		\int_{\mathbb{T}^{dN}} \bar{\rho}_N  \prod_{\nu = 1}^{2 m}  \phi (x_{i_{ \nu }}, x_{j_{\nu }})\mathd X^N \leqslant \| \phi \|^{2m}_{L^{\infty}}
	e^{9 m}  .
	\end{align*}
	Combining the two cases $4 \leqslant 4 m \leqslant N$ and $4m > N$, it
	follows that
	\begin{align*}
		& \int_{\mathbb{T}^{dN}} \bar{\rho}_N \exp \left( N  |\langle \phi,\mu_N \otimes \mu_N \rangle|^2  \right) \mathd X^N \leqslant   1 + \sum_{m = 1}^{\left\lfloor \frac{N}{4}
			\right\rfloor} \| \phi \|^{2m}_{L^{\infty}}
 e^{9 m} + \sum^{\infty}_{\left\lfloor \frac{N}{4}
			\right\rfloor + 1} m^{- \frac{1}{2}} 4^{ m} e^{ m} \| \phi \|^{2m}_{L^{\infty}}.
	\end{align*}
	Recall that 	
	\begin{align}
		\alpha_0 = e^{9}  \| \phi
		\|^2_{L^{\infty}}  < 1, \quad \beta_0= 4 e
 \| \phi \|^2_{L^{\infty}} < 1.
		\nonumber
	\end{align}
	The proof is thus completed by noticing that
	\begin{align*}
		\sum_{m = 1}^{\left\lfloor \frac{N}{4} \right\rfloor}  \| \phi \|^{2m}_{L^{\infty}}  e^{9m}
		\leqslant  \sum_{m = 1}^{\infty} \alpha_0^{m }  = \frac{ \alpha_0}{1 - \alpha_0}
	\end{align*}
	and
	\begin{align}
		\sum^{\infty}_{\left\lfloor \frac{N}{4} \right\rfloor + 1} m^{-
			\frac{1}{2}} 4^{ m} e^{m}  \| \phi
		\|^{2m}_{L^{\infty}}\leqslant \sum_{m = 1}^{\infty} \beta_0^m
		= \frac{\beta_0}{1 - \beta_0} . \nonumber
	\end{align}
\end{proof}
\begin{remark} Lemma \ref{lemma us} and Lemma \ref{lemma jw} can be generalized in several aspects.  Firstly, the  space $\mathbb{T}^d$ could be replaced by any measurable spaces. Also, when $\phi$ is vector-valued, the result still holds with a slight  modification in the proof as follows
	\begin{align*}
		\prod_{\nu = 1}^{2 m}  \phi (x_{i_{\nu }}, x_{j_{ \nu }}) \xrightarrow{\text{replaced by }} \prod_{\nu = 1}^{ m}  \phi (x_{i_{\nu }}, x_{j_{ \nu }}) \cdot\phi (x_{k_{\nu }}, x_{l_{ \nu }}).
	\end{align*}
Indeed,  given $\phi$ a vector-valued function, the modification only comes from the expanding
\begin{equation*}
|\langle \phi, \mu_N \otimes \mu_N \rangle |^2 = \Big| \frac{1}{N^2}  \sum_{i, j=1}^N \phi (x_i , x_j)\Big|^2 = \frac{1}{N^4} \sum_{i, j, k, l =1}^N  \phi(x_i, x_j) \cdot \phi(x_k, x_l) .
\end{equation*}

 \end{remark}
\begin{remark}\label{re:uu}
	In Lemma \ref{lemma jw} and Lemma \ref{lemma us}, we proved that  two exponential integrals which are in the form of  $\mathbb{E}e^{NU_N}$ are uniformly bounded with respect to  $N$. In the first case as in \cite{jabin2018quantitative},
	\begin{equation*}
		U_N^1 = \frac{1}{N^2}\sum_{i ,j= 1}^N\phi(X_i,X_j),
	\end{equation*}
	while in the second case, $U_N^2 =  (U_N^1)^2$, or  $U_N^2$ may be   expressed as
	\begin{align*}
		U_N^2= \frac{1}{N^4}\sum_{i_1,i_2,i_3,i_4=1}^{N}\phi(X_{i_1},X_{i_2})\phi(X_{i_3},X_{i_4}).
	\end{align*}
	 Those $U_N$ in both cases are $U$-statistics, which are symmetric functions of   $N$ i.i.d  random variables. The degree of $U_N$ is said to be $k$ if $U_N$ is a symmetric version of an  interaction  function between $k$ variables. So $U_N$ in Lemma \ref{lemma jw}  is of degree $2$ and $U_N$ in Lemma \ref{lemma us} is of degree $4$.  Here the cancellation properties actually imply the first order degeneracy of $U$-statistics, which  together with the boundedness condition gives the weak convergence of the law of $NU^i_N$, $i=1,2$. We refer to \cite{lee2019u} for more details.
\end{remark}
\subsection{Uniform estimates}\label{sec:uni}Now we are in the position to state and prove the uniform estimates.

The first estimate concerns on the convergence from $\mu_N (t)$ to
$\bar{\rho}_t$ in $H^{- \alpha}$, for $\alpha > d / 2$.
\begin{lemma}
	\label{lemma uni mu}For each $\alpha > d / 2$, there exists a constant
	$C_{\alpha}$ such that, for all $t \in [0, T]$,
	\begin{equation*}
		\mathbb{E} \| \mu_N (t) - \bar{\rho}_t \|^2_{H^{- \alpha}} \leqslant
		\frac{C_{\alpha}}{N} (\nocomma H_t (\rho_N | \bar{\rho}_N) + 1),
	\end{equation*}
where we recall $\mu_N(t) = \frac{1}{N} \sum_{i=1}^N \delta_{X_i(t)}$ and the expectation is taken according to the joint distribution $\rho_N(t, \cdot)$ of the particle system \eqref{equation pa} and $\bar \rho_N(t, \cdot) = \bar \rho(t)^{\otimes N}$.
\end{lemma}

This lemma has a direct consequence.  Recall that the fluctuation measure $\eta^N(t) = \sqrt{N} (\mu_N(t) - \bar \rho_t )$. Under Assumption {\bf (A3)}, i.e. $\sup_{t \in [0, T]} \sup_{N } H_t(\rho_N\vert \bar \rho_N ) \lesssim 1 $, one can then immediately obtain
\begin{equation*}
	\sup_{t \in [0, T]} \sup_{N }   	\mathbb{E} \|  \eta^N (t) \|^2_{H^{- \alpha}} \lesssim 1, \quad \mbox{for } \, \alpha > d/2.
	\end{equation*}

\begin{proof}
	Since the Dirac measure belongs to $H^{- \alpha} (\mathbb{T}^d)$
	for every $\alpha > d / 2$, it follows that $\mu_N (t) - \bar{\rho}_t \in
	H^{- \alpha} (\mathbb{T}^d)$.
		Then by {\eqref{trans}}
	, we find for any $\kappa > 0$,
	\begin{equation}
		\begin{split}
		\mathbb{E} \| \mu_N (t) - \bar{\rho}_t \|^2_{H^{- \alpha}}  &  = \int_{\mathbb{T}^{dN }} \|\mu_N - \bar \rho_t  \|_{H^{- \alpha }}^2 \rho_N(t, X^N )  \mathd  X^N   \\
	 & 	\leqslant \frac{1}{\kappa N}  \left( \nocomma H_t (\rho_N | \bar{\rho}_N) +
		\log \int_{\mathbb{T}^{dN}} \exp\({\kappa N\|\mu_N-\bar{\rho}_t\|_{H^{-\alpha}}^2} \) \bar{\rho}_N \mathd X^N
		\right) . \label{pp}
		\end{split}
	\end{equation}
	Recalling $\{e_k \}_{k \in
		\mathbb{Z}^d}$ is the Fourier basis, and
	\begin{align}
		 \| \mu_N  - \bar{\rho}_t \|^2_{H^{- \alpha}} =  &
		\sum_{k \in \mathbb{Z}^d} \langle k \rangle^{- 2 \alpha} |\langle e_k,
		\mu_N  - \bar{\rho}_t \rangle  |^2\nonumber.
	\end{align}
Since the exponential function is convex, using Jensen's inequality gives that
\begin{align}
\int_{\mathbb{T}^{dN}} \exp\({\kappa N\|\mu_N-\bar{\rho}_t\|_{H^{-\alpha}}^2} \) \bar{\rho}_N \mathd X^N
	\leq& \frac{1}{C}\sum_{k\in \mathbb{Z}^d}\langle k \rangle^{- 2 \alpha}	\int_{\mathbb{T}^{dN}}\exp \(\kappa NC|\langle e_k,
		\mu_N - \bar{\rho}_t \rangle  |^2\) \bar{\rho}_N \mathd X^N\nonumber
		\\\leq &\sup_{k\in \mathbb{Z}^d}\int_{\mathbb{T}^{dN}}\exp \(\kappa NC|\langle e_k,
		\mu_N  - \bar{\rho}_t \rangle  |^2\) \bar{\rho}_N \mathd X^N,\label{pppp}
\end{align}
	where the constant $C=\sum_{k\in \mathbb{Z}^d}\langle k \rangle^{- 2 \alpha}$ depends only on $\alpha$ and is finite  since $\alpha > d / 2$.

	We define
	\begin{align}
		\phi_1 (t, k,x, y) \assign &  [e_k (x) - \langle e_k, \bar{\rho}_t \rangle]  [e_{- k} (y) -
		\langle e_{- k}, \bar{\rho}_t \rangle], \nonumber
	\end{align}
therefore
\begin{align*}
	\int_{\mathbb{T}^{dN}}\exp \(\kappa NC|\langle e_k,
	\mu_N - \bar{\rho}_t \rangle  |^2\) \bar{\rho}_N \mathd X^N=\int_{\mathbb{T}^{dN}}\exp \(\kappa NC\left\langle \phi_1(t,k,\cdot,\cdot),\mu_N\otimes\mu_N\right\rangle \) \bar{\rho}_N \mathd X^N.
\end{align*}
Since $\bar{\rho}$ is a probability measure,  $\|\phi_1\|_{L^{\infty}}$ is bounded uniformly in $t$ and $k$. One can also easily check that
	\begin{align}
		\int_{\mathbb{T}^d} \phi_1 (t, k, x, y)  \bar{\rho}_t (x) \mathd
		\nocomma x = 0 \quad \forall y, &\quad  \int_{\mathbb{T}^d}  \phi_1 (t, k,
		x, y)  \bar{\rho}_t (y) \mathd \nocomma y = 0 \quad \forall x. \nonumber
	\end{align}
	Then by Lemma \ref{lemma jw} with $\kappa$ (depending on $\alpha$) small enough, we deduce that
	\begin{align}
		&\sup_N\sup_{k\in \mathbb{Z}^d}\int_{\mathbb{T}^{dN}}\exp \(\kappa NC|\langle e_k,
		\mu_N  - \bar{\rho}_t \rangle  |^2\) \bar{\rho}_N \mathd X^N \nonumber
		\\& =\sup_N\sup_{k\in \mathbb{Z}^d}
			\int_{\mathbb{T}^{dN}} \bar{\rho}_N \exp (\kappa NC \langle \phi_1 (t, k,
			\cdummy, \cdummy), \mu_N \otimes \mu_N \rangle) \mathd X^N\nonumber\\
			& =\sup_N\sup_{k\in \mathbb{Z}^d} \int_{\mathbb{T}^{dN}} \bar{\rho}_N \exp (\kappa N C\langle \tmop{Re}
			\phi_1 (t, k,\cdummy, \cdummy), \mu_N \otimes \mu_N \rangle) \mathd X^N <\infty,\label{pp0p}
	\end{align}
where the equalities follows by
\begin{align*}
	|\langle e_k,
	\mu_N  - \bar{\rho}_t \rangle  |^2=\langle \phi_1 (t, k,
	\cdummy, \cdummy), \mu_N \otimes \mu_N \rangle\in \mathbb{R}.
\end{align*}

Combining  {\eqref{pp}}-\eqref{pp0p} yields	
	\begin{align}
		\mathbb{E} \| \mu_N (t) - \bar{\rho}_t \|^2_{H^{- \alpha}} \leqslant &
		\frac{1}{\kappa  N}  (\nocomma H_t (\rho_N | \bar{\rho}_N) +
		C_{\alpha}), \nonumber
	\end{align}	
	where $C_{\alpha}$ is a constant depending only on $\alpha$. We thus arrive
	at the result.
\end{proof}
In particular, Lemma \ref{lemma uni mu} gives  the tightness of laws of $\{\eta^N(0)\}$ on $H^{-\alpha}$ under the condition that $H(\rho_{ N}(0)|\bar{\rho}_N(0))$ is finite, which together with Assumption {\bf{(A1)}} yields the convergence of $\{\eta^N(0)\}$ in the negative Sobolev spaces.
\begin{corollary}\label{coro22}For every $\alpha>d/2$,
	$\eta^N_0$ converges in distribution to $\eta_0$ given by {\bf{(A1)}} in $H^{-\alpha}$.
\end{corollary}

The next lemma concerns on the interaction part in the decomposition {\eqref{equation fluctuation}}.
\begin{lemma}
	\label{lemma uni nonlinear}If the  kernel $K $ satisfies Assumption {\bf{(A2)}},
then for each $\alpha >
	d / 2 + 2$, there exists a constant $C_{\alpha}$ such that, for all $t \in
	[0, T]$,
	\begin{equation}
		\mathbb{E} \| \nabla \cdummy [K \ast \mu_N (t) \mu_N (t) - \bar{\rho}_t K
		\ast \bar{\rho}_t] \|^2_{H^{- \alpha}} \leqslant \frac{C_{\alpha}}{N}
		(\nocomma H_t (\rho_N | \bar{\rho}_N) + 1) ,\nonumber
	\end{equation}
where the expectation is taken according to the joint distribution $\rho_N(t, \cdot)$ of the particle system \eqref{equation pa}.
\end{lemma}

\begin{proof}
	The proof is similar to Lemma \ref{lemma uni mu}. First, by {\eqref{trans}} we find for any $\kappa > 0$,
	
	\begin{align}
	&	\mathbb{E} \| \nabla \cdummy [K \ast \mu_N (t) \mu_N (t) - \bar{\rho}_t K
		\ast \bar{\rho}_t] \|^2_{H^{- \alpha}} \nonumber
		\\&\leqslant
		\frac{1}{\kappa N}  \left( \nocomma H_t (\rho_N | \bar{\rho}_N) + \log
		\int_{\mathbb{T}^{dN}} \bar{\rho}_N\exp \left( \kappa N  \| \nabla \cdummy [K \ast \mu_N  \mu_N  - \bar{\rho}_t K
		\ast \bar{\rho}_t] \|^2_{H^{- \alpha}}\right)\mathd X^N \right)
		. \label{ppp}
	\end{align}
Next, we find that
\begin{align}
	\| \nabla \cdummy [K \ast \mu_N  \mu_N  - \bar{\rho}_t K
	\ast \bar{\rho}_t] \|^2_{H^{- \alpha}}=& \sum_{k \in
		\mathbb{Z}^d} \langle k \rangle^{- 2 \alpha} | \langle \nabla e_k, K \ast
	\mu_N  \mu_N  - \bar{\rho}_t K \ast \bar{\rho}_t \rangle |^2\nonumber
	\\ \leq &  \sum_{k \in \mathbb{Z}^d\setminus \{ 0\} } \langle k
	\rangle^{- 2 \alpha} |k|^2 | \langle e_k, K \ast \mu_N  \mu_N-
	\bar{\rho}_t K \ast \bar{\rho}_t \rangle |^2.\nonumber
\end{align}
	For the case $|x| \nocomma K (x) \in L^{\infty}$ and $K (x) = - K (- x)$, we do a symmetrization trick. That is,  for any $\varphi \in C^\infty (\mathbb{T}^d)$ and a probability measure $\mu$,
	\begin{equation*}
		\begin{split}
		\int_{\mathbb{T}^d} \varphi(x) K * \mu(x) \mu(\mathd x )  & = \int_{\mathbb{T}^{2d}} \varphi(x) K(x-y) \mu^{\otimes 2 }(\mathd x \mathd y ) \\ & = \frac 1 2 \int_{\mathbb{T}^{2d}} (\varphi (x) - \varphi (y)) \cdot K (x-y) \mu^{\otimes 2 } (\mathd x \mathd y ).
		\end{split}
	\end{equation*}
We define that \begin{align*}
	\mathbb{K}_{\varphi} (x, y) : = \frac{1}{2}K (x - y) [\varphi (x) - \varphi
	(y)], \quad \forall\varphi \in C^{\infty} (\mathbb{T}^d).
\end{align*}
Thus in this case, $\|\mathbb{K}_{\varphi }\|_{L^\infty} \lesssim  \|\nabla \varphi \|_{L^\infty} \| |x| K \|_{L^\infty}$.  Consequently, since
\begin{equation*}
	\langle e_k, K \ast \mu_N  \mu_N-
	\bar{\rho}_t K \ast \bar{\rho}_t \rangle  =  	\langle  \mathbb{K}_{e_k}(\cdot,\cdot), \mu_N^{\otimes 2 } - \bar \rho_t^{\otimes 2} 	\rangle,
	\end{equation*}
and  $\| \mathbb{K}_{e_k}\| \lesssim |k|$, one proceeds as 	
	\begin{align}
		\| \nabla \cdummy [K \ast \mu_N  \mu_N  - \bar{\rho}_t K
	\ast \bar{\rho}_t] \|^2_{H^{- \alpha}} = &  \sum_{k \in \mathbb{Z}^d\setminus\{0\}} \langle k \rangle^{-
			2 \alpha} |k|^2 | \langle \mathbb{K}_{e_k},  \mu_N \otimes  \mu_N -
		\bar{\rho}_t \otimes \bar{\rho}_t \rangle |^2\nonumber
		\\=&\sum_{k \in \mathbb{Z}^d\setminus \{ 0\}} \langle k \rangle^{-
			2 \alpha} |k|^4 | \langle\frac{\mathbb{K}_{e_k}}{|k|} ,  \mu_N \otimes  \mu_N -
		\bar{\rho}_t \otimes \bar{\rho}_t \rangle |^2, \nonumber
	\end{align}
	where now $\frac{\mathbb{K}_{e_k}}{|k|}$ is bounded.
	
\noindent	For the case that $K \in L^\infty$ but $K$ is not necessarily anti-symmetric, we directly write
	\begin{equation*}
		\langle e_k, K \ast \mu_N  \mu_N-
		\bar{\rho}_t K \ast \bar{\rho}_t \rangle  =  	\langle  e_k(x) K(x-y), \mu_N^{\otimes 2 } - \bar \rho_t^{\otimes 2} 	\rangle.
		\end{equation*}
	
\noindent To sum it up, we define $\phi_2 : [0, T] \times\{\mathbb{Z}^d\setminus \{0\}\}\times \mathbb{T}^{2 d \nocomma} \rightarrow
	\mathbb{R}^d$  by
	\begin{align}
		\phi_2 (t,k, x, y) : =
		\begin{cases}  K (x - y) e_k (x) -  \langle e_k, \bar{\rho}_t K
			\ast \bar{\rho}_t \rangle,&  \mbox{if } \, K \in L^{\infty} \nonumber\\\frac{\mathbb{K}_{e_k}(x,y)}{|k|} - \langle
			\frac{\mathbb{K}_{e_k}}{|k|}, \bar{\rho}_t \otimes \bar{\rho}_t \rangle, & \mbox{if }\, |x| \nocomma K (x) \in L^{\infty}, K (x) = - K (- x).
		\end{cases}
	\end{align}
Using Jensen's inequality, for both cases we have
\begin{align}
	&\int_{\mathbb{T}^{dN}} \bar{\rho}_N\exp \left( \kappa N  \| \nabla \cdummy [K \ast \mu_N  \mu_N  - \bar{\rho}_t K
	\ast \bar{\rho}_t] \|^2_{H^{- \alpha}}\right)\mathd X^N \nonumber
	\\&=\int_{\mathbb{T}^{dN}} \bar{\rho}_N\exp \bigg( \kappa N  \sum_{k \in \mathbb{Z}^d\setminus \{0\} } \langle k
	\rangle^{- 2 \alpha} |k|^2 | \langle e_k, K \ast \mu_N  \mu_N-
	\bar{\rho}_t K \ast \bar{\rho}_t \rangle |^2\bigg)\mathd X^N \nonumber
	\\&\leq \int_{\mathbb{T}^{dN}} \bar{\rho}_N\exp \bigg( \kappa N  \sum_{k \in \mathbb{Z}^d\setminus \{0\}} \langle k
	\rangle^{- 2 \alpha+4} | \langle \phi_2(t,k,\cdot,\cdot),  \mu_N \otimes \mu_N\rangle |^2\bigg)\mathd X^N \nonumber
	\\&\leq \frac{1}{C}\sum_{k \in \mathbb{Z}^d\setminus \{0\}} \langle k
	\rangle^{- 2 \alpha+4}\int_{\mathbb{T}^{dN}} \bar{\rho}_N\exp \left( \kappa N  C| \langle \phi_2(t,k,\cdot,\cdot),  \mu_N \otimes \mu_N\rangle |^2\right)\mathd X^N \nonumber
	\\&\leq \sup_{k\in \mathbb{Z}^d\setminus \{0\}}\int_{\mathbb{T}^{dN}} \bar{\rho}_N\exp \left( \kappa N  C| \langle \phi_2(t,k,\cdot,\cdot),  \mu_N \otimes \mu_N\rangle |^2\right)\mathd X^N ,\label{ppp1}
\end{align}
	where the constant $C:=\sum_{k\in \mathbb{Z}^d\setminus0}\langle k \rangle^{- 2 \alpha+4}$ depends only on $\alpha$, and is finite since $\alpha > d / 2+2$.

Furthermore, since $\phi_2$ is complex-valued, we find
\begin{align}
&\sup_N	\sup_{k\in \mathbb{Z}^d\setminus \{0\}}\int_{\mathbb{T}^{dN}} \bar{\rho}_N\exp \left( \kappa N  C| \langle \phi_2(t,k,\cdot,\cdot),  \mu_N \otimes \mu_N\rangle |^2\right)\mathd X^N \nonumber
\\&\leq \frac{1}{2}\sup_N	\sup_{k\in \mathbb{Z}^d\setminus \{0\}}\int_{\mathbb{T}^{dN}} \bar{\rho}_N\exp \left( 2\kappa N  C|\langle \tmop{Re}\phi_2(t,k,\cdot,\cdot),  \mu_N \otimes \mu_N\rangle|^2\right)\mathd X^N \nonumber
\\&\quad + \frac{1}{2}\sup_N	\sup_{k\in \mathbb{Z}^d\setminus \{0\}}\int_{\mathbb{T}^{dN}} \bar{\rho}_N\exp \left( 2\kappa N  C|\langle \tmop{Im}\phi_2(t,k,\cdot,\cdot),  \mu_N \otimes \mu_N\rangle|^2  \right)\mathd X^N ,\label{ppp2}
\end{align}
where the inequality follows by Jensen's inequality and the fact that
\begin{align*}
	 |\langle \phi_2(t,k,\cdot,\cdot),  \mu_N \otimes \mu_N\rangle |^2=|\langle \tmop{Re}\phi_2(t,k,\cdot,\cdot),  \mu_N \otimes \mu_N\rangle|^2 +|\langle \tmop{Im}\phi_2(t,k,\cdot,\cdot),  \mu_N \otimes \mu_N\rangle|^2.
\end{align*}
One can easily find that $\|\phi_2\|_{L^{\infty}}$ is bounded uniformly in $(t,k)$, and satisfies the cancellation
	\begin{align}
		\int_{\mathbb{T}^{2 d}} \phi_2 (t,k,x, y)  \bar{\rho}_t (x)  \bar{\rho}_t (y)
		\mathd \nocomma x \nocomma \mathd \nocomma y = 0,&  \nonumber
	\end{align}
and so do the real and imaginary part of $\phi_2$.
	
	Choosing $\kappa $ (depending on $\alpha$) sufficiently small, then we are able to apply
	Lemma \ref{lemma us} to obtain that
	\begin{align}
	\sup_N	\sup_{k\in \mathbb{Z}^d\setminus0}\int_{\mathbb{T}^{dN}} \bar{\rho}_N\exp \left( \kappa N  C| \langle \phi_2(t,k,\cdot,\cdot),  \mu_N \otimes \mu_N\rangle |^2\right)\mathd X^N 	
		\leqslant & C_{\alpha} , \nonumber
	\end{align}
where the universal  constant $C_{\alpha}$  only depends  on $\alpha$. The proof is then  completed by combining this with  {\eqref{ppp}} and \eqref{ppp1}.
\end{proof}

The last estimate in this section  plays a crucial role in identifying the limit in Section \ref{sec:chara}.

\begin{lemma}
	\label{lemma uni product} If the  kernel $K $ satisfies  assumption {\bf{(A2)}}, 
then for each $\varphi \in
	C^1$, there exists a  universal constant $C$ such that, for all $t \in [0, T]$,
	\begin{equation}
		\mathbb{E} | \langle \varphi K \ast (\mu_N (t) - \bar{\rho}_t), \mu_N (t)
		- \bar{\rho}_t \rangle | \leqslant \frac{C}{ N}  (H_t (\rho_N |
		\bar{\rho}_N) + 1) ,\nonumber
	\end{equation}
where the expectation is taken according to the joint distribution $\rho_N(t, \cdot)$ of the particle system \eqref{equation pa}.
\end{lemma}

\begin{proof}
	We first write the quantity in the following form
	\begin{equation}
		\mathbb{E} | \langle \varphi K \ast (\mu_N (t) - \bar{\rho}_t), \mu_N (t)
		- \bar{\rho}_t \rangle | =\mathbb{E} | \Phi (t, X_t^N) | =
		\int_{\mathbb{T}^{dN}} | \Phi (t, X^N) | \rho_N \mathd X^N,
		\label{product 2}
	\end{equation}
	where $\Phi$ is defined by
	
	\begin{align*}
		\Phi  (t, X^N) & = \langle \varphi K \ast ( \mu_N- \bar{\rho}_t),  \mu_N
		- \bar{\rho}_t \rangle .
	\end{align*}
	
	For the  case $K \in L^{\infty}$, we find
	\begin{align}
		\Phi  (t, X^N) & = \langle \phi_3 (t, \cdummy, \cdummy),  \mu_N \otimes
		 \mu_N \rangle, \nonumber
	\end{align}
		with $\phi_3$ defined by
	\begin{align}
		\phi_3 (t, x, y) := & K (x - y) \varphi (x) - \varphi (x) K \ast
		\bar{\rho}_t (x) - \langle K (\cdummy - y) \varphi, \bar{\rho}_t \rangle +
		\langle \varphi K \ast \bar{\rho}_t, \bar{\rho}_t \rangle . \nonumber
	\end{align}
	
		For the case $|x| \nocomma K (x) \in L^{\infty}$ and $K (x) = - K (- x)$, 
	we do a symmetrization for $\Phi$ as in the proof of Lemma \ref{lemma uni nonlinear}, i.e.
	\begin{align}
		\Phi  (t, X^N) & = \langle \mathbb{K}_{\varphi}, ( \mu_N - \bar{\rho}_t)
		\otimes  \mu_N - \bar{\rho}_t \rangle = \langle \phi_3 (t, \cdummy,
		\cdummy),  \mu_N \otimes  \mu_N \rangle, \nonumber
	\end{align}
		with $\phi_3$ defined by
	\begin{align}
		\phi (t, x, y) := & \mathbb{K}_{\varphi} (x, y) -  \langle
		\mathbb{K}_{\varphi} (x, \cdummy), \bar{\rho}_t \rangle -\langle
		\mathbb{K}_{\varphi} ( \cdummy,y), \bar{\rho}_t \rangle+ \langle
		\mathbb{K}_{\varphi}, \bar{\rho}_t \otimes \bar{\rho}_t \rangle .
		\nonumber
	\end{align}
		By {\eqref{trans}}, it holds for any $\kappa > 0$ that
	\begin{equation}
		\int_{\mathbb{T}^{dN}} | \Phi (t, X^N) | \rho_N \mathd X^N
		\leqslant \frac{1}{\kappa N}  \left( \nocomma H (\rho_N | \bar{\rho}_N) +
		\log \int_{\mathbb{T}^{dN}} \bar{\rho}_N e^{\kappa N | \Phi |} \mathd
		X^N \right) . \label{product 1}
	\end{equation}
	On the other hand, one can easily check the following cancellations
	
	\begin{align}
		\int_{\mathbb{T}^d} \phi_3 (t, x, y)  \bar{\rho}_t (x) \mathd \nocomma x =
		0 \quad \forall y, \quad \int_{\mathbb{T}^d} \phi_3 (t, x, y)
		\bar{\rho}_t (y) \mathd \nocomma y = 0 \quad \forall x. \nonumber
	\end{align}
	Since in both cases, $\phi_3$ is bounded uniformly in $t$, we can choose
	$\kappa$ such that $\sqrt{\kappa} \| \phi_3 \|_{L^{\infty}}$ sufficiently
	small. Letting $\bar{\rho}_t$ and $\sqrt{\kappa} \phi_3$ play the roles of
	$\bar{\rho}$ and $\phi$ in Remark \ref{remark jw} respectively, we deduce
	that
	\begin{align}
		\int_{\mathbb{T}^{dN}} \bar{\rho}_N e^{\kappa N | \Phi |} \mathd X^N
		\leqslant & C, \nonumber
	\end{align}
	where $C$ is a constant depending only on $\varphi$. Combining this with
	{\eqref{product 2}} and {\eqref{product 1}}, we thus arrive at the result.
\end{proof}
\section{The SPDE Limit}\label{sec:SPDE}

The aim of this section is to analyze fluctuation behavior of the empirical measure $\mu^N$ for the non-degenerate case, i.e. $\sigma >  0$.
It will be shown that $\eta^N = \sqrt{N} (\mu_N  - \bar{\rho})$
converges in distribution to the unique solution $\eta$ to the  linear SPDE
{\eqref{equation li}}. We shall start with proving that the sequence of   $(\eta^N
)_{N \geqslant 1}$ is tight. 
 Then each tight limit of the subsequence from
$(\eta^N)_{N \geqslant 1}$ will be identified as a
martingale solution to the  equation {\eqref{equation li}}. The next
step is to show  pathwise uniqueness of {\eqref{equation li}}, which
allows us to conclude the proof of Theorem \ref{thm:1}.  In  Section \ref{sec:op}, we prove the optimal regularity of solutions to the  limit SPDE \eqref{equation li}. Finally, the proof of Proposition \ref{prop:gauss}, which gives the Gaussianity of the unique limit of fluctuation measures, is given in Section \ref{sec:gauss}.

\subsection{Tightness }\label{sec:tight}

Before proving tightness, we introduce pathwise realization of the martingale
part in the  decomposition {\eqref{equation fluctuation}}. Recall that the martingale part is given by
\begin{equation}
	\frac{\sqrt{2 \sigma_N}}{\sqrt{N}} \sum_{i = 1}^N \int^t_0 \nabla \varphi
	(X_i) \mathd \nocomma B_s^i,\nonumber
\end{equation}
for each $\varphi \in C^{\infty} (\mathbb{T}^d)$. Formally, one could define a
random operator $\mathcal{M}_t^N : \Omega \times H^{\alpha} \rightarrow \mathbb{R}$, $\alpha>d/2+1$, 
 for each $t \in [0, T]$ through
\begin{equation}
	\mathcal{M}_t^N (\varphi) = \frac{\sqrt{2 \sigma_N}}{\sqrt{N}}  \sum_{i =
		1}^N \int^t_0 \nabla \varphi (X_i) \mathd \nocomma B_s^i, \quad \mathbb{P}-
	a.s. \label{equation realization}
\end{equation}
However, the measurability of $\mathcal{M}_t^N : \Omega \rightarrow H^{-
	\alpha}$ is nontrivial due to the fact that the above stochastic integral is
defined as a $\mathbb{P}$-equivalence class for each $\varphi$. Finding a
measurable map $\mathcal{M}_t^N$ from $\Omega$ to $H^{-\alpha}$ requires a pathwise meaning of the map
$\mathcal{M}_t^N (\varphi)$, such that $\mathcal{M}_t^N (\varphi)$ is
continuous with respect to $\varphi$ for almost every $\omega \in \Omega$.

Pathwise realization has been studied for different function spaces, for instance  in \cite[Theorem 3.1]{ito1983distribution}, {\cite{flandoli2005stochastic}},  \cite[Lemma 9]{mourrat2017global},  etc.
Adapting the
idea of investigating stochastic currents in {\cite{flandoli2005stochastic}} by  Flandoli, Gubinelli, Giaquinta, and Tortorelli,
a pathwise realization $\mathcal{M}^N$ with values  in Hilbert spaces can be obtained with a relatively simple
proof, which is postponed into Appendix \ref{sec:appb}. We state the result as follows. 

\begin{lemma}\label{lemma reali}
	For each $N$, there exists a progressively measurable process
	$\mathcal{M}^N$ with values in $H^{- \alpha}$, for any $\alpha > d / 2 +
	1$, such that {\eqref{equation realization}} holds almost surely for all $t
	\in [0, T]$ and $\varphi \in C^{\infty}$.
\end{lemma}

In the following, we are going to prove tightness of $(\eta^N,
\mathcal{M}^N)_{N \geqslant 1}$.
To start, recall the following tightness criterion
given by Arzela-Ascoli theorem {\cite[Theorem
	7.17]{kelley2017topology}}.
 Suppose that $(u^N)_{N \geqslant 1}$ is a class of random
variables in $C ([0, T], E)$ with a given Polish space $E$. The sequence of $(u^N)_{N \geqslant 1}$
is tight in $C ([0, T], E)$ if and only if the following conditions hold:
\begin{enumerate}
	\item For each $\epsilon > 0$ and each $t \in [0, T]$, there is a compact set
	$A \subset E$ (possibly depending on $t$) such that
	\begin{align}
		\sup_N \mathbb{P} (u^N_t \in A) > 1 - \epsilon . \nonumber
	\end{align}
	\item For each $\epsilon > 0$,
	
	\begin{align}
		\lim_{h \rightarrow 0} \sup_N \mathbb{P} (\sup_{s, t \in [0, T]} \sup_{|t -
			s| \leq h} \|u^N_t - u^N_s \|_E > \epsilon) = 0. \nonumber
	\end{align}
\end{enumerate}
Since the embedding $H^{- \alpha'} \hookrightarrow H^{- \alpha}$ is compact if
$\alpha' < \alpha$ (see {\cite[Proposition 4.6]{triebel2006theory}}), using
Chebyshev's inequality, one can get the following sufficient conditions for
tightness in $C ([0, T], H^{- \alpha})$
	\begin{enumerate}[label=(\roman*)]
	\item For each $t \in [0, T]$, there exists some $\alpha' < \alpha$ such
	that
	\begin{equation}
		\sup_N \mathbb{E} \| u^N_t \|_{H^{- \alpha'}} < \infty . \label{aa1}
	\end{equation}
	\item There exists $\theta > 0$ such that
	\begin{equation}
		\sup_N \mathbb{E} \| u^N_t \|_{C^{\theta} ([0, T], H^{- \alpha})} = \sup_N
		\mathbb{E} \left( \sup_{ 0 \leq s <  t \leq T } \frac{\| u_t^N - u_s^N \|_{H^{-
					\alpha}}}{(t - s)^{\theta}} \right) < \infty . \label{aa2}
	\end{equation}
\end{enumerate}
Therefore to obtain tightness of $\{\eta^N,\mathcal{M}^N\}_{N\in \mathbb{N}}$ it suffices to justify (i) and (ii) with $\mathcal{M}^N$ and $\eta^N$ playing the role of $u^N$.

The following lemma gives tightness of the martingale part.
\begin{lemma}\label{lemma m}
	For every $\alpha > d / 2 + 1$, the sequence of  $(\mathcal{M}^N)_{N
		\geqslant 1}$ is tight in the space $C ([0, T], H^{- \alpha})$.
\end{lemma}
\begin{proof}By the above tightness criterion, it is indeed  sufficient to prove that: for each $\alpha > d / 2 + 1$ and $\theta'\in (0,\frac{1}{2})$,  it holds that
	\begin{equation*}
		\sup_N \mathbb{E} (\| \mathcal{M}^N
		\|_{C^{\theta'} ([0, T], H^{- \alpha})}^2) < \infty.
	\end{equation*}
First, for the Fourier basis $\{ e_k \}_{k \in \mathbb{Z}^d}$ and $t \in
	[0, T]$, we find
	\begin{equation}
		\mathcal{M}_t^N (e_k) = \frac{\sqrt{2\sigma_N}}{\sqrt{N}}  \sum_{i = 1}^N \int^t_0 \nabla e_k (X_{i}) \cdot
		\mathd B_s^i =  \sqrt{- 1} \frac{\sqrt{2\sigma_N}}{\sqrt{N}}   \sum_{i = 1}^N \int^t_0
		e_k (X_{i}) k \cdummy \mathd B_s^i . \label{aux1}
	\end{equation}
For any $\theta > 1$ and $0
	\leqslant s < t \leqslant T$, we deduce from H\"older's inequality that
	\begin{align*}
		\sup_N \mathbb{E} (\| \mathcal{M}_t^N -\mathcal{M}_s^N \|_{H^{-
				\alpha}}^{2 \theta}) = & \sup_N \mathbb{E} \bigg[ \bigg( \sum_{k \in
			\mathbb{Z}^d} \langle k \rangle^{- 2 \alpha} | \mathcal{M}_t^N (e_k)
		-\mathcal{M}_s^N (e_k) |^2 \bigg)^{\theta} \bigg]\\
		 \leqslant & \sup_N \mathbb{E} \bigg( \sum_{k \in \mathbb{Z}^d} \langle
		k \rangle^{- \alpha_1 \theta} | \mathcal{M}_t^N (e_k) -\mathcal{M}_s^N
		(e_k) |^{2 \theta} \bigg) \bigg( \sum_{k \in \mathbb{Z}^d} \langle k
		\rangle^{- \alpha_2 \frac{\theta}{\theta - 1}} \bigg)^{\theta - 1},
	\end{align*}
	where $\alpha_1 + \alpha_2 = 2 \alpha$ and $\alpha_1, \alpha_2 > 0$. Further choosing $\alpha_1=\alpha+1-\frac{d}{2}+\frac{d}{\theta}$ and $\alpha_2=\alpha-1+\frac{d}{2}-\frac{d}{\theta}$, then we have $(\alpha_1-2)\theta>d$ and  $\alpha_2 \frac{\theta}{\theta - 1} > d$, due to  $\theta > 1$ and the condition $\alpha > d / 2 + 1$. Hence
	the summation $\sum_{k \in \mathbb{Z}^d} \langle k \rangle^{- \alpha_2
		\frac{\theta}{\theta - 1}}$ is finite. Moreover, using the equality
	{\eqref{aux1}} gives
	\begin{eqnarray}
		\sup_N \mathbb{E} (\| \mathcal{M}_t^N -\mathcal{M}_s^N \|_{H^{-
				\alpha}}^{2 \theta}) & \lesssim_{\alpha_2, \theta} & \sup_N \sum_{k \in
			\mathbb{Z}^d} \langle k \rangle^{- \alpha_1 \theta} \mathbb{E} (|
		\mathcal{M}_t^N (e_k) -\mathcal{M}_s^N (e_k) |^{2 \theta}) \nonumber\\
		& \lesssim_{\alpha_2, \theta} & \sum_{k \in \mathbb{Z}^d} \langle k
		\rangle^{- \alpha_1 \theta + 2 \theta} \sup_N \mathbb{E} \left|  \frac{\sqrt{2\sigma_N}}{\sqrt{N}} \sum_{i = 1}^N \int^t_s e_k (X_i) \mathd B_r^i
		\right|^{2 \theta} \nonumber\\
		& \lesssim_{\alpha_2, \theta} & \sum_{k \in \mathbb{Z}^d} \langle k
		\rangle^{- \alpha_1 \theta + 2 \theta} \sup_N \mathbb{E} \left( \int^t_s
		\frac{2\sigma_N}{N}\sum_{i = 1}^N | e_k (X_i) |^2 \mathd r \right)^{\theta}
		\nonumber\\
		& \lesssim_{\alpha_2, \theta} & (t - s)^{\theta} \sum_{k \in
			\mathbb{Z}^d} \langle k \rangle^{- \alpha_1 \theta + 2 \theta}\lesssim_{\alpha_1,\alpha_2, \theta}(t-s)^{\theta} ,
		\label{tight m 1}
	\end{eqnarray}
	where the third inequality follows by the Burkholder-Davis-Gundy's inequality.
Therefore, {\eqref{tight m 1}} allows us to apply the  Kolmogorov
	continuity theorem {\cite[Theorem 2.3.11]{hofmanova2018spde}}, and we
	find
	\begin{equation}
		\sup_N \mathbb{E} (\| \mathcal{M}^N \|_{C_{}^{\theta'} ([0, T], H^{-
				\alpha})}^{2 \theta}) < \infty, \label{tight m 2}
	\end{equation}
	for any $0 < \theta' < \frac{\theta - 1}{2 \theta}$, $\theta > 1$, and
	$\alpha > d / 2 + 1$. The result follows by arbitrary $\theta>1$.
\end{proof}

Next, we need  the tightness of the fluctuation measures.

\begin{lemma}
	\label{lemma tight mu}Under the assumptions {\bf{(A2)-(A4)}}, for every
	$\alpha > d / 2 + 2$, the sequence of $(\eta^N)_{N \geqslant 1}$ is
	tight in the space $C ([0, T], H^{- \alpha})$.
\end{lemma}

\begin{proof}
	First, by Assumption {\bf{(A3)}} and Lemma \ref{lemma uni
		mu}, one can easily deduce \eqref{aa1} with $\eta^N$ playing the role of
	$u^N$ for any $\alpha > d / 2 + 2$. Indeed, taking $\mu_N (\cdummy) -
	\bar{\rho} = \frac{1}{\sqrt{N}} \eta^N$ into Lemma \ref{lemma uni mu}
	immediately gives
	\begin{equation}
	\sup_{ t \in [0, T]}	\sup_N \mathbb{E} \| \eta_t^N \|^2_{H^{- \alpha + 2}} \lesssim \sup_{ t \in [0, T]}\sup_N H
		(\rho_N | \bar{\rho}_N \nobracket) (t) + 1. \label{tight mu 1}
	\end{equation}
	Then {\bf{(A3)}} implies that the  right
	hand side of {\eqref{tight mu 1}} is finite. Thus \eqref{aa1} follows by
	$\alpha - 2$ playing the role of $\alpha'$.
	
	As for \eqref{aa2}, it suffices to prove the case $\alpha - 2 \in
	(d / 2, \beta)$, where $\beta$ is given in Assumption {\bf{(A5)}}. 
	 Recall the decomposition {\eqref{equation fluctuation}},
	$\| \eta_t^N - \eta^N_s \|_{H^{- \alpha}}$, $0 \leqslant s < t < T$, is
	controlled via the following relation
	\begin{equation}
		\| \eta_t^N - \eta^N_s \|^2_{H^{- \alpha}} \lesssim  \sum_{i = 1}^5 J_{s,
			t}^i \label{tight mu 2},
	\end{equation}
	where $J_{s, t}^i$, $i = 1, \ldots, 5$, are defined by
	\begin{align}
		J_{s, t}^1 \assign & \left\| \sigma \int^t_s \Delta \eta_r^N \mathd r
		\right\|^2_{H^{- \alpha}}, \quad\quad\quad
		J^2_{s, t} \assign  \left\| \int^t_s \mathcal{K}^N_r \mathd r
		\right\|^2_{H^{- \alpha}}, \nonumber\\
		J^3_{s, t} \assign & \left\| \int^t_s \nabla \cdummy (F \eta_r^N) \mathd r
		\right\|^2_{H^{- \alpha}}, \quad
		J^4_{s, t} \assign  \left\| \int^t_s \sqrt{N} (\sigma_N - \sigma) \Delta
		\mu_N (r) \mathd r \right\|^2_{H^{- \alpha}}, \nonumber\\
		J^5_{s, t} \assign & \| \mathcal{M}_t^N -\mathcal{M}_s^N \|^2_{H^{-
				\alpha}} . \nonumber
	\end{align}
	For $J_{s,
		t}^1$, applying H\"older's inequality gives
	\begin{equation}
		\sup_N \mathbb{E} \left( \sup_{0 \leqslant s < t \leqslant T} \frac{J_{s,
				t}^1}{t - s} \right) \lesssim \sup_N \mathbb{E} \int^T_0 \| \Delta
		\eta_t^N \|^2_{H^{- \alpha}} \mathd t \lesssim \sup_N \sup_{t \in [0, T]}
		\mathbb{E} \| \eta_t^N \|_{H^{- \alpha + 2}}^2 < \infty ,\
	\end{equation}
	where we used {\eqref{tight mu 1}} in the last step.
	
\noindent For $J_{s,
		t}^2$, similarly, applying H\"older's inequality gives
	
	\begin{align}
		\sup_N \mathbb{E} \left( \sup_{0 \leqslant s < t \leqslant T} \frac{J_{s,
				t}^2}{t - s} \right) \lesssim & \sup_N \mathbb{E} \left( \int^T_0
		\|\mathcal{K}^N_t \|^2_{H^{- \alpha}} \mathd t \right) \nonumber\\
		\lesssim & \sup_N \sup_{t \in [0, T]} \mathbb{E} \|\mathcal{K}^N_t
		\|^2_{H^{- \alpha}} . \nonumber
	\end{align}
	Recall that $\mathcal{K}^N_t = \sqrt{N} \nabla \cdummy [K \ast \mu_N (t)
	\mu_N (t) - \bar{\rho}_t K \ast \bar{\rho}_t]$, and thus Lemma \ref{lemma uni
		nonlinear} and the assumptions {\bf{(A2)-(A3)}} deduces that
	\begin{equation}
		\sup_N \mathbb{E} \left( \sup_{0 \leqslant s < t \leqslant T} \frac{J_{s,
				t}^2}{t - s} \right) \lesssim \sup_N \sup_{t \in [0, T]} H (\rho_N |
		\bar{\rho}_N \nobracket) + 1 < \infty .
	\end{equation}
	Similarly, we obtain
	\begin{align}
		\sup_N \mathbb{E} \left( \sup_{0 \leqslant s < t \leqslant T} \frac{J_{s,
				t}^3}{t - s} \right) \lesssim & \sup_N \sup_{t \in [0, T]} \mathbb{E} \|F
		\eta_t^N \|_{H^{- \alpha + 2}}^2, \nonumber
	\end{align}
		and
	\begin{align}
		\sup_N \mathbb{E} \left( \sup_{0 \leqslant s < t \leqslant T} \frac{J_{s,
				t}^4}{t - s} \right) \lesssim & \sup_N \sup_{t \in [0, T]} N | \sigma_N -
		\sigma |^2 \mathbb{E} \| \mu_N (t)\|_{H^{- \alpha + 2}}^2 . \nonumber
	\end{align}
	Furthermore, Lemma \ref{lemma embedding} together with  Lemma \ref{lemma triebel} shows that
	\begin{equation*}
		\|F\eta\|_{H^{-\alpha+2}}\lesssim \|F\|_{C^{\beta}} \| \eta^N
	\|_{H^{- \alpha + 2}}.
	\end{equation*}
Hence using  {\eqref{tight mu 1}} and Assumption {\bf{(A4)}} gives
	\begin{align}
		\sup_N \mathbb{E} \left( \sup_{0 \leqslant s < t \leqslant T} \frac{J_{s,
				t}^3}{t - s} \right) \lesssim & \sup_N \sup_{t \in [0, T]}
		{\|F\|_{C^{\beta}}^2} \mathbb{E} \| \eta_t^N
		\|_{H^{- \alpha + 2}}^2 < \infty .  \label{tight mu 3}
	\end{align}
On the other hand, Assumption {\bf{(A4)}}, $\mu_N=N^{-1/2}\eta^N+\bar\rho$ and  {\eqref{tight mu 1}} imply that
	\begin{align}
		\sup_N \mathbb{E} \left( \sup_{0 \leqslant s < t \leqslant T} \frac{J_{s,
				t}^4}{t - s} \right) \rightarrow 0. \label{tight mu 4}
	\end{align}
		For $J_{s,
		t}^5$, we deduce from {\eqref{tight m 2}} that for any $\theta \in
	\left( 0, \frac{1}{2} \right)$,
	\begin{equation}
		\sup_N \mathbb{E} \left( \sup_{0 \leqslant s < t \leqslant T} \frac{J_{s,
				t}^5}{(t - s)^{2 \theta}} \right) = \sup_N \mathbb{E} (\| \mathcal{M}^N
		\|_{C_{}^{\theta} ([0, T], H^{- \alpha})}^2) < \infty . \label{tight mu 5}
	\end{equation}
	We are in a position to conclude \eqref{aa2} with $\eta^N$
	playing the role of $u^N$ for any $\alpha > d / 2 + 2$, and tightness of
	the sequence $(\eta^N)_{N \geqslant 1}$ follows. Indeed, combining {\eqref{tight
			mu 2}}-{\eqref{tight mu 5}} yields that
	
	\begin{align}
		\sup_N \mathbb{E} (\| \eta^N \|_{C^{\theta} ([0, T], H^{-
				\alpha})}^2) &=  \sup_N \mathbb{E} \left( \sup_{0 \leqslant s < t
			\leqslant T} \frac{\| \eta_t^N - \eta^N_s \|^2_{H^{- \alpha}}}{(t - s)^{2
				\theta}} \right) \nonumber\\
	&	\lesssim_T   \sup_N  \sum^5_{i = 1} \mathbb{E} \left( \sup_{0 \leqslant s <
			t \leqslant T} \frac{J_{s, t}^i}{(t - s)^{2 \theta}} \right) < \infty,
		\nonumber
	\end{align}
	for any $\theta \in \left( 0, \frac{1}{2} \right)$. The result then follows.
\end{proof}

\begin{remark}
Careful readers may find that it suffices to assume  $|\sigma_N- \sigma | = o(\frac{1}{\sqrt{N}})$ in order to obtain the tightness of $(\eta^N)$. But we still adopt the assumption that $ |\sigma_N- \sigma| = \mathcal{O} \left( \frac{1}{N} \right) $ in  Assumption {\bf{(A4)}} and Assumption {\bf{(A5)}} since this is one of the assumptions used in \cite{jabin2018quantitative} to obtain the uniform  bound for $H(\rho_N \vert \bar \rho^{\otimes N })$, i.e. our Assumption {\bf{(A3)}}.
\end{remark}

Define the topological space $\mathcal{X}$:
\begin{align*}
	\mathcal{X} \assign  \left\{\bigcap_{k\in \mathbb{N}}  \left[ C ([0, T],
H^{- \frac{d}{2}-2-\frac{1}{k}}) \cap L^2([0,T],H^{-\frac{d}{2}-\frac{1}{k}})\right] \right\} \times
\left\lbrace \bigcap_{k\in \mathbb{N}}C ([0, T], H^{- \frac{d}{2}-1-\frac{1}{k}})\right\rbrace  .
\end{align*}
The space $Y:=\cap_{k\in \mathbb{N}}Y_k $ with $C ([0, T],
H^{- \frac{d}{2}-2-\frac{1}{k}}) \cap L^2([0,T],H^{-\frac{d}{2}-\frac{1}{k}})$ or $ C ([0, T], H^{- \frac{d}{2}-1-\frac{1}{k}})$ playing the role of  $Y_k$ is endowed with the metric
$	d_{Y }(f,g)=\sum_{k=1}^{\infty}2^{-k}(1\wedge\|f-g\|_{Y_k})$.
Thus the convergence in $Y $ is equivalent to the convergence in $Y_k$ for every $k\in \mathbb{N}$. Moreover, $\mathcal{X}$ is a Polish space.

We then deduce the following result by the Skorokhod theorem.
\begin{theorem}\label{thm:skorokhod}
There exists a subsequence of $(\eta^N,
	\mathcal{M}^N)_{N \geqslant 1}$, still denoted by $(\eta^{N}, \mathcal{M}^{N})$ for simplicity,
	and a probability space $(\tilde{\Omega}, \tilde{\mathcal{F}},
	\tilde{\mathbb{P}})$ with $\mathcal{X}$-valued random variables
	$(\tilde{\eta}^{N_{}}, \tilde{\mathcal{M}}^N)_{N \geqslant 1}$ and
	$(\tilde{\eta}, \tilde{M})$ such that
	\begin{enumerate}
		\item For each $N \in \mathbb{N}$, the law of $(\tilde{\eta}^N,
		\tilde{\mathcal{M}}^N)$ coincides with the law of $(\eta^{N},
		\mathcal{M}^{N})$.
		\item The sequence of $\mathcal{X}$-valued random variables $(\tilde{\eta}^N,
		\tilde{\mathcal{M}}^N)_{N \geqslant 1}$ converges to $(\tilde{\eta},
		\tilde{\mathcal{M}})$ in $\mathcal{X}$  $\tilde{	\mathbb{P}}$-a.s.
	\end{enumerate}
\end{theorem}

\begin{proof}By  the Skorokhod theorem, the result follows by
 justifying the fact that the joint law of $(\eta^N,
\mathcal{M}^N)_{N \geqslant 1}$ is tight on $\mathcal{X}$.

We start with proving the set $A$ defined below is relatively compact in the space $C ([0, T],
H^{- \alpha-2}) \cap L^2([0,T],H^{-\alpha})$ for each $\alpha>d/2$,
\begin{align*}
	A\assign \left \{ u \in K;\int_0^T\|u(t)\|_{H^{-\frac{2\alpha+d}{4}}}^2\mathd t\leq M,  \right\},
\end{align*}
where $K$  is relatively compact in $C ([0, T],
	H^{- \alpha-2})$.
Suppose a sequence $\{u_n\}\subset A$, then there is a subsequence $\{u_{n_m}\}$  converging in $C ([0, T],
H^{- \alpha-2})$. On the other hand, by  the Sobolev interpolation theorem
{\cite[Proposition 1.52]{chemin2011fourier}}, we find for $n,n'\in \mN$ and $-\alpha-2<-\alpha<-\frac{2\alpha+d}4$
\begin{align}
&	\int^T_0 \| u_{n}(t) -u_{n'}(t)\|^2_{H^{- \alpha}} \mathd t \nonumber
	\\&\leqslant
	\int^T_0 \| u_{n}(t) -u_{n'}(t) \|^{2 \theta}_{H^{- \frac{2\alpha+d}{4}}}  \| u_{n}(t) -u_{n'}(t) \|^{2 (1 - \theta)}_{H^{- \alpha-2}} \mathd t \nonumber\\
&\leqslant  \left( \int^T_0 \| u_{n}(t) -u_{n'}(t) \|^2_{H^{-  \frac{2\alpha+d}{4}}}
\mathd t \right)^{\theta} \left( \int^T_0 \| u_{n}(t) -u_{n'}(t) \|^2_{H^{-
		\alpha-2}} \mathd t \right)^{1 - \theta} . \nonumber\\
&\leqslant \left( \int^T_0 \| u_{n}(t) -u_{n'}(t) \|^2_{H^{-  \frac{2\alpha+d}{4}}} \mathd t \right)^{ \theta} \left( T \nocomma \sup_{t \in [0, T]} \| u_{n}(t) -u_{n'}(t)
\|^2_{H^{- \alpha-2}} \right)^{1-\theta} , \nonumber
\end{align}
where the interpolation constant $\theta\in (0,1)$ depends  on $\alpha$ and $d$.  This implies the convergence of the subsequence $\{u_{n_m}\}$ in $ L^2([0,T],H^{-\alpha})$ for each $\alpha>d/2$, and $A$ is thus relatively compact in $C ([0, T],
H^{- \alpha-2}) \cap L^2([0,T],H^{-\alpha})$. For each $\epsilon>0$,  by \eqref{tight mu 1} and Lemma \ref{lemma tight mu},  one can find  $M$ sufficiently large and a compact set $K$ in $C ([0, T],
H^{- \alpha-2}) $ such that
\begin{align*}
	\mathbb{P}(\eta^N\notin A)\leq& P(\eta^N\notin K) +\mathbb{P}\(\int_0^T\|\eta^N(t)\|_{H^{-\frac{2\alpha+d}{4}}}^2\mathd t> M \)
	\\\leq& P(\eta^N\notin K) +\frac{T}{M}\sup_{ t \in [0, T]}\sup_N\mathbb{E}\|\eta^N(t)\|_{H^{-\frac{2\alpha+d}{4}}}^2<\epsilon,
\end{align*}
where the second line follows by Chebyshev's inequality. Therefore  the sequence of laws of $(
\mathcal{\eta}^N)_{N \geqslant 1}$ is tight on $C ([0, T],
H^{- \alpha-2}) \cap L^2([0,T],H^{-\alpha})$ for every $\alpha>d/2$.

Furthermore, recall that  Lemma \ref{lemma m} gives that the sequence of laws of $(
\mathcal{M}^N)_{N \geqslant 1}$ is tight on $C([0,T],H^{-\alpha-1} )$ for every $\alpha>d/2$.
For each $\epsilon>0$ and  $k\in \mathbb{N}$, choose compact sets $A_k^{\epsilon}$ and $B_k^{\epsilon}$ in $ C ([0, T],
H^{- \frac{d}{2}-\frac{1}{k}-2}) \cap L^2([0,T],H^{-\frac{d}{2}-\frac{1}{k}})$ and $ C([0,T],H^{-\frac{d}{2}-\frac{1}{k}-1} )$, respectively, such that
\begin{align*}
	\mathbb{P}(\eta^N\notin A_k^{\epsilon})< \epsilon 2^{-k},\quad 	\mathbb{P}(\mathcal{M}^N\notin B_k^{\epsilon})< \epsilon 2^{-k},\quad \forall N\in \mathbb{N}.
\end{align*}
Thus the set $A^{\epsilon}\times B^{\epsilon}$ in $\mathcal{X}$ defined by
\begin{align*}
	A^{\epsilon}\times B^{\epsilon}\assign \left(\bigcap_{k\in \mathbb{N}}A_k^{\epsilon}\right)\times\left(\bigcap_{k\in \mathbb{N}}B_k^{\epsilon}\right)
\end{align*}
is relatively compact 
and satisfies
\begin{align*}
	\mathbb{P}\((\eta^N,
	\mathcal{M}^N)\notin A^{\epsilon}\times B^{\epsilon}\)\leq \sum_{k\in\mathbb{N}}\mathbb{P}(\eta^N\notin A_k^{\epsilon})+\mathbb{P}(\mathcal{M}^N\notin B_k^{\epsilon})<2\epsilon,\quad \forall N\in \mathbb{N},
\end{align*}
which shows the tightness of   $(\eta^N,
\mathcal{M}^N)_{N \geqslant 1}$ in $\mathcal{X}$.
\end{proof}

\begin{corollary}
	\label{prop converge}For every $\alpha > d / 2$, it holds that
	\begin{equation}
		\tilde{\mathbb{E}} \int^T_0 \| \tilde{\eta}_t^N - \tilde{\eta}_t \|_{H^{- \alpha}} \mathd t
		\xrightarrow{N \rightarrow \infty} 0. \label{coro 1}
	\end{equation}
\end{corollary}

\begin{proof}
	Notice that
	\begin{equation}
	\tilde{\mathbb{E}}\int_0^T \| \tilde{\eta}_t \|^2_{H^{- \alpha}} \mathd t\leq 	\sup_N\tilde{\mathbb{E}}\int_0^T \| \tilde{\eta}_t^N \|^2_{H^{- \alpha}} \mathd t\leqslant
		T\sup_{t \in [0, T]} \sup_N \tilde{\mathbb{E}} \| \tilde{\eta}_t^N \|^2_{H^{- \alpha}} <
		\infty, \quad \forall \alpha > \frac{d}{2},\nonumber
	\end{equation}
which provides the uniform (in $[0, T]
	\times \Omega$) integrability of $\| \tilde{\eta}_t^N - \tilde{\eta}_t \|_{H^{- \alpha}}$,
	thus  the  convergence of $\| \tilde{\eta}_t^N -
	\tilde{\eta}_t \|_{H^{- \alpha}}$  $\mathd t \times \mathd\tilde{\mathbb{P}}$-a.e. leads to {\eqref{coro 1}}.
\end{proof}

For each $N$, let $(\tilde{\mathcal{F}}_t^N)_{t \geqslant 0}$ and
$(\tilde{\mathcal{F}}_t)_{t \geqslant 0}$ be the normal filtration generated
by $(\tilde{\eta}^N, \tilde{\mathcal{M}}^N)$ and $(\tilde{\eta},
\tilde{\mathcal{M}})$, respectively. Then we have
\begin{align}
	\tilde{\mathcal{M}}^N_t=\tilde{\eta}_{t}^N-\tilde{\eta}_{0}^N-\sigma\int_0^t\Delta\tilde{\eta}_s^N\mathd s+\int_0^t\tilde{\mathcal{K}}_s^N\mathd s+\int_0^t\nabla(F\tilde{\eta}_s^N)\mathd s-\tilde{R}_t^N,\label{equation mn}
\end{align}
where $\tilde{\mathcal{K}}^N$ and  $\tilde{R}^N$ are defined  with $\tilde{\mu}_N:=\bar{\rho}+\frac{1}{\sqrt{N}}\tilde{\eta}^N$ and
\begin{align*}
	\tilde{\mathcal{K}}_t^N\assign\sqrt{N}\nabla \cdot \(K*\tilde{\mu}_N(t)\tilde{\mu}_N(t)-K*\bar{\rho}_t\bar{\rho}_t\),\quad \tilde{R}^N_t\assign\sqrt{N}(\sigma_N-\sigma)\int_0^t\Delta\tilde{\mu}_N(s)\mathd s.
\end{align*}
Here $\tilde{\mathcal{K}}^N$ is well-defined since $\tilde{\mu}_N$  is  linear combination of Dirac measure and $K*\tilde{\mu}_N\tilde{\mu}_N$ is understood as
\begin{align*}
	\<K*\tilde{\mu}_N\tilde{\mu}_N,\varphi\>=\int_{\mathbb{T}^{d}\times \mathbb{T}^{d}}K(x-y)\varphi(x)\tilde{\mu}_N(\mathd x)\tilde{\mu}_N(\mathd y),
\end{align*}
for $\varphi\in C^1$.





\subsection{Characterization of the limit}\label{sec:chara}

In this section, we conclude that the original sequence $(\eta^N)_{N \geqslant 1}$ converges in distribution to the  equation
{\eqref{equation li}}. Recall that the sequence $(\tilde{\eta}^N)_{N \geqslant 1}$ converges in $C([0,T], H^{-\alpha-2})\cap L^2([0,T],H^{-\alpha})$   $\tilde{	\mathbb{P}}$-a.s.  for $\alpha>d/2$ and shares the
same distribution with a subsequence of $(\eta^N)_{N \geqslant
	1}$. Hence it is sufficient to justify two facts. One is that each limit
$\tilde{\eta}$ is a martingale solution to {\eqref{equation li}}.
The other  is that the
law of the solution to {\eqref{equation li}} is unique, which would
follow by pathwise uniqueness and the Yamada-Watanabe theorem.

Throughout  this section, we always assume {\bf{(A1)-(A4)}}.

Identifying the limit of the interacting term $\tilde{\mathcal{K}}^N$ is one of the main difficulties in this article, it derseves to be  treated separately from other terms in the decomposition \eqref{equation fluctuation}. The following lemma identifies  the limit of the interacting term $\tilde{\mathcal{K}}^N$. The idea of  the proof
is to  split the interacting term into some regular part
and a term in the form of  a function of $\mu_N  - \bar{\rho}$,
which can be controlled in Lemma \ref{lemma uni product} by the techniques developed in Section \ref{sec:estimates}.

\begin{lemma}
	\label{lemma limit no}For each $\varphi \in C^{\infty} (\mathbb{T}^d)$, it
	holds that
	
	\begin{align}
		\tilde{\mathbb{E}} \left( \sup_{t \in [0, T]} \left| \int^t_0 \tilde{\mathcal{K}}^N_s
		(\varphi) - \langle \bar{\rho}_s K \ast \tilde{\eta}_s + \tilde{\eta}_s K \ast
		\bar{\rho}_s, \nabla \varphi \rangle \mathd \nocomma s \right| \right)
		\xrightarrow{N \rightarrow \infty} 0. \nonumber
	\end{align}
\end{lemma}

\begin{proof} Direct computations give the following identity
	\begin{equation*}
		\sqrt{N} \big(\tilde\mu_N K * \tilde\mu_N - \bar\rho K * \bar \rho \big) = \bar \rho K * \eta^N + \eta^N K * \bar \rho + \frac{1}{\sqrt{N}} \, \eta^N K * \eta^N.
	\end{equation*}
Consequently,  for each $\varphi \in C^{\infty}$,
	\begin{equation}
		\sup_{t \in [0, T]} \left| \int^t_0 \tilde{\mathcal{K}}^N_s (\varphi) - \langle
		\bar{\rho}_s K \ast \tilde{\eta}_s + \tilde{\eta}_s K \ast \bar{\rho}_s, \nabla \varphi
		\rangle \mathd \nocomma s \right| \leqslant J_1^N (\varphi) + J_2^N
		(\varphi), \label{limit no 1}
	\end{equation}
	where
	
	\begin{align}
		J_1^N (\varphi) \assign & \sqrt{N} \int^T_0 | \langle \nabla \varphi K \ast
		(\tilde{\mu}_N (t) - \bar{\rho}_t), \tilde{\mu}_N (t) - \bar{\rho}_t \rangle | \mathd t,
		\nonumber\\
		J_2^N (\varphi) \assign& \int^T_0 | \langle \bar{\rho}_t K \ast \tilde{\eta}^N_t +
		\tilde{\eta}^N_t K \ast \bar{\rho}_t, \nabla \varphi \rangle - \langle
		\bar{\rho}_t K \ast \tilde{\eta}_t + \tilde{\eta}_t K \ast \bar{\rho}_t, \nabla \varphi
		\rangle | \mathd t. \nonumber
	\end{align}
	On one hand, we deduce from Lemma \ref{lemma uni product} that
	\begin{align*}
		\tilde{\mathbb{E}}J_1^N (\varphi) & \leqslant T \nocomma \sqrt{N} \sup_{t \in [0,
			T]} \tilde{\mathbb{E}} | \langle \nabla \varphi K \ast (\tilde{\mu}_N (t) - \bar{\rho}_t),
		\tilde{\mu}_N (t) - \bar{\rho}_t \rangle |
		\\ &= T \sqrt{N} \sup_{t \in [0,
			T]} {\mathbb{E}} | \langle \nabla \varphi K \ast ({\mu}_N (t) - \bar{\rho}_t),
	{\mu}_N (t) - \bar{\rho}_t \rangle
		\\
		& \lesssim N^{- \frac{1}{2}} \nocomma \sup_N \sup_{t \in [0, T]} (H
		(\rho_N | \bar{\rho}_N \nobracket) + 1) \xrightarrow{N \rightarrow \infty}
		0,
	\end{align*}
		where the limit follows by Assumption {\bf{(A3)}}. On the
	other hand, we find
	\begin{equation}
		\tilde{\mathbb{E}}J_2^N (\varphi) \leqslant \tilde{\mathbb{E}} \int^T_0 | \langle
		\bar{\rho}_t K \ast (\tilde{\eta}^N_t - \tilde{\eta}_t), \nabla \varphi \rangle | + |
		\langle (\tilde{\eta}^N_t - \tilde{\eta}_t) K \ast \bar{\rho}_t, \nabla \varphi \rangle |
		\mathd t. \label{limit no2}
	\end{equation}
	For each $t \in [0, T]$, it holds for every $\alpha \in( d / 2,\beta)$ that
	\begin{align}
		| \langle \bar{\rho}_t K \ast (\tilde{\eta}^N_t - \tilde{\eta}_t), \nabla \varphi \rangle
		| & = | \langle K(-\cdot) \ast (\bar{\rho}_t \nabla \varphi), \tilde{\eta}^N_t -
		\tilde{\eta}_t \rangle | \leqslant \| \tilde{\eta}^N_t - \tilde{\eta}_t \|_{H^{- \alpha}}  \|
		K(-\cdot) \ast (\bar{\rho}_t \nabla \varphi)\|_{H^{\alpha}}, \nonumber\\
		| \langle (\tilde{\eta}^N_t - \tilde{\eta}_t) K \ast \bar{\rho}_t, \nabla \varphi \rangle
		| & \leqslant \| \tilde{\eta}^N_t - \tilde{\eta}_t \|_{H^{- \alpha}}  \| \nabla \varphi
		\cdummy K \ast \bar{\rho}_t \|_{H^{\alpha}},\nonumber
	\end{align}
where
\begin{align}
	K(-\cdot) * g (x) := \int K(y-x) g(y) \mathd y.\label{notation:k}
\end{align}
Applying
	Lemma \ref{lemma convolution} with $p = p_1 = q = 2$ and Lemma \ref{lemma:infity} yields that
	\begin{align}
		| \langle \bar{\rho}_t K \ast (\tilde{\eta}^N_t - \tilde{\eta}_t), \nabla \varphi \rangle
		| & \lesssim \| \tilde{\eta}^N_t - \tilde{\eta}_t \|_{H^{- \alpha}} \|K\|_{L^1}  (\|
		\bar{\rho}_t \|_{H^{\alpha}} \| \nabla \varphi \|_{L^{\infty}} +\|
		\bar{\rho}_t \|_{L^{\infty}} \| \nabla \varphi \|_{H^{\alpha}}),
		\nonumber\\
		| \langle (\tilde{\eta}^N_t - \tilde{\eta}_t) K \ast \bar{\rho}_t, \nabla \varphi \rangle
		| & \lesssim \| \tilde{\eta}^N_t - \tilde{\eta}_t \|_{H^{- \alpha}} \|K\|_{L^1}  (\|
		\bar{\rho}_t \|_{H^{\alpha}} \| \nabla \varphi \|_{L^{\infty}} +\|
		\bar{\rho}_t \|_{L^{\infty}} \| \nabla \varphi \|_{H^{\alpha}}) .
		\nonumber
	\end{align}
Here the fact $K\in L^1$ follows by Assumption {\bf{(A2)}}.
	Then taking these two estimates into {\eqref{limit no2}}, applying Sobolev embedding
	$H^{\alpha} \hookrightarrow L^{\infty}$ with $\alpha>d/2$, we thus arrive at
	\begin{equation}
	\tilde{	\mathbb{E}}J_2^N (\varphi) \lesssim_{\varphi} \|K\|_{L^1}
	{\sup_{t \in [0, T]} \| \bar{\rho}_t \|_{H^{\alpha}}}
		\mathbb{E} \int^T_0 \| \tilde{\eta}^N_t - \tilde{\eta}_t \|_{H^{- \alpha}} \mathd t
		\xrightarrow{N \rightarrow \infty} 0,\nonumber
	\end{equation}
	where the limit follows by {\eqref{coro 1}}. Using inequality {\eqref{limit
			no 1}} and $\mathbb{E}J_1^N (\varphi) \rightarrow 0$, the proof is
	completed.
\end{proof}

\begin{remark}
	One may easily find that in the ``identifying the limit" part, we only need to assume that the relative entropy grow slower than the order $\sqrt{N}$, i.e. $H(\rho_N \vert \bar \rho_N) = o(\sqrt{N})$ as $N \to \infty$. However, in the tightness part we need a stronger assumption, namely our Assumption {\bf{(A3)}}. As a separate question, it would be interesting to show whether or not there exists some symmetric probability measure $\rho_N \in \mathcal{P}_{Sym}(S^N)$ such that $H(\rho_N \vert \bar \rho^{\otimes N}) = N^{\theta}$ with $\theta \in (0, 1)$, where $\bar \rho $ is a given probability measure on the Polish space $S$.
\end{remark}

Now we are in the position to conclude that $\tilde{\eta}$ solves \eqref{equation li}.
\begin{theorem}
	\label{theorem existence}The limit $\tilde{\eta}$ is a
	martingale solution to {\eqref{equation li}} in the sense of Definition \ref{def:mar}.
\end{theorem}

\begin{proof}
	We deduce from \eqref{equation mn} that
	\begin{align}
	\tilde{\mathcal{M}}_t^N (\varphi) =&	\langle \tilde{\eta}^N_t, \varphi \rangle - \langle \tilde{\eta}^N_0, \varphi \rangle -
		\sigma \int_0^t \langle \Delta \varphi, \tilde{\eta}^N_s \rangle \mathd \nocomma s
		- \int_0^t \tilde{\mathcal{K}}_s^N (\varphi) \mathd \nocomma s - \int_0^t \langle
		\nabla \varphi, F \tilde{\eta}_s^N \rangle \mathd s \nonumber\\
		& - \sqrt{N}  (\sigma_N - \sigma)  \int^t_0
		\langle \Delta \varphi, \tilde{\mu}_N (s) \rangle \mathd s, \nonumber
	\end{align}	
	for each $\varphi \in C^{\infty} (\mathbb{T}^d)$ and $t \in [0, T]$. 
	By
	Lemma \ref{lemma limit no}, $\sigma_N-\sigma=\mathcal{O} \left( \frac{1}{N}
	\right)$,  and the fact that $\tilde{\eta}^N$
	converges to $\tilde{\eta}$ in $C ([0, T], H^{- \alpha - 2})\cap L^2 ([0, T], H^{- \alpha})$ for every $\alpha > d / 2$ $\tilde{\mathbb{P}}$-a.s., one can take limit of  every term above on both sides and have
	\begin{align}
\tilde{\mathcal{M}}_t(\varphi)=\left\langle \tilde{\eta}_{t},\varphi\right\rangle -\left\langle \tilde{\eta}_{0},\varphi\right\rangle -\sigma\int_0^t\left\langle \Delta\varphi,\tilde{\eta}_s\right\rangle \mathd s-\int_0^t\left\langle \nabla\varphi, \bar{\rho}_sK*\tilde{\eta}_s+\tilde{\eta}_sK*\bar{\rho}_s+F\tilde{\eta}_s\right\rangle \mathd s,\quad\tilde{\mathbb{P}}-a.s. \nonumber
 \end{align}
	To indentify $\tilde{\eta}$ is a martingale solution, we need to justify properties of $\tilde{\mathcal{M}}$. Since $\tilde{\mathcal{M}}^N$ are martingales w.r.t. the normal filtration $(\tilde \cF_t^N)$ generated by $\tilde \eta^N$, and using Theorem \ref{thm:skorokhod}, the limit $\tilde{\mathcal{M}}$ is a martingale 
	with values in $H^{- \alpha - 1}$ for every $\alpha > d / 2$ w.r.t. the filtration generated by $\tilde \eta$. More precisely we have for $t\geq s\geq0$ and any bounded continuous function on $C ([0, s], H^{- \alpha - 2})\cap L^2 ([0, s], H^{- \alpha})$
	\begin{align*}
		\tilde {\mathbb{E}}[(\tilde{\mathcal{M}}_t(\varphi)-\tilde{\mathcal{M}}_s(\varphi))g(\tilde \eta|_{[0,s]})]=\lim_{N\to\infty}\tilde {\mathbb{E}}[(\tilde{\mathcal{M}}^N_t(\varphi)-\tilde{\mathcal{M}}_s^N(\varphi))g(\tilde \eta^N|_{[0,s]})]=0.
	\end{align*}
	
	
	
As for the covariance functions, on one hand, applying Burkholder-Davis-Gundy's inequality, we
	have for each $1 < \theta \leqslant 2$
	\begin{align*}
		\sup_N\tilde{\mathbb{E}} [\sup_{t \in [0, T]} | \tilde{\mathcal{M}}_t^N (\varphi) |^{2
			\theta}] & =\sup_N	\mathbb{E} [\sup_{t \in [0, T]} | {\mathcal{M}}_t^N (\varphi) |^{2
			\theta}]
		\\&\lesssim \sup_N\mathbb{E} \left( \int_0^T \sum_{i = 1}^N
		\frac{\sigma_N}{N}  | \nabla \varphi (X_i) |^2 \mathrm{~ d} t
		\right)^{\theta} =\sup_N\mathbb{E} \left( \int_0^T \sigma_N  \langle | \nabla
		\varphi |^2, \mu_N (t) \rangle \mathrm{d} t \right)^{\theta}\\
		& \lesssim_{\varphi} \sup_N\sigma_N^{\theta} T^\theta \| \nabla \varphi \|_{L^\infty}^{2 \theta } < \infty.
	\end{align*}
This implies uniform integrability of $| \tilde{\mathcal{M}}_t^N
	(\varphi) |^2$ for each $t \in [0, T]$.
	

On the other hand, we have that  $| \tilde{\mathcal{M}}^N_t
(\varphi_1) \tilde{\mathcal{M}}^N_s (\varphi_2) |$ converges to $| \tilde{\mathcal{M}}_t
(\varphi_1)\tilde{ \mathcal{M}}_s (\varphi_2) |$ $\tilde{	\mathbb{P}}$-a.s. for $s, t \in [0,
T]$ and $\varphi_1, \varphi_2 \in C^{\infty}$. Thus by the uniform integrability of $| \tilde{\mathcal{M}}^N_t
(\varphi)|^2$ and
\begin{equation}
	|\tilde{ \mathcal{M}}^N_t (\varphi_1)\tilde{ \mathcal{M}}^N_s (\varphi_2) | \leqslant |
	\tilde{\mathcal{M}}^N_t (\varphi_1) |^2 + |\tilde{ \mathcal{M}}^N_s (\varphi_2) |^2,\nonumber
\end{equation}
we arrive at
\begin{align}\tilde{ \mathbb{E}} [\tilde{\mathcal{M}}_t (\varphi_1) \tilde{\mathcal{M}}_s (\varphi_2)] =
	\lim_{N \rightarrow \infty} \tilde{\mathbb{E}} [\tilde{\mathcal{M}}^N_t (\varphi_1)
	\tilde{\mathcal{M}}^N_s (\varphi_2)] =		\lim_{N \rightarrow \infty}\mathbb{E} [\mathcal{M}^N_t (\varphi_1) \mathcal{M}^N_s (\varphi_2)]. \nonumber\end{align}
Furthermore, using  \eqref{equation realization} and Ito's isometry
we obtain that

\begin{align}
	\mathbb{E} [\mathcal{M}^N_t (\varphi_1) \mathcal{M}^N_s (\varphi_2)] = &
	\frac{2 \sigma_N}{N} \mathbb{E} \left[ \left( \sum_{i = 1}^N \int^t_0
	\nabla \varphi_1 (X_i) \mathd B_r^i \right) \left( \sum_{i = 1}^N \int^s_0
	\nabla \varphi_2 (X_i) \mathd B_r^i \right) \right] \nonumber\\
	= & \frac{2 \sigma_N}{N} \mathbb{E} \left[ \sum_{i = 1}^N \int^{s \wedge
		t}_0 \nabla \varphi_1 (X_i) \nabla \varphi_2 (X_i) \mathd r \right]
	\nonumber\\
	= & 2 \sigma_N \int^{s \wedge t}_0 \mathbb{E} \langle \nabla \varphi_1 \cdot
	\nabla \varphi_2, \mu_N (r) \rangle \mathd r. \nonumber
\end{align}
Since
\begin{align} 2 \sigma_N \int^{s \wedge t}_0 \mathbb{E} \langle \nabla \varphi_1  \cdot \nabla
	\varphi_2, \mu_N (r) \rangle \mathd r \xrightarrow{N \rightarrow \infty}
	2 \sigma \int^{s \wedge t}_0 \langle \nabla \varphi_1  \cdot \nabla \varphi_2,
	\bar{\rho}_r \rangle \mathd r, \nonumber\end{align}
	we obtain 
	\begin{align*}\tilde{ \mathbb{E}} [\tilde{\mathcal{M}}_t (\varphi_1) \tilde{\mathcal{M}}_s (\varphi_2)]=	2 \sigma \int^{s \wedge t}_0 \langle \nabla \varphi_1  \cdot \nabla \varphi_2,
		\bar{\rho}_r \rangle \mathd r,
		\end{align*}
		which  implies that $\tilde M$ is a Gaussian process. The proof is complete.  
\end{proof}

The rest of this subsection is devoted to obtain the well-posedness of the SPDE \eqref{equation li}, and finish the proof of Theorem \ref{thm:1}.

Let us first introduce an equivalent definition of martingale  solutions to \eqref{equation li}, which is used in the  proof of Theorem \ref{thm:1}.  For notations'  simplicity, we omit the tildes in the following.
\begin{definition}
	We call $(\eta,\mathcal{M})$ a probabilistically weak solution to \eqref{equation li} on stochastic basis $(\Omega, \mathcal{F},
	\mathcal{F}_t, \mathbb{P})$ with initial data $\eta_0$ if
	\begin{enumerate}
		\item $\eta$  is a continuous $(\mathcal{F}_t)$-adapted process with values
		in $H^{- \alpha - 2}$  and  $\eta\in L^2 ([0, T], H^{-
			\alpha}) $  for every $\alpha > d / 2$,  $\mathbb{P}$-a.s.
		\item$\mathcal{M}$ is a continuous $(\mathcal{F}_t)$-adapted centered Gaussian process
		with values  in $H^{- \alpha - 1}$ for every $\alpha > d / 2$, with
		covariance given by
		\begin{equation}
			\mathbb{E} [\mathcal{M}_t (\varphi_1) \mathcal{M}_s (\varphi_2)] = 2
			\sigma \int^{s \wedge t}_0 \langle \nabla \varphi_1 \cdot\nabla \varphi_2,
			\bar{\rho}_r \rangle \mathd r,\label{eqt:mco}
		\end{equation}
	for each $\varphi_1,\varphi_2\in C^{\infty}$ and $s,t\in [0,T]$.
		\item For each $\varphi \in C^{\infty} (\mathbb{T}^d)$ and $t \in [0, T]$, it
		holds that
		
		\begin{align}
			\mathcal{M}_t (\varphi)=&\langle \eta_t, \varphi \rangle - \langle \eta_0, \varphi \rangle -
			\int^t_0 \langle \sigma \Delta \varphi, \eta \rangle \mathd \nocomma s -
			\int^t_0 \langle \nabla \varphi, \bar{\rho} K \ast \eta \rangle \mathd
			\nocomma s - \int^t_0 \langle \nabla \varphi, \eta K \ast \bar{\rho}
			\rangle \mathd \nocomma s \nonumber\\
			& - \int_0^t \langle \nabla \varphi, F \eta \rangle \mathd s
			. \nonumber
		\end{align}
	\end{enumerate}

Furthermore, given a centered Gaussian process $\mathcal{M}$ on stochastic basis $(\Omega, \mathcal{F},
\mathcal{F}_t, \mathbb{P})$  with covariance characterized by \eqref{eqt:mco}, we call $\eta$ is a probabilistically strong solution to \eqref{equation me} if $(\eta,\mathcal{M})$ is a probabilisttically weak solution and $\eta$ is adapted to the normal filtration generated by $\mathcal{M}$.
\end{definition}

Uniqueness  in law of the solutions to \eqref{equation li} usually follows by the Yamada-Watanabe theorem, which requires existence of probabilistically weak solutions and pathwise uniqueness.  Since the martingale solutions and the probabilistically weak solutions are equivalent, Theorem \ref{theorem existence} means that there exists a stochastic  basis $(\Omega, \mathcal{F},
\mathcal{F}_t, \mathbb{P})$  such that  $(\eta,\mathcal{M})$ is  a probabilistically weak solution to \eqref{equation li}, it thus  suffices to prove the pathwise uniqueness.

We now briefly explain the concept of pathwise uniqueness 
of  probabilistically weak solutions introduced before.  Equation \eqref{equation li} can be viewed as a system, for which the information of the initial data and the nosie  is given (i.e. the distribution of $(\mathcal{M},\eta_{0})$ is fixed), and 
$(\mathcal{M},\eta_{0})$ can be seen as the input and $\eta$ is the output.
Pathwise uniqueness means that if on some fixed stochastic basis there exist two  outputs  $\eta$ and $\tilde{\eta}$  with  given $\eta_{0}$ and $\mathcal{M}$, then $\eta$ coincides with $\tilde{\eta}$ $\mathbb{P}$-a.s..

Notice that the covariance function of $\mathcal{M}$ 
and Assumption {\bf{(A1)}} have determined the  distribution of $(\mathcal{M},\eta_{0})$.  Since   equation \eqref{equation li} is linear and is driven by additive noise, pathwise uniqueness of solutions to the equation \eqref{equation li} follows from  uniqueness of solutions to the following PDE
\begin{equation}
	\partial_t u = \sigma \Delta u - \nabla \cdot (\bar{\rho} K \ast u) - \nabla
	\cdot (uK \ast \bar{\rho}) - \nabla \cdot (Fu) ,\quad u_0=0. \label{equation pde 1}
\end{equation}
\begin{lemma}
	\label{lemma unique pde}Under the assumptions  {\bf{(A2)}} and {\bf{(A4)}} 
	with parameter $\beta$, for each $\alpha \in (d / 2, \beta)$, $u
	\equiv 0$ is the only solution with zero initial value to {\eqref{equation pde 1}} in the sense that
	\begin{enumerate}
		\item $u\in L^2 ([0, T], H^{- \alpha})\cap C([0,T],H^{-\alpha-2})$.
		\item For each $\varphi \in C^{\infty}$ and  $t\in [0,T]$,
		\begin{equation}
			\langle u_t, \varphi \rangle =  \int^t_0
			\langle \sigma u_s, \Delta \varphi \rangle \mathd s + \int^t_0 \langle
			\bar{\rho}_s K \ast u_s + u_s K \ast \bar{\rho}_s + F \nocomma u_s, \nabla
			\varphi \rangle \mathd s.\nonumber
		\end{equation}
	\end{enumerate}
\end{lemma}

\begin{proof}
Testing $u$ with the
	Fourier basis $\{e_k \}_{k \in \mathbb{Z}^d}$, then we find for every $t \in
	[0, T]$ and $k \in \mathbb{Z}^d$,
	
	\begin{align}
		\partial_t | \langle u_t, e_k \rangle |^2 = & - 2 \sigma |k|^2
		\langle u_t, e_k \rangle \langle u_t, e_{- k} \rangle + \langle u_t, e_{-
			k} \rangle  [J_t^1 (k) + J_t^2 (k)] + \langle u_t, e_k \rangle  [J_t^1 (-
		k) + J_t^2 (- k)] \nonumber\\
		& + \sqrt{- 1} k \langle u_t, e_{- k} \rangle \langle Fu_t, e_k
		\rangle- \sqrt{- 1} k \langle u_t, e_{ k} \rangle \langle Fu_t, e_{-k}
		\rangle,  \label{unique 2}
		\end{align}
		where $J_t^1 (k)$ and $J_t^2 (k)$, for each $k \in \mathbb{Z}^d$, are
	defined by
	
	\begin{align*}
		J_t^1 (k) \assign & \sqrt{- 1} k \int_{\mathbb{T}^d} K \ast u_t (x)
		e_k (x)  \bar{\rho}_t (x) \mathd \nocomma x,\\
		J^2_t (k) \assign &  \sqrt{- 1} k \int_{\mathbb{T}^d} K \ast
		\bar{\rho}_t (x) e_k (x) u_t (x) \mathd \nocomma x.
	\end{align*}
	Integrating  \eqref{unique 2} over time, summing up over $k$ with
	weight $\langle k \rangle^{- 2 \alpha - 2}$, and applying Young's inequality
	yields that there exists a constant $C_{\epsilon}$ for each $\epsilon > 0$ such that
	
	\begin{align}
		& \sum_{k \in \mathbb{Z}^d} \langle k \rangle^{- 2 \alpha - 2} | \langle
		\nocomma u_t, e_k \rangle |^2 + 2\sigma  \sum_{k \in \mathbb{Z}^d}
		\langle k \rangle^{- 2 \alpha}  \int_0^t | \langle u_t, e_k \rangle |^2
		\mathd s \nonumber\\
		\leqslant & C_{\epsilon}  \sum_{k \in \mathbb{Z}^d} \langle k \rangle^{- 2
			\alpha - 2}  \int_0^t | \langle u_s, e_k \rangle |^2 \mathd s + \epsilon
		\sum_{k \in \mathbb{Z}^d} \langle k \rangle^{- 2 \alpha - 2}  \int_0^t
		|J_s^1 (- k) + J_s^2 (- k) |^2 \mathd s \nonumber\\
		& + \epsilon \sum_{k \in \mathbb{Z}^d} \langle k \rangle^{- 2 \alpha - 2}
		|k|^2  \int_0^t | \langle Fu_s, e_k \rangle |^2 \mathd s.  \label{unique
			1}
	\end{align}
	To make {\eqref{unique 1}} suitable for applying Gronwall's lemma, we first
	find estimates related to $J_t^{1} (k)$ and $J_t^2 (k)$,
	
	\begin{align}
		\sum_{k \in \mathbb{Z}^d} \langle k \rangle^{- 2 \alpha - 2} |J^1 (k) |^2
		& =  \sum_{k \in \mathbb{Z}^d} \langle k \rangle^{- 2 \alpha - 2}
		|k|^2 \langle K \ast u \bar{\rho}, e_k \rangle \langle K \ast u
		\bar{\rho}, e_{- k} \rangle \leqslant   \|K \ast u \bar{\rho}
		\|^2_{H^{- \alpha}}, \nonumber
	\end{align}
	and
	\begin{align}
		\sum_{k \in \mathbb{Z}^d} \langle k \rangle^{- 2 \alpha - 2} |J^2 (k) |^2
		& =  \sum_{k \in \mathbb{Z}^d} \langle k \rangle^{- 2 \alpha - 2}
		|k|^2 \langle K \ast \bar{\rho} u, e_k \rangle \langle K \ast \bar{\rho}
		u, e_{- k} \rangle \leqslant   \|K \ast \bar{\rho} u\|^2_{H^{-
				\alpha}} . \nonumber
	\end{align}
	Then applying Lemmas \ref{lemma embedding}  and \ref{lemma triebel} gives that
	
	\begin{align}
		\|K \ast u \bar{\rho} \|_{H^{- \alpha}} \leqslant C_{\alpha}  \|K \ast
		u\|_{H^{- \alpha}} \| \bar{\rho} \|_{C^{\beta}}, \quad \|K
		\ast \bar{\rho} u\|_{H^{- \alpha}} \leqslant C_{\alpha} \|u\|_{H^{-
				\alpha}}  \|K \ast \bar{\rho} \|_{C^{\beta}}. \nonumber
	\end{align}
	Furthermore, by Lemma \ref{lemma convolution}, we deduce
	
	\begin{align}
		\sum_{k \in \mathbb{Z}^d} \langle k \rangle^{- 2 \alpha - 2} |J^1 (k) |^2
		+ \sum_{k \in \mathbb{Z}^d} \langle k \rangle^{- 2 \alpha - 2} |J^2 (k)
		|^2 \leqslant C_{\alpha} \|u\|^2_{H^{- \alpha}} \|K\|^2_{L^1} \|
		\bar{\rho} \|^2_{C^{\beta}} .  \label{unique 5}
	\end{align}
	Similarly, we obtain
	
	\begin{align}
		\sum_{k \in \mathbb{Z}^d} \langle k \rangle^{- 2 \alpha - 2} |k|^2 |
		\langle Fu, e_k \rangle |^2 \leq \|Fu\|^2_{H^{- \alpha}} \leq C_{\alpha}
		{\color[HTML]{000000}\|F\|^2_{C^{\beta}}} \|u\|^2_{H^{-
				\alpha}} .  \label{unique 3}
	\end{align}
	Since $u \in L^2 ([0, T], H^{- \alpha})$, we obtain $\p_tu\in L^2([0, T], H^{- \alpha-2})$, which by Lions-Magenes Lemma implies $u\in C([0,T],H^{-\alpha-1})$.
	 Combining {\eqref{unique
			1}}-{\eqref{unique 3}} leads to
	
	\begin{align}
		& \|u_t \|^2_{H^{- \alpha - 1}} + 2 \sigma   \int_0^t \|u_s
		\|^2_{H^{- \alpha}} \mathd s \nonumber\\
		\leqslant & C_{\epsilon}  \int_0^t \|u_s \|^2_{H^{- \alpha - 1}} \mathd s
		+ \epsilon C_{\alpha}  \int_0^t \left( {\nobracket
			\|K\|^2_{L^1} \| \bar{\rho}_s \|^2_{C^{\beta}}
			+\|F\|^2_{C^{\beta}})} \|u_s \|^2_{H^{- \alpha}} \mathd s.
	 \right. \label{eqt:energy}
	\end{align}
	Choosing $\epsilon$ such that
	
	\begin{align}
		\epsilon C_{\alpha} \|K\|^2_{L^1} \(\sup_{ s \in [0, t]}\| \bar{\rho}_s \|^2_{C^{\beta}} + \|F\|^2_{C^{\beta}}\) < & 2\sigma ,
		\nonumber
	\end{align}
		then using Gronwall's inequality gives
	
	\begin{align}
		\|u_t \|^2_{H^{- \alpha - 1}} + \int_0^t \|u_s \|^2_{H^{- \alpha}} \mathd
		s & = 0. \nonumber
	\end{align}
	This completes the proof.
\end{proof}

\begin{proof}[Proof of Theorem \ref{thm:1} ]\label{proof:thm}
	We have proved that the sequence of laws of $\{\eta^N \}_{N \in \mathbb{N}}$
is tight and every tight limit 
is a
martingale solution to {\eqref{equation li}} ( Theorem \ref{theorem existence}). As a result,  existence of  martingale solutions (equivalently probabilistically weak solutions) follows.
On the other hand, Lemma
\ref{lemma unique pde} together with Corollary \ref{coro22} implies pathwise uniqueness of probabilistically weak
solutions. Then applying the general Yamada-Watanabe theorem  {\cite[Theorem
	1.5]{kurtz2014yamada}} gives that the law of martingale solutions starting from the same initial distribution  is
unique, and  every probabilistically weak solution  is a probabilistically strong solution. 
Therefore $\eta^N$ converges in distribution to the unique (in distribution) martingale solution $\eta$.
\end{proof}

\br
From the proof of Theorem \ref{thm:1}, we also obtain the well-posedness of probabilistically strong solutions to the SPDE \eqref{equation li}.
\er


\subsection{Optimal regularity}\label{sec:op}
In this  subsection we improve the regularity of $\eta$ by using the mild formulation and the smooth effect of the heat kernel.

Recall that   $\mathcal{M}$ is a centered Gaussian process with covariance given by
\begin{equation*}
	\mathbb{E} [\mathcal{M}_t (\varphi_1) \mathcal{M}_s (\varphi_2)] = 2
	\sigma \int^{s \wedge t}_0 \langle \nabla \varphi_1 \cdot\nabla \varphi_2,
	\bar{\rho}_r \rangle \mathd r,
\end{equation*}
for  $\varphi_1, \varphi_2 \in C^{\infty}(\mathbb{T}^d)$. Therefore, the distribution of $\mathcal{M}$ is uniquely determined, and one can regard $\mathcal{M}$ as $\nabla\cdot \int_0^{\cdot} \sqrt{\bar{\rho}}\xi(\dif s,\dif x)$ with $\xi=(\xi^i)_{i=1}^d$ being vector valued  space-time white noise  on $\mR^+\times \mT^d$.  
In fact, for every $\varphi\in C^{\infty}$, 
\begin{equation}
	\mathcal{M}_t(\varphi)\overset{d}{=}-\sqrt{2\sigma}\int_0^t\int_{\mathbb{T}^d}\nabla\varphi(x)\sqrt{\bar{\rho}_s(x)}\xi(\dif s,\dif x),\label{eqt:mar}
\end{equation}
where $\overset{d}{=}$ means equal in distribution and  we omit the inner product in $\mR^d$  between $\xi$ and $\nabla \varphi$.
We start with investigating the regularity of a stochastic integral, which will be the stochastic term in the  mild form of equation  \eqref{equation li}.
Define  a stochastic process $Z$ as
\begin{equation}\label{eq:Z}
	Z_t:=\int_0^t\int_{\mathbb{T}^d}\nabla\Gamma_{t-s}(\cdot-y)\sqrt{\bar{\rho}_s(y)}\xi(\dif s,\dif y),
\end{equation}
where $\Gamma $ is the heat kernel of $\sigma\Delta$ on $\mathbb{T}^d$.

Recall that $\{\chi_n\}_{n\geq -1}$ is the Littlewood-Paley partition functions and $\chi_n(\cdot)=\chi_0(2^{-n}\cdot)$ for $n\geq 0$ (see Appendix \ref{sec:appa}). Denote $\psi_n $ be the inverse Fourier transform of $\chi_n$ for every $n$, we then have the following result.
\begin{lemma}\label{lemma kol}For each $\kappa>0$ and every $n\geq -1$, it holds for all $s,t\in [0,T]$ and $x\in \mathbb{T}^d$ that
	\begin{align*}
		\mathbb{E}\left|\left\langle Z_t, \psi_n(\cdot-x)\right\rangle \right|^2&\lesssim 2^{{dn}},
		\\	\mathbb{E}\left|\left\langle Z_t, \psi_n(\cdot-x)\right\rangle-\left\langle Z_s, \psi_n(\cdot-x)\right\rangle \right|^2&\lesssim 2^{{dn}+2\kappa n}(t-s)^{\kappa},
	\end{align*}
where the proportional constants depend on $\|\bar{\rho}\|_{C_TL^{\infty}}$.
\end{lemma}
\begin{proof} For simplicity we set $\sigma=1$ in the proof.
	First, we use Fourier transform to represent $	\left\langle Z_t, \psi_n(\cdot-x)\right\rangle$ as follows,
	\begin{align}
			\left\langle Z_t, \psi_n(\cdot-x)\right\rangle &=\int_0^t\int_{\mathbb{T}^d}\int_{\mathbb{T}^d}\nabla\Gamma_{t-r}(y-z)\sqrt{\bar{\rho}_r(z)}\psi_n(y-x)\mathd y\xi(\dif r,\dif z)\nonumber
			\\&=\int_0^t\int_{\mathbb{T}^d}\left\langle  \nabla \Gamma_{t-r}(\cdot-z), \psi_n(\cdot-x)\right\rangle\sqrt{\bar{\rho}_r(z)}\xi(\dif r,\dif z)\nonumber
			\\&=\int_0^t\int_{\mathbb{T}^d}\sum_{k\in \mathbb{Z}^d}G_{t-r}(k) e_{-k}(z)\sqrt{\bar{\rho}_r(z)}\xi(\dif r,\dif z),\label{equation kol1}
	\end{align}
where $G_t(k)$  is defined by
\begin{align*}
	G_t(k):=\int_{\mathbb{T}^d}\left\langle  \nabla \Gamma_{t}(\cdot-z'), \psi_n(\cdot-x)\right\rangle e_{k}(z')\mathd z'
\end{align*}
and we used $\langle  \nabla \Gamma_{t-r}(\cdot-z), \psi_n(\cdot-x)\rangle\in L^2(\mT^d)$ and the sum in \eqref{equation kol1} converges in $L^2(\mT^d)$. 
Furthermore, noticing that $G_t(-k)$ is the  complex conjugate  of $G_t(k)$, we thus have
\begin{align}
		\mathbb{E}\left|\left\langle Z_t, \psi_n(\cdot-x)\right\rangle\right|^2=& \int_0^t\int_{\mathbb{T}^d}\left|\sum_{k \in \mathbb{Z}^d}G_{t-r}(k)e_{-k}(z) \right|^2 \bar{\rho}_r(z)\mathd z\mathd r\nonumber
		\\\lesssim&\|\bar{\rho}\|_{C_TL^{\infty}}\sum_{k_1\in \mathbb{Z}^d}\sum_{k_2 \in \mathbb{Z}^d}\int_0^t\int_{\mathbb{T}^d}G_{t-r}(k_1)G_{t-r}(-k_2)e_{-k_1}(z)e_{k_2}(z)\mathd z \mathd r\nonumber
			\\\lesssim&\|\bar{\rho}\|_{C_TL^{\infty}}\sum_{k\in \mathbb{Z}^d}\int_0^tG_{t-r}(k)G_{t-r}(-k) \mathd r,\label{equation kol2}
\end{align}
where  the last inequality follows by  $\int_{\mathbb{T}^d} e_{k_2-k_1}(z)\mathd z=C_d \delta_{k_2=k_1}$,  $C_d$ is the volume of $\mathbb{T}^d$. 

For each $k\in \mathbb{Z}^d$ and $t\in [0,T]$, we find that
\begin{align}
	G_{t}(k)=&\int_{\mathbb{T}^d}\int_{\mathbb{T}^d}\nabla\Gamma_{t}(y-z')\psi_n(y-x)e_{k}(z')\mathd y\mathd z'\nonumber
	\\=& \left\langle \nabla\Gamma_{t},e_{-k}\right\rangle\left\langle \psi_n,e_{k} \right\rangle e_{k}(x)\nonumber
	\\=&-\sqrt{-1}ke^{-t|k|^2}\chi_n(-k)e_{k}(x).\label{equation kol3}
\end{align}
Here we used   the facts that $\Gamma$ is the heat kernel on $\mathbb{T}^d$ and  $\left\langle \psi_n,e_{-k}\right\rangle=\chi_n(k) $.

Combining \eqref{equation kol2} and \eqref{equation kol3} yields that
\begin{align*}
	\mathbb{E}\left|\left\langle Z_t, \psi_n(\cdot-x)\right\rangle \right|^2\lesssim\|\bar{\rho}\|_{C_TL^{\infty}}\sum_{k \in \mathbb{Z}^d}\int_0^t|k|^2e^{-2|k|^2(t-r)}\chi_n(-k)\chi_n(k) \mathd r.
\end{align*}
Notice that
\begin{align*}
	\int_0^t|k|^2e^{-2|k|^2(t-r)}\mathd r=\frac{1}{2} \(1-e^{-2|k|^2t}\),
\end{align*}
which implies that for $n\geq 0$
\begin{align}
		\mathbb{E}\left|\left\langle Z_t, \psi_n(\cdot-x)\right\rangle \right|^2&\lesssim_{\bar{\rho}}\sum_{k \in \mathbb{Z}^d}\chi_n(-k)\chi_n(k)\(1-e^{-2|k|^2t}\)\nonumber
		\\&\lesssim_{\bar{\rho}}2^{dn}\int\chi_0(-k')\chi_0(k')\(1-e^{-2^{2n+1}|k'|^2t}\)\mathd k'\lesssim_{\bar{\rho}}2^{dn}.\label{equation kol4}
\end{align} 
Here we used $\chi_n(\cdot)=\chi_0(2^{-n}\cdot)$ and the fact that  $\chi_0$ is of compact support. The case  $n=-1$ is similar.

Next, we deduce  by \eqref{equation kol1} that
\begin{align*}
	&\left\langle Z_t, \psi_n(\cdot-x)\right\rangle-\left\langle Z_s, \psi_n(\cdot-x)\right\rangle
	\\=&\int_0^t\int_{\mathbb{T}^d}\sum_{k\in \mathbb{Z}^d}G_{t-r}(k) e_{-k}(z)\sqrt{\bar{\rho}_r(z)}\xi(\dif r,\dif z)-\int_0^s\int_{\mathbb{T}^d}\sum_{k\in \mathbb{Z}^d}G_{s-r}(k) e_{-k}(z)\sqrt{\bar{\rho}_r(z)}\xi(\dif r,\dif z)
	\\=&\int_s^t\int_{\mathbb{T}^d}\sum_{k\in \mathbb{Z}^d}G_{t-r}(k) e_{-k}(z)\sqrt{\bar{\rho}_r(z)}\xi(\dif r,\dif z)
	+\int_0^s\int_{\mathbb{T}^d}\sum_{k\in \mathbb{Z}^d}\[G_{t-r}(k)-G_{s-r}(k)\] e_{-k}(z)\sqrt{\bar{\rho}_r(z)}\xi(\dif r,\dif z)
	\\\assign&J_{1}^n+J_2^n.
\end{align*}
Moreover, we have
\begin{align}
	\mathbb{E}\left|\left\langle Z_t, \psi_n(\cdot-x)\right\rangle-\left\langle Z_s, \psi_n(\cdot-x)\right\rangle \right|^2\leq  2\mathbb{E}|J_1^n|^2+2\mathbb{E}|J_2^n|^2.\label{equation kol5}
\end{align}
Again, it suffices to check the cases with $n\geq 0$. Similar as in  $\eqref{equation kol4}$, we have
\begin{align*}
	\mathbb{E}|J_1^n|^2
	\lesssim_{\bar{\rho}}&\sum_{k \in \mathbb{Z}^d}\int_s^t|k|^2e^{-2|k|^2(t-r)}\chi_n(k)\chi_n(-k)\mathd r
	\\\lesssim_{\bar{\rho}}&\sum_{k \in \mathbb{Z}^d}\(1-e^{-2|k|^2(t-s)}\)\chi_n(k)\chi_n(-k)
	\\\lesssim_{\bar{\rho}}&2^{dn}\int\(1-e^{-2^{2n+1}|k'|^2(t-s)}\)\chi_0(k')\chi_0(-k')\mathd k'\
	\lesssim_{\bar{\rho}}2^{dn}(1-e^{-C2^{2n}(t-s)}),
\end{align*}
where $C>0$ is a universal constant determinded by the support of $\chi_0$. Notice that for each  $\kappa>0$, it holds  that $1-e^{-a}\lesssim a^{\kappa}$ for  $a\geq 0$. Therefore, for each $\kappa>0$, let $C2^{2n}(t-s)$ in the  above inequality play the  role of $a$, we arrive at
\begin{align}
	\mathbb{E}|J_1^n|^2\lesssim_{\bar{\rho}} 2^{dn+2\kappa n}(t-s)^{\kappa}.\label{equation kol6}
\end{align}
Similarly, one can study  $J_2^n$ and find
\begin{align*}
\mathbb{E}|J_2^n|^2\lesssim_{\bar{\rho}}&\sum_{k \in \mathbb{Z}^d}\int_0^s|k|^2\chi_n(k)\chi_n(-k)\(e^{-|k|^2(s-r)}-e^{-|k|^2(t-r)}\)^2\mathd r
\\\lesssim_{\bar{\rho}}&\sum_{k \in \mathbb{Z}^d}\(1-e^{-2|k|^2s}\)\chi_n(k)\chi_n(-k)\(1-e^{-|k|^2(t-s)}\)^2
\\\lesssim_{\bar{\rho}}&2^{dn}\int_{k' \in \mathbb{R}^d}\chi_0(k')\chi_0(-k')\(1-e^{-2^{2n}|k'|^2(t-s)}\)^2\mathd k'
\lesssim_{\bar{\rho}}2^{dn+2\kappa n}(t-s)^{\kappa},
\end{align*}
for each $\kappa >0$. This together with \eqref{equation kol5} and \eqref{equation kol6} leads to
\begin{equation*}
		\mathbb{E}\left|\left\langle Z_t, \psi_n(\cdot-x)\right\rangle-\left\langle Z_s, \psi_n(\cdot-x)\right\rangle \right|^2\lesssim_{\bar{\rho}}2^{dn+2\kappa n}(t-s)^{\kappa}.
\end{equation*}
The proof is thus completed.	\end{proof}

    We now apply the above result to study regularity of the process $Z$.
\begin{lemma}\label{lemma z} Suppose that $\bar{\rho}\in C([0,T],L^{\infty})$. It holds that $Z\in C([0,T],C^{-\alpha})$ $\mP$-a.s. for every $\alpha>d/2$. Moreover,  for all $p>2$,
	\begin{equation*}
	\mathbb{E}\sup_{t \in [0, T]}\|Z_t\|^p_{C^{-\alpha}}<\infty.
	\end{equation*}
\end{lemma}

\begin{proof}
	Since $Z$ is a centered Gaussian process,
	 Lemma \ref{lemma kol} together with  the hypercontractivity property \cite[Theorem 1.4.1]{nualart2006malliavin} implies that
	\begin{align*}
	\mathbb{E}\left|\left\langle Z_t, \psi_n(\cdot-x)\right\rangle \right|^{p}&\lesssim \(\mathbb{E}\left|\left\langle Z_t, \psi_n(\cdot-x)\right\rangle \right|^{2}\)^{\frac{p}{2}}\lesssim2^{\frac{dnp}{2}},
	\\	\mathbb{E}\left|\left\langle Z_t, \psi_n(\cdot-x)\right\rangle-\left\langle Z_s, \psi_n(\cdot-x)\right\rangle \right|^p&\lesssim\(\mathbb{E}\left|\left\langle Z_t, \psi_n(\cdot-x)\right\rangle-\left\langle Z_s, \psi_n(\cdot-x)\right\rangle \right|^2\)^{\frac{p}{2}}\lesssim 2^{(\frac{d}{2}+\kappa) np}(t-s)^{\frac{\kappa p}{2}},
\end{align*}
for each $\kappa>0$, $p>2$, and every $n\geq -1$. This allows us to apply the Kolmogorov criterion  \cite[Lemma 10]{mourrat2017global} to conclude that $Z\in C([0,T], B^{-\alpha}_{p,p})$ $\mP$-a.s., for each $p>2$ and every $\alpha>d/2+2/p$.  Moreover,
\begin{equation*}
	\mathbb{E}\sup_{t \in [0, T]}\|Z_t\|_{B^{-\alpha}_{p,p}}^p<\infty.
\end{equation*}
The result  follows by  the embedding  $B^{-\alpha}_{p,p}\hookrightarrow B^{-\beta}_{\infty,\infty}$ for $\beta>\alpha+d/p$ (see Lemma \ref{lemma embedding}). 
	\end{proof}

Next we  rewrite the unique solution $\eta$ to \eqref{equation li}, which has been obtained in Section \ref{sec:chara},   in the  mild form.
\begin{proposition} \label{lemma mild}Under the assumptions {\bf{(A1)-(A4)}}, the unique solution $\eta $ to \eqref{equation li} satisfies
	\begin{equation*}
		\eta_t=\Gamma_t*\eta_0-\int_0^t\nabla\Gamma_{t-s}*(\bar{\rho} K \ast \eta+\eta K \ast \bar{\rho}+F \eta)\mathd s-\sqrt{2\sigma}\tilde{Z}_t,\quad \mathbb{P}-a.s,
	\end{equation*}
where $\tilde{Z}$ has the same distribution as  $Z$. 
\end{proposition}
\begin{proof}
	We start with the following statement: for every function $\varphi$ of class
	$C^1 ([0, t], C^{\infty}(\mathbb{T}^d))$ and $t \in [0, T]$, it  holds that,
\begin{align}
	\left\langle \eta_{t},\varphi(t)\right\rangle -\left\langle \eta_{0},\varphi(0)\right\rangle =&\int_0^t  \left\langle\eta_s, \partial_s\varphi+\sigma \Delta \varphi \right\rangle    \mathd s+\int^t_0 \langle
	\bar{\rho}_s K \ast \eta_s + \eta_s K \ast \bar{\rho}_s + F \nocomma \eta_s, \nabla
	\varphi \rangle \mathd s\nonumber
	\\&+\sqrt{2\sigma}\int_0^t\int_{\mathbb{T}^d}\nabla\varphi(x)\sqrt{\bar{\rho}_s(x)}\xi(\dif s,\dif x).\label{eqt:weakform}
\end{align}
It is	straightforward to check the statement for finite linear combinations of
	functions $\varphi$ of the form $\varphi (s, x) = \varphi_1 (s) \varphi_2
	(x)$, where $\varphi_1 \in C^{\infty} ([0, t])$ and $\varphi_2 \in
	C^{\infty} (\mathbb{T}^d)$. Then one can uniformly approximate functions in
	$C^1 ([0, t], C^{\infty} (\mathbb{T}^d))$ with such combinations and find \eqref{eqt:weakform}.
	
	For every  $\varphi_0\in C^{\infty}(\mT^d)$ and $0\leq s\leq t$, define $\varphi(s):=\Gamma_{t-s}*\varphi_0$,  then $\partial_s\varphi(s)= -\sigma\Delta\varphi(s)$ and $\varphi(t)=\varphi_0$. By \eqref{eqt:weakform}, we find
	\begin{align*}
			\left\langle \eta_{t},\varphi_0\right\rangle -\left\langle \eta_{0},\Gamma_{t}*\varphi_{0}\right\rangle =&\int^t_0 \langle
			\bar{\rho}_s K \ast \eta_s + \eta_s K \ast \bar{\rho}_s + F \nocomma \eta_s, \nabla \Gamma_{t-s}*
			\varphi_0 \rangle \mathd s
			\\&+\sqrt{2\sigma}\int_0^t\int_{\mathbb{T}^d}(\nabla\Gamma_{t-s}*
			\varphi_0)(x)\sqrt{\bar{\rho}_s(z)}\xi(\dif s,\dif x)
			\\=&
			-\int_0^t\left\langle \nabla\Gamma_{t-s}*(\bar{\rho} K \ast \eta+\eta K \ast \bar{\rho}+F \eta),\varphi_0\right\rangle \mathd s
			\\&-\sqrt{2\sigma}\int_{\mathbb{T}^d}\varphi_{0}(x)\left( \int_0^t\int_{\mathbb{T}^d} \nabla\Gamma_{t-s}(x-y)\sqrt{\bar{\rho}(y)}\xi(\mathd s,\mathd y)\right) \mathd x,
	\end{align*}
where we used symmetry of $\Gamma$ at the last inequality. The result then follows by arbitrary $\varphi_0\in C^{\infty }(\mathbb{T}^d)$ and the definition of $Z$.
\end{proof}

This result gives rise to the definition of mild solutions to \eqref{equation li} on a stochastic basis $(\Omega, \mathcal{F}_t,  \mathcal{F},\mP)$. We set $Z$ given by \eqref{eq:Z} with $\xi$ being vector-valued space-time white noise on $(\Omega, \mathcal{F}_t,  \mathcal{F},\mP)$.

\begin{definition}\label{def:mild} Assume  that $K\in L^1$, $\bar{\rho}\in C([0,T],C^{\beta})$, and $F\in C^{\beta}$ for some $\beta>d/2$.
We call $\eta\in C([0,T],\mathcal{S}'(\mathbb{T}^d))\cap L^2([0,T], B^{-\alpha}_{p,q})$ with $\alpha<\beta$, $p,q\in [1,\infty]$ is a mild solution to \eqref{equation li} with initial condition $\eta_0$  if for every $\varphi\in C^{\infty}$
	\begin{align*}
	\left\langle \eta_{t},\varphi\right\rangle =\left\langle \Gamma_{t}*\eta_{0},\varphi\right\rangle
	-\int_0^t\left\langle \nabla\Gamma_{t-s}*(\bar{\rho} K \ast \eta+\eta K \ast \bar{\rho}+F \eta),\varphi\right\rangle \mathd s -\sqrt{2\sigma}\left\langle Z_t,\varphi\right\rangle  .
\end{align*}
\end{definition}

\br By Proposition \ref{lemma mild}, we know, under the assumptions {\bf{(A1)-(A4)}},  the solutions  to \eqref{equation li} obtained from Theorem \ref{theorem existence}  have the same law as the mild solutions.

To make sense of $\bar{\rho}K*\eta$, $\eta K*\bar{\rho}$, and $F\eta$ in the definition of mild solutions,  we need the condition $\eta_t\in B^{-\alpha}_{p,q}$ for a.e. $t\in [0,T]$  with $\alpha<\beta$ and $p,q\in [1,\infty]$, where $\beta$ is from {\bf{(A4)}}.
 \er

 The following result based on the smoothing effect of heat kernel (see Lemma \ref{lemma schau}) gives the optimal regularity of $\eta$.
\begin{proposition}\label{pro:op} Suppose that Assumption {\bf{(A4)}} holds with parameter $\beta>d/2$
	and $\eta$ is a mild solution to \eqref{equation li},  and assume $\eta_{0}\in L^r(\Omega, B^{-\alpha}_{p,q})$ for some  $r>2$, $\alpha\in ( d/2,\beta)$, and $p,q\in [1,\infty]$. Then  $\eta\in C([0,T],B^{-\alpha}_{p,q} )$  almost surely. Moreover,
	\begin{equation*}
		\mathbb{E}\sup_{ t \in [0, T]}\|\eta_t\|^r_{B^{-\alpha}_{p,q}}<\infty.
	\end{equation*}
\end{proposition}
\begin{proof}  
	Firstly,
applying Lemma \ref{lemma schau}, we have
\begin{align*}
	\|\eta_{t}\|_{B^{-\alpha}_{p,q}}\lesssim& \|\eta_{0}\|_{B^{-\alpha}_{p,q}}+\int_0^t(t-s)^{-\frac{1}{2}}\[\|\bar{\rho}K*\eta\|_{B_{p,q}^{-\alpha}}+\|\eta K*\bar{\rho}\|_{B_{p,q}^{-\alpha}}+\|F\eta \|_{B_{p,q}^{-\alpha}}\]\mathd s
	\\&+\|{Z}_t\|_{B_{p,q}^{-\alpha}}.
\end{align*}
To further estimate the right hand side of the above inequality, by $\alpha<\beta$,   applying Lemmas \ref{lemma embedding}-\ref{lemma convolution} gives that
\begin{align*}
	\|\bar{\rho}K*\eta\|_{B_{p,q}^{-\alpha}}+	\|\eta K*\bar{\rho}\|_{B_{p,q}^{-\alpha}}\lesssim\|\bar{\rho}\|_{C^{\beta}}\|K\|_{L^1}\|\eta\|_{B^{-\alpha}_{p,q}}, \quad	\|F\eta\|_{B^{-\alpha}_{p,q}}\lesssim\|F\|_{C^{\beta}}\|\eta\|_{B^{-\alpha}_{p,q}}.
		\end{align*}
Hence,
\begin{align}
		\|\eta_{t}\|_{B^{-\alpha}_{p,q}}\lesssim& \|\eta_{0}\|_{B^{-\alpha}_{p,q}}+\int_0^t(t-s)^{-\frac{1}{2}}\|\eta_s\|_{B^{-\alpha}_{p,q}}\mathd s+\|{Z}_t\|_{B_{p,q}^{-\alpha}}.\label{equation op1}
\end{align}
By H\"older's inequality, we find
\begin{align}
	\int_0^t(t-s)^{-\frac{1}{2}}\|\eta_s\|_{B^{-\alpha}_{p,q}}\mathd s&\lesssim \(\int_0^t\|\eta_s\|^r_{B^{-\alpha}_{p,q}}\mathd s\)^{\frac{1}{r}}\(\int_0^t(t-s)^{-\frac{r}{2(r-1)}}\mathd s\)^{\frac{r-1}{r}}\nonumber
	\\&\lesssim \(\int_0^t\|\eta_s\|^r_{B^{-\alpha}_{p,q}}\mathd s\)^{\frac{1}{r}}t^{\frac{1}{2}-\frac{1}{r}},\label{equation op2}
\end{align}
where in the last step we used $r>2$ to have $\frac{r}{2(r-1)}<1$.
Furthermore, applying Gronwall's inequality  to \eqref{equation op1} yields that
	\begin{equation*}
	\mathbb{E}\sup_{ t \in [0, T]}\|\eta_t\|_{B^{-\alpha}_{p,q}}^r\lesssim \mathbb{E} \|\eta_{0}\|^r_{B^{-\alpha}_{p,q}}+\mathbb{E}\|{Z}\|_{C_TB_{p,q}^{-\alpha}}^r+1.
\end{equation*}
By the assumption on $\eta_{0}$ and Lemma \ref{lemma z},  the right hand side of the above inequality is thus finite.

 Using  \eqref{equation op1}, the continuity of $\eta $ on $[0,T]$  follows by the continuity of  ${Z}$ and continuity of $\Gamma_t$ from Lemma \ref{lemma schau}.  The proof is thus completed.
	\end{proof}
\begin{remark}  
	 By \cite{hairer2014theory} we know  the space-time white noise $\xi\in C_{t,x}^{-\frac{d}2-1-\eps}$  $\mP$-a.s. for every $\eps>0$, where $ C_{t,x}^{-\frac{d}2-1-\eps}$ is endowed with suitable parabolic time space scaling. Hence by Schauder estimates  $\eta\in C([0,T],C^{-\alpha})$ for $\alpha>d/2$ gives the best regularity by taking $p,q=\infty$ in Proposition \ref{pro:op}.
\end{remark}
\begin{remark}\label{re:00}
	By the optimal regularity of $\eta$  and Lemma \ref{lemma triebel},  $K*\bar{\rho}$ needs to stay in $C^{\beta}$ with $\beta >d/2$ so that  the term $K*\bar{\rho}\eta$ appearing in  SPDE \eqref{equation li} is well-defined.  Applying Lemma \ref{lemma convolution} and noticing $K\in L^1$ for $K$ satisfying {\bf{(A2)}},  the  assumption about $\bar{\rho}$ in {\bf{(A4)}} is thus a sufficient condition for $K*\bar{\rho}\in C^{\beta}$.  Moreover, $\beta>d/2$ is optimal in general.  With  appropriate modifications of the proof in  Section \ref{sec:chara} and by the  convolution inequality in Lemma \ref{lemma convolution},   we could also weaken the  condition of $\bar\rho$ to $\bar{\rho}\in C^{\beta-\beta_1}$  with $\beta>d/2$ and $\beta_1\in (0,\frac{d}{2})$, at  the cost of stronger condition on the interacting kernel: $K\in C^{\beta_1}$.
\end{remark}

\subsection{Gaussianity}\label{sec:gauss}
This section is devoted to the proof of Proposition \ref{prop:gauss}. As mentioned in the introduction, we  need a class of  time evolution operators $\{Q_{s,t}\}$ in order to rewrite $\eta$ as the  generalized Ornstein-Uhlenbeck process \eqref{eqt:ou},  which would be given by the following result.
\begin{lemma} \label{lemma:back}Assume that $\bar{\rho}\in C([0,T],C^{\beta+1}(\mathbb{T}^d))$ and $F\in C^{\beta+1} (\mathbb{T}^d)$ with $\beta >d/2$, for each $\varphi\in C^{\infty}(\mathbb{T}^d)$ and $t\in [0,T]$, there exists a unique solution $f\in  L^2([0,t], H^{\beta+2})\cap C([0,t],H^{\beta+1})$ with $\partial_sf\in  L^2([0,t],H^{\beta})$
	to the following backward equation:
	\begin{equation}
	 f_s =\varphi +\sigma \int_s^t\Delta f_r\mathd r+\int_s^t\[K*\bar{\rho}_r\cdot\nabla f_r+K(-\cdot)*(\nabla f_r\bar{\rho}_r) +F\cdot\nabla f_r \] \mathd r,\quad s\in [0,t],\label{eqt:back}
	\end{equation}
where $K(-\cdot)*g$ is given in \eqref{notation:k}.
\end{lemma}
\begin{proof}Similar to \eqref{eqt:energy}, we obtain the following a priori  energy estimate for any $\epsilon>0$
	\begin{align*}
		 &\|f_s \|^2_{H^{\beta+1}} + 2 \sigma   \int_s^t \|f_r
		\|^2_{H^{\beta+2}} \mathd r
		\\\leqslant & \|\varphi\|_{H^{\beta+1}}^2+C_{\epsilon}  \int_s^t \|f_r \|^2_{H^{\beta+1}} \mathd r
		+ \epsilon ( {\nobracket
			\|K\|^2_{L^1} \| \bar{\rho} \|^2_{C_TC^{\beta+1}}
			+\|F\|^2_{C^{\beta+1}}}) \int_s^t  \|f_r \|^2_{H^{\beta+2}} \mathd r.
	\end{align*}
Choosing $\epsilon>0$ sufficiently small,  $f\in  L^2([0,t], H^{\beta+2})$ follows from the Gronwall's inequality.
 Furthermore, by Lemma \ref{lemma triebel} and Lemma \ref{lemma convolution}, we find
\begin{align*}
	&\|K*\bar{\rho}\cdot\nabla f\|_{H^{\beta}}\lesssim \|K*\bar{\rho}\|_{C^{\beta}}\|f\|_{H^{\beta+1}}\lesssim\|K\|_{L^1}\|\bar{\rho}\|_{C^{\beta}}\|f\|_{H^{\beta+1}},\quad 	\|F\cdot\nabla f\|_{H^{\beta}}\lesssim\|F\|_{C^{\beta}}\|f\|_{H^{\beta +1}},
	\\&	\|K(-\cdot)*(\bar{\rho}\nabla f)\|_{H^{\beta}}\lesssim \|K\|_{L^1}\|\bar{\rho}\nabla f\|_{H^{\beta}}\lesssim\|K\|_{L^1}\|\bar{\rho}\|_{C^{\beta}}\|f\|_{H^{\beta+1}}.
\end{align*}
 Hence we deduce from equation \eqref{eqt:back} that $\partial_tf\in  L^2([0,t],H^{\beta})$, which combined with $f\in  L^2([0,t], H^{\beta+2})$ implies that $f\in C([0,t],H^{\beta+1})$ by Lions-Magenes Lemma. When $\varphi=0$, the above energy estimate implies that $f=0$. This fact together with linearity of equation implies the uniqueness of solutions. On the other hand, one can obtain the existence of solutions to \eqref{eqt:back} by classical Galerkin method (cf. \cite[Chapter 7]{evans1998partial}).
\end{proof}
Define the space $\mathcal{X}_t^{\beta}$ and time  evolution operators $\{Q_{\cdot,t}\}_{0\leq t\leq T}:C^{\infty}(\mathbb{T}^d)\rightarrow \mathcal{X}_t^{\beta}$ as
\begin{gather*}
	\mathcal{X}_t^{\beta}\assign\left\lbrace f\in  L^2([0,t], H^{\beta+2})\cap C([0,t],H^{\beta+1}); \partial_sf\in  L^2([0,t],H^{\beta}) \right\rbrace,
	\\Q_{s,t}\varphi\assign f(s).
\end{gather*}
where $f$ is the unique solution to \eqref{eqt:back} with terminal value $\varphi$ at time $t$ and is given by Lemma \ref{lemma:back}.

Now we are in the position to justify the Gaussianity of the unique (in distribution) limit of fluctuation measures.
\begin{proof}[Proof of Proposition \ref{prop:gauss}] Recall that  $\eta\in L^2([0,T],H^{-\alpha})$ for any $\alpha>d/2$. Then for each test function $f\in \mathcal{X}_t^{\beta}$ with $\beta >d/2$, by Lemma \ref{lemma triebel} and Lemma \ref{lemma convolution}, we have
	\begin{align*}
	&	\|\langle \eta_s,\partial_sf\rangle \|_{L_T^1}\lesssim\|\eta\|_{L_T^2H^{-\beta}}\|\partial_sf\|_{L_T^2H^{\beta }},  \quad
		\|\langle \eta_s,\Delta f_s\rangle \|_{L_T^1}\lesssim\|\eta\|_{L_T^2H^{-\beta}}\|f\|_{L_T^2H^{\beta+2 }},
		\\&\| \langle	\bar{\rho}_s K \ast \eta_s + \eta_s K \ast \bar{\rho}_s + F \nocomma \eta_s, \nabla f_s\rangle \|_{L_T^1}\lesssim \(\|K\|_{L^1}\|\bar{\rho}\|_{C_TC^{\beta}}+\|F\|_{C^{\beta}}\)\|\eta\|_{L_T^2H^{-\beta+\epsilon}}\|f\|_{L_T^2H^{\beta+1}}\lesssim1,
	\end{align*}
where $\epsilon>0$ is sufficiently small, so that the weak formulation \eqref{eqt:weakform}  extends to all the test fucntions  $f\in \mathcal{X}_t^{\beta}$ with $\beta >d/2$. 	For each $\varphi\in C^{\infty}$, choosing $Q_{\cdot,t}\varphi$ as the test function in \eqref{eqt:weakform}, we find
	\begin{align*}
		\left\langle \eta_{t},\varphi\right\rangle -\left\langle \eta_{0},Q_{0,t}\varphi\right\rangle =&\int_0^t  \left\langle\eta_s, \partial_sQ_{s,t}\varphi+\sigma \Delta Q_{s,t}\varphi \right\rangle    \mathd s+\int^t_0 \langle
		\bar{\rho}_s K \ast \eta_s + \eta_s K \ast \bar{\rho}_s + F \nocomma \eta_s, \nabla
		Q_{s,t}\varphi\rangle \mathd s
		\\&+\sqrt{2\sigma}\int_0^t\int_{\mathbb{T}^d}\nabla Q_{s,t}\varphi(x)\sqrt{\bar{\rho}_s(x)}\xi(\dif s,\dif x)
		\\=&\sqrt{2\sigma}\int_0^t\int_{\mathbb{T}^d}\nabla Q_{s,t}\varphi(x)\sqrt{\bar{\rho}_s(x)}\xi(\dif s,\dif x),
	\end{align*}
	where we used Lemma \ref{lemma:back} with $Q_{s,t}\varphi=f(s)$.   The result follows by the assumption on $\eta_0$ and the fact that the stochastic integral is a centered  Gaussian process with quadratic variation
	\begin{equation*}
	2\sigma\int_0^t \left\langle |\nabla Q_{s,t}\varphi|^2,\bar{\rho}_s \right\rangle \mathd s.
	\end{equation*}
\end{proof}

\section{The Vanishing Diffusion Case}\label{sec:vani}

In this section, we study particle systems with vanishing diffusion, i.e. $\sigma=0$. 
We  denote  the  fluctuation measures by $\eta^N
\assign \sqrt{N} (\mu^N- \bar{\rho})$ as well. Instead of  the SPDE limit \eqref{equation li} in the case when $\sigma >  0$,   now  the fluctuation measures converge to a deterministic  first order nonlocal PDE with random initial value $\eta_0$, which reads
\begin{equation}
	\partial_t \eta = - \nabla \cdot (\bar{\rho} K \ast \eta) - \nabla \cdot
	(\eta K \ast \bar{\rho}) - \nabla \cdot (F \eta) . \label{equation
		degenerate pde}
\end{equation}

With  the  same proof as in Section \ref{sec:tight} and Section \ref{sec:chara}, we can  deduce that under  the assumptions {\bf{(A1)-(A3), (A5)}}, we have
\begin{enumerate}
	\item The sequence of laws of $(\eta^N)_{N \geqslant 1}$ is tight in the space
	$L^2 ([0, T], H^{- \alpha}) \cap C ([0, T], H^{- \alpha - 2})$, for every
	$\alpha > d / 2$.
	
	\item Any limit $\eta$ of converging (in distribution) subsequence of $(\eta^N)_{N
		\geqslant 1}$ is an analytic weak solution to {\eqref{equation degenerate
			pde}} in the sense that
	\begin{align}
	\left\langle \eta_{t},\varphi(t) \right\rangle  =\left\langle \eta_0,\varphi(0)\right\rangle + \int^t_0
	\int_{\mathbb{T}^d} {\eta_s} [\partial_s \varphi + K(-\cdot) \ast
	(\bar{\rho} \nabla \varphi) + K \ast \bar{\rho} \cdummy \nabla \varphi + F
	\nabla \varphi] \mathd x \nocomma \mathd s, \label{equation test}
	\end{align}	
for every $\varphi \in C^{1} ([0,t],C^{\beta+1})$ with $\beta>d/2$, $\mathbb{P}$-a.s. and $K(-\cdot)*g$ is given in \eqref{notation:k}.
\end{enumerate}

However, for the case with vanishing diffusion, we cannot deduce uniqueness of the solutions to \eqref{equation degenerate pde} by the proof of  Lemma \ref{lemma unique pde}, due to the lack of  the energy inequality  \eqref{eqt:energy}. The following uniqueness result follows by the  method of characteristics. We also recall the definition of flow from \cite[Chapter 4]{kunita1997stochastic}\footnote{Although the framework in \cite{kunita1997stochastic} is for the flow on $\mR^d$, it also holds for the periodic case since the  functions on $\mT^d$ could be viewed as periodic functions on $\mR^d$ and the framework in \cite[Chapter 4]{kunita1997stochastic} has been extended  to Riemannian manifold. }, which is used in the following proof. $\phi_{s,t}$ be a continuous map from $\mT^d$ into itself for any $s,t\in [0,T]$ is called a flow if it satisfies the following property
\begin{enumerate}
	\item $\phi_{s,u}=\phi_{t,u}\circ \phi_{s,t}$ holds for all $s,t, u$, where $\circ$ denotes the composition of maps;
	\item $\phi_{s,s}=\mathrm{Id}$;
	\item $\phi_{s,t}:\mT^d\to \mT^d$ is an onto homeomorphism for all $s,t$.
\end{enumerate}
We refer to \cite[Chapter 4]{kunita1997stochastic} for the relation between  flows and ODEs. In general the solutions to ODEs with regular coefficients could generate a flow.

\begin{proposition}\label{lemma unique degenerate}
	Under the assumptions  {\bf{(A2)}} and {\bf{(A5)}},  $\eta
	\equiv 0$ is the only solution with zero initial value to {\eqref{equation degenerate pde}} in the space 	$L^2 ([0, T], H^{- \alpha}) \cap C ([0, T], H^{- \alpha - 2})$ for $\alpha\in (d/2,\beta)$, the parameter $\beta$ is from {\bf{(A5)}}.
\end{proposition}

\begin{proof}
	We first claim that a similar result to Lemma \ref{lemma:back} holds for $\sigma=0$. That is, there exists a unique solution $\varphi \in C^1 ([0, t],
	C^{\beta'})$ with $\beta'\in (\alpha+1,\beta +1)$ for any $\varphi (t, x) = \psi (x) \in C^{\infty}$ to the
	following backward equation
	\begin{equation}
		\partial_s \varphi + K(-\cdot) \ast (\bar{\rho} \nabla \varphi) + K \ast
		\bar{\rho} \cdummy \nabla \varphi + F \cdummy \nabla \varphi = 0, \quad \forall s
		\in [0, t] . \label{equation back}
	\end{equation}
	Suppose that the claim holds. Then for every $\psi \in C^{\infty}$ and $t \in [0, T]$, we use
	{\eqref{equation test}} with the test function given by the solution $\varphi$ to \eqref{equation back} and $\eta_{0}=0$. Then we conclude that $\int_{\mathbb{T}^d} {\eta} (t, x) \psi (x)
	\mathd x = 0$, which implies the result. It thus suffices to justify the claim.
	
	In the following we verify the claim by
	considering the backward flow $(\phi_{t,s} )_{0\leq s\leq t\leq T}$ 
	generated by
	\begin{equation}\label{eq:phi}
		\phi_{t,s} = x +\int^t_s (K \ast \bar{\rho}_r (\phi_{t,r} ) + F (\phi_{t,r}
		)) \mathd r, \quad x \in \mathbb{T}^d , s\in [0,t].
	\end{equation}
	Define the forward flow $\phi_{s,t}\assign\phi_{t,s}^{-1}$, $0\leq s\leq t\leq T$. Since $F \in C^{\beta+1}$ and $K \ast \bar{\rho} \in
	C^1([0,T],C^{\beta +2})$ by Assumption {\bf{(A5)}} and Lemma \ref{lemma convolution}, the existence of the flow  $(\phi_{t,s} )_{t,s\in [0,T]}$  in $C^{\beta'}$ for any $\beta'<\beta+1$ follows from \cite[Theorem 4.6.5]{kunita1997stochastic}.  For fixed $t\in [0,T]$, denote $\phi_s\assign \phi_{t,s}$ and $\phi_{-s}\assign\phi_{s,t}$ for $s\in [0,t]$, then $\phi_{- s}\circ \phi_s=\tmop{Id}$.
	
	The next step is  the  one-to-one correspondence between the solutions to {\eqref{equation back}}  and the solutions to the following equation
		\begin{align}
		g_s(x)= &\psi(x)-\int_s^t \[K(-\cdot)*(g_r\circ \phi_{-r}\nabla\bar{\rho})\]\circ \phi_r(x) \mathd r -\int_s^t\[\div K(-\cdot)*(g_r\circ \phi_{-r}\bar{\rho})\]\circ \phi_r(x)\mathd r\nonumber
		\\\assign&\psi(x)+\Phi_s(g), \quad s\in [0,t]\label{eqt:ode}.
	\end{align}
	Here the notation $K(-\cdot)*f$ is given in \eqref{notation:k}. We also write $g(s,x)\assign g_s(x)$. We will prove that $g(s,x):=\varphi(s,\phi_s)$ satisfies \eqref{eqt:ode}. Indeed, suppose that $\varphi\in C^1 ([0, t], C^{\beta'})$ with $\beta'>\alpha+1$ solves  {\eqref{equation back}}, then by the chain rule, we have
	\begin{align}
		\partial_s  \{\varphi (s, \phi_s (x))\} = & \partial_s \varphi (s, \phi_s (x)) + \nabla
		\varphi (s, \phi_s (x)) \cdummy \partial_s \phi_s (x)=- {K}(-\cdot) \ast (\bar{\rho} \nabla \varphi) (\phi_s (x)) \nonumber\\
		= & {K}(-\cdot) \ast (\varphi \nabla \bar{\rho} ) (\phi_s (x)) + \tmop{div} \nocomma
		{K}(-\cdot) \ast (\varphi \bar{\rho} ) (\phi_s (x)) \nonumber\\
		= &\[K(-\cdot)*\{(\varphi \circ \phi_{s}\circ \phi_{-s})\nabla\bar{\rho}\}\]\circ \phi_s(x)+
		\[\div K(-\cdot)*\{(\varphi\circ\phi_s\circ \phi_{-s})\bar{\rho}\}\]\circ \phi_s(x),
		\nonumber
	\end{align}
	where we used integration by parts formula in the third step.
Therefore $\varphi(s,\phi_s)$ satisfies \eqref{eqt:ode}.

 Conversely, if $g \in C^1 ([0, t],
	C^{\beta'})$ is a solution to equation \eqref{eqt:ode}, let $\varphi(s,x)\assign g (s, \phi_{- s} (x))$, then $g(s,x)=\varphi(s,\phi_s(x))$. Similarly, we have
	\begin{align*}
		\partial_s g(s,x)= &\partial_s  \{\varphi (s, \phi_s (x))\} =  \partial_s \varphi (s, \phi_s (x)) + \nabla
		\varphi (s, \phi_s (x)) \cdummy \partial_s \phi_s (x),
		\\	\partial_s g(s,x)= &{K}(-\cdot) \ast (\varphi \nabla \bar{\rho} ) (\phi_s (x)) + \tmop{div} \nocomma
		{K}(-\cdot) \ast (\varphi \bar{\rho} ) (\phi_s (x)) =- {K}(-\cdot) \ast (\bar{\rho} \nabla \varphi) (\phi_s (x)),
	\end{align*}
	where the first line follows by the chain rule, while the second line follows by integration by parts. Substituting \eqref{eq:phi} into the first line, we obtain that  $ \varphi$ is a solution to
	{\eqref{equation back}}. Hence justifying the claim is turned into obtaining  the well-posedness of backward equation \eqref{eqt:ode} in $ C^1([0,t],C^{\beta'})$.
	
	Notice that  $\Phi_{\cdot}:C([s,t],C^{\beta'})\rightarrow C([s,t],C^{\beta'})$ satisfies
		\begin{align}
	\sup_{ r\in [s,t]}	\|\Phi_r(g)\|_{C^{\beta'}}\leq &\int_s^t\left\| \[K(-\cdot)*(g_r\circ \phi_{-r}\nabla\bar{\rho})\]\circ \phi_r\right\|_{C^{\beta '}} +\left\| \[\div K(-\cdot)*(g_r\circ \phi_{-r}\bar{\rho})\]\circ \phi_r\right\|_{C^{\beta '}} \mathd r\nonumber
		\\\lesssim &\int_s^t\|K\|_{L^1}\|g_r\|_{C^{\beta '}}\(1+\|\phi_{-r}\|_{C^{\beta '}}^{\beta'}\)\|\bar{\rho}\|_{C^{\beta '+1}}(1+\|\phi_r\|_{C^{\beta'}}^{\beta'})\mathd r\nonumber
\\&+\int_s^t\|\div K\|_{L^1}\|g_r\|_{C^{\beta '}}\(1+\|\phi_{-r}\|_{C^{\beta '}}^{\beta'}\)\|\bar{\rho}\|_{C^{\beta '}}\(1+\|\phi_r\|_{C^{\beta'}}^{\beta'}\)\mathd r,\label{xxxxxx}
	\end{align}
where we used Lemma \ref{lemma triebel}, Lemma \ref{lemma convolution}, and  the fact that $\|f_1\circ f_2\|_{C^{\alpha}}\lesssim\|f_1\|_{C^{\alpha}}(1+\| f_2\|_{C^{\alpha}}^{\alpha})$ when $\alpha\geq 1$.
	 Recalling Assumption {\bf{(A5)}}	and $\phi\in C([0,t], C^{\beta'})$,  we find
	 \begin{align*}
	 	\sup_{ r\in [s,t]}\|\Phi_r(g)\|_{C^{\beta'}}	\lesssim(t-s)\sup_{ r\in [s,t]}\|g_r\|_{C^{\beta'}}.
	 \end{align*}
 Choosing $s$ close to $t$ enough, by the linearity of $\Phi$, we find $g\mapsto\psi+\Phi_{\cdot}(g)$ is a contraction mapping on  $C([s,t],C^{\beta'})$, hence it has a unique fixed point solving  {\eqref{eqt:ode}} on $[s,t]$. Applying this argument a finite number of times, we remove the constraint on $s$. Therefore, there exists a unique solution   $g\in C([0,t],C^{\beta'})$ to the ODE \eqref{eqt:ode}. Finally, we deduce  $g\in C^1([0,t],C^{\beta'})$ from $\partial_s\Phi_s\in C([0,t],C^{\beta'})$, which follows by  the  calculation in \eqref{xxxxxx}.
 Now we obtain the global well-posedness of \eqref{eqt:ode} in $ C^1([0,t],C^{\beta'})$, which by the one-to-one correspondence between $\varphi$ and $g$ concludes the result.
\end{proof}

Now we are able  to conclude the result for vanishing diffusion
cases similar to Theorem \ref{thm:1}.

\begin{theorem}\label{thm:de}
	Under the assumptions {\bf{(A1)-(A3), (A5)}}, let $\eta$ be the unique solution to  {\eqref{equation
			degenerate pde}} on the same  stochastic basis with the particle system \eqref{equation pa}, the sequence $(\eta^N)_{N \geqslant
		1}$ converges in probability to $\eta$ in  $L^2 ([0, T], H^{- \alpha}) \cap C ([0, T],
	H^{- \alpha - 2})$, for every $\alpha > d / 2$.
\end{theorem}

\begin{proof}
By the facts that  $(\eta^N)_{N \geqslant1}$ are tight, the tight limits of converging subsequences solve {\eqref{equation
		degenerate pde}},  and  there exists a unique analytic weak solution to the   equation {\eqref{equation
		degenerate pde}}, which is ensured by Proposition  \ref{lemma unique degenerate}, it follows immediately that the sequence $(\eta^N)_{N \geqslant1}$  converges in distribution to the unique solution $\eta$ .
	
	Similar to $(\eta^N)_{N \geqslant1}$, one can first obtain tightness of laws of $(\eta^{l},\eta^{m})_{l,m\in \mathbb{N}}$. Without loss of generality, we assume $(\eta^{l},\eta^{m})_{l,m\in \mathbb{N}}$ to be two converging subsequences. Then using Skorohod theorem and identifying the limit  we deduce that  $(\eta^{l})$ and $(\eta^{m})$ converge in distribution to $\eta$ and $\eta'$, which both solve  random PDE \eqref{equation degenerate pde} with the same initial value $\eta_0$. Furthermore, Proposition \ref{lemma unique degenerate} leads to $\eta=\eta'$ $\mP$-a.s.. Therefore, we can deduce by Lemma \ref{lemma gk} below  that  $(\eta^N)_{N \geqslant1}$  converges in probability to the unique solution $\eta$.
\end{proof}

\begin{lemma}[Gy\"ongy and Krylov \cite{gyongy1996existence}]\label{lemma gk}
Let $(Z^N)_{N\geq 1}$ be  a sequence of random elements in a Polish space $E$ equipped with the Borel $\sigma$-algebra. Then $(Z^N)_{N\geq 1}$  converges in probability
		to an $E$-valued random element if and only if for every pair of subsequences $(Z^l)$ and $(Z^m)$ there exists a subsequence $u^k:=(Z^{l(k)},Z^{m(k)})$ converging in distribution to
		a random element $u$ supported on the diagonal $\{(x,y)\in E\times E:x=y\}$.
\end{lemma}

\begin{proof}[Proof of Theorem \ref{thm:2}]
	The proof is similar  to Proposition \ref{prop:gauss}. In addition to the convergence obtained in Theorem \ref{thm:de}, we need to  check the Gaussianity of the unique solution $\eta$ to {\eqref{equation
			degenerate pde}} with Gaussian initial value $\eta_{0}$.
		
		Define the time evolution operators  $\{Q_{s,t}\}_{0\leq s\leq t\leq T}$  by $Q_{s,t}\varphi\assign f(s)$, where $f\in C^1([0,t],C^{\beta'})$, $\beta'\in (d/2+1,\beta+1)$, is the unique solution to \eqref{equation back} with terminal value $\varphi $ at time $t$. Then for each $\varphi\in C^{\infty}$ and $t\in [0,T]$, let $Q_{\cdot,t}\varphi$ play the role of test function in \eqref{equation test}. We find
		\begin{equation*}
			\left\langle \eta_t,\varphi \right\rangle=\left\langle\eta_0,Q_{0,t}\varphi \right\rangle .
		\end{equation*}
		Finally, the result follows by the assumption on $\eta_{0}$.
\end{proof}
\section{Applications}\label{sec:exm}

In this section we finish the proof of Theorem \ref{th:main} and then give a similar result for the particle system \eqref{equation pa} with $C^1$ kernels. At last, we prove a central limit theorem for the initial values, which gives a sufficient condition to {\bf{(A1)}}.

\subsection{Examples}Let us start with proving Theorem \ref{th:main}, which concerns on the most important example of this article: the Biot-Savart law.
\begin{proof}[Proof of Theorem \ref{th:main}]
By our main results Theorem \ref{thm:1}, Proposition \ref{prop:gauss}, and  Theorem \ref{thm:2}, it suffices to check the assumptions {\bf{(A1)-(A5)}}.

 {\bf{(A1)}} is automatically satisfied since  the point vortex model \eqref{eq:in} is of i.i.d initial data, and one can easily check that the Biot-Savart law satisfies the second case of Assumption {\bf{(A2)}}. Moreover, by \cite[Theorem 2]{jabin2018quantitative}, the following condition \eqref{condition1} yields {\bf{(A3)}},
\begin{align}
	\bar{\rho}\in C([0,T],C^3)\quad  \text{ and } \quad \inf \bar{\rho}>0 .\label{condition1}
\end{align}
Now we check  \eqref{condition1}. On one hand, when $\sigma>0$ we deduce  $	\bar{\rho}\in C([0,T],C^3)$  under the assumption that $\bar{\rho}_0\in C^3$ by \cite[Theorem A]{ben1994global}. The fact that $\bar{\rho}\in C([0,T],C^3)$ for the case $\sigma=0$ follows by  \cite[Theorem 2.4.1]{marchioro2012mathematical}.  On the other hand,  $\bar{\rho}\in C([0,T],C^3)$ and Lemma \ref{lemma convolution} implies that $K*\bar{\rho} $ is bounded and Lipschitz continuous, which yields the global well-posedness for the Cauchy problem to the following SDE:
\begin{align}\label{eq:varphi}
\varphi_t=x+\int_0^tK*\bar{\rho}_s(\varphi_s)\mathd s+\sqrt{2\sigma}B_t.
\end{align}
where $B$ is a $d$-dimensional Brownian motion.
Let $\{\varphi_t\}_{t\in [0,T]}$ be the unique solution to \eqref{eq:varphi}, and notice that  $\bar{\rho}$ is the density of the time marginal law of solution to  \eqref{eq:varphi} with initial value $\bar{\rho}_0$-distributed.  

Since  $K$ is divergence free and $\sigma$ is a constant, the flow $\{\varphi_t\}_{t\in [-T,T]}$ ( recall that $\varphi_{-t}:=\varphi_{t}^{-1}$)  is measure preserving (see  \cite[Lemma 4.3.1]{kunita1997stochastic}).  Then for any bounded measurable function $f$,  we have
\begin{align*}
	\left\langle f,\bar{\rho}_t\right\rangle=\mathbb{E}\int_{\mathbb{T}^{d}} f(\varphi_t(x))\bar{\rho}_0(x)\mathd x =\mathbb{E}\int_{\mathbb{T}^{d}} f(x)\bar{\rho}_0(\varphi_{-t}(x))\mathd x,
\end{align*}
which implies $\inf \bar{\rho}>0$ since $\inf \bar{\rho}_0>0$. Thus we obtain \eqref{condition1} and thus obviously  Assumption {\bf (A4)} for the 2D Navier-Stokes equations.   Lastly,  for the 2D Euler equations, we deduce Assumption {\bf{(A5)}} from \cite[Theorem 2.4.1]{marchioro2012mathematical}.
\end{proof}

For general cases, our main result Theorem \ref{thm:1} only requires bounded kernels. However, in order to check Assumption {\bf{(A3)}}  using \cite{jabin2018quantitative} , the extra condtion $\div K\in L^{\infty}$ is necessary. Thus we consider the system \eqref{equation pa} with $C^1$ kernels below, and give a complete result with the only assumption on the  initial data and the confined potential $F$. Nevertheless, the following result considerably relaxes  the condition on kernels in the classical work by  Fernandez and M\'el\'eard \cite{fernandez1997hilbertian}, where the kernel $K$ they considered  should be regular enough,  for instance in  $C^{ 2 + d/2}$.

\begin{Examples}[$C^1$ kernels] Consider the particle system  \eqref{equation pa}  on $\mathbb{T}^d$ and a sequence of independent initial random variables $\{X_i(0)\}_{i\in \mathbb{N}}$ with identical probability density $\bar{\rho}_0$.
Assume that $\sigma>0$, $K\in C^1$, $F,\rho_0\in C^{\beta}$ for some $\beta>2\vee d/2$, and $\inf \bar{\rho}_0>0$.
Then the assumptions {\bf{(A1)-(A4)}} hold.  In particular,   Theorem \ref{thm:1} and Proposition \ref{prop:gauss} hold in this case.
\begin{proof}
	The assumptions {\bf{(A1)-(A2)}}  follow immediately. For simplicity, we set $\sigma=1$ and prove the required regularity for $\bar{\rho}$.  Consider the McKean-Vlasov equation:
	\begin{align}
		\mathd X_t=K*\bar{\rho}_t(X_t)\mathd t+F(X_t)\mathd t+\sqrt{2}\mathd B_t,\quad \bar{\rho}_t=\mathcal{L}(X_t),\label{eqt:mvv}
	\end{align}
where $B$ is a $d$-dimensional Brownian motion. Since $K$ and $F$ are bounded and Lipschitz continuous,
 one can obtain the well-posedness of \eqref{eqt:mvv}, c.f.  \cite[Theorem 3.3]{coghigess2019}. Furthermore,  applying It\^{o}'s formula and superposition principle (see \cite{bogachev2021ambrosio,trevisan2016well}) we have the one-to-one correspondance between the solutions to the Mckean-Vlasov equation \eqref{eqt:mvv} and the solutions to the mean field equation \eqref{equation me}, which implies  the global well-posedness of the mean field equation in the space $C([0,T],\mathcal{P}(\mathbb{R}^d))$.     
Recall that the mild form of the mean-field equation \eqref{equation me} is stated as
\begin{equation*}
	\bar{\rho}_t=\Gamma_{t}*\bar{\rho}_0+\int_0^t\nabla\Gamma_{t-s}*(K*\bar{\rho}_s\bar{\rho}_s+F\bar{\rho}_s)\mathd s,
\end{equation*}
where $\Gamma$ is the heat kernel of $\Delta$. Since $K$ is bounded and Lipschitz continuous, $K*\rho_t$ is  bounded and Lipschitz continuous as well, hence belongs to the Besov space  $B^1_{\infty,\infty}$. Next we consider the following linearized equation 
\begin{equation}
{\rho}_t=\Gamma_{t}*\bar{\rho}_0+\int_0^t\nabla\Gamma_{t-s}*(K*\bar{\rho}_s{\rho}_s+F	{\rho}_s)\mathd s,\label{eqt:frozen}
\end{equation}
where $\bar{\rho}$ is the unique solution to the mean field equation in $C([0,T],\mathcal{P}(\mathbb{R}^d))$. We are going to exploit the regularity of the kernel to improve the regularity of $\bar{\rho}$.  Notice first that $	\Gamma_{\cdot}*\bar{\rho}_0\in C([0,T],C^{\beta})$ and we use Lemma \ref{lemma schau} to have
\begin{align*}
\left\| \int_0^t\nabla\Gamma_{t-s}*(K*\bar{\rho}_s	{\rho}_s+F	{\rho}_s)\mathd s\right\| _{B^1_{\infty,\infty}}^p\lesssim \left(\int_0^t(t-s)^{-\frac{1}{2}}\left\|K*\bar{\rho}_s	{\rho}+F	{\rho}\right\| _{B^1_{\infty,\infty}}\mathd s\right)^p.
\end{align*}
Furthermore,  let $p>2$, then Lemma \ref{lemma triebel} and H\"older's inequality yield that
\begin{align*}
\left\| \int_0^t\nabla\Gamma_{t-s}*(K*\bar{\rho}_s{\rho}_s+F	{\rho}_s)\mathd s\right\| _{B^1_{\infty,\infty}}^p&\lesssim \left(\int_0^t(t-s)^{-\frac{1}{2}}\left\|K*\bar{\rho}_s+F	\right\| _{B^1_{\infty,\infty}}\|{\rho}\|_{B^1_{\infty,\infty}}\mathd s\right)^p
\\&\lesssim_{K,\bar{\rho}, F}\left(\int_0^t(t-s)^{-\frac{p}{2(p-1)}}\mathd s \right)^{\frac{p-1}{p}}\int_0^t \|{\rho}_s\|_{B^1_{\infty,\infty}}^p\mathd s
\\&\lesssim_{K,\bar{\rho}, F}t^{\frac{p-2}{2p}}\int_0^t \|{\rho}_s\|_{B^1_{\infty,\infty}}^p\mathd s.
\end{align*}
The constant omitted here 
depends on $\|F\|_{C^{\beta}}$ and $\sup_{ t \in [0, T]}\|K*\bar{\rho}_s\|_{B^1_{\infty,\infty}}$. Thus we have
\begin{align*}
	\|{\rho}_t\|_{B^1_{\infty,\infty}}^p\lesssim \|\bar{\rho}_0\|_{B^1_{\infty,\infty}}^p+ \int_0^t \|{\rho}_s\|_{B^1_{\infty,\infty}}^p\mathd s,
\end{align*}
for any ${\rho}$ satisfies the linearized equation \eqref{eqt:frozen}. Since the solution $\bar{\rho}$ to the mean field equation also satisfies   \eqref{eqt:frozen}, we find $\bar{\rho}\in C([0,T],B^1_{\infty,\infty})$ by Gronwall's inequality.

 Recall that  we first obtained probability measure-valued solution to \eqref{equation me}. As $\bar{\rho}\in C([0,T],B^1_{\infty,\infty})$, the coefficient $K*\bar{\rho}$ has better regularity, which provides the possibility to improve the regularity of $\bar{\rho}$. In fact, by $K\in B^1_{\infty,\infty}$,  Lemma \ref{lemma triebel} and Lemma \ref{lemma convolution}, we deduce  $K*\bar{\rho}_t\in B^{\alpha+1-\epsilon}_{\infty,\infty}$ whenever $\bar{\rho}\in B^{\alpha}_{\infty,\infty}$, for sufficiently small  $\epsilon>0$. Therefore, $\bar{\rho}$ helps improving the regularity of the coefficient to the linearized equation \eqref{eqt:frozen}. As a result, we could repeat the above estimates with $B^1_{\infty,\infty}$-norm replaced by $B^{2-\epsilon}_{\infty,\infty}$-norm for some $\epsilon>0$ and conclude $\bar{\rho}\in C([0,T],B^{2-\epsilon}_{\infty,\infty})$. We iterate this procedure again and we get $\bar{\rho}\in C([0,T],B^{\beta}_{\infty,\infty})$ for $\beta>2\vee d/2$, which implies Assumption {\bf{(A4)}}.

 As to Assumption {\bf{(A3)}}, by \cite[Theorem 2]{jabin2018quantitative} and $\bar{\rho}\in C([0,T],B^{\beta}_{\infty,\infty})$, it is sufficient to check $\inf\bar{\rho}>0$. Similar to the proof of Theorem \ref{th:main}, we need the  auxiliary SDE
 \begin{align*}
 	\varphi_t=x+\int_0^tK*\bar{\rho}_s(\varphi_s)\mathd s+\int_0^tF(\varphi_s)\mathd s+\sqrt{2}B_t.
 \end{align*}Then for any nonnegative measurable function $f$ on $\mT^d$, we have
\begin{align*}
	\left\langle f,\bar{\rho}_t\right\rangle&=\mathbb{E}\int_{\mathbb{T}^{d}} f(\varphi_t(x))\bar{\rho}_0(x)\mathd x =\mathbb{E}\int_{\mathbb{T}^{d}} f(x)\bar{\rho}_0(\varphi_{-t}(x))|\det\partial\varphi_{-s}(x)|\mathd x,
	\\&\geq\inf\bar{\rho}_0\int_{\mathbb{T}^{d}} f(x)e^{-\int_0^s\|\div(K*\bar{\rho}_r+F)\|_{L^{\infty}}dr}\mathd x
	\\&\geq \inf\bar{\rho}_0e^{-T(\|K\|_{C^1}+\|F\|_{C^1})}\int_{\mathbb{T}^{d}} f(x)\mathd x,
\end{align*}where the first inequality follows by the representation of the determinant for the Jacobian matrix $\det\partial\varphi_{-s}$, see \cite[Lemma 4.3.1]{kunita1997stochastic}.
This implies  {\bf{(A3)}}.
\end{proof}
\end{Examples}

\subsection{A sufficient condition for {\bf{(A1)}}}\label{sec:5.2}
Motivated by the main Assumption {\bf(A3)} in this article,  we give a sufficient condition for central limit theorem for random variables in terms of relative entropy. These random variables may  be  neither independent nor identical distributed. This result could be applied to check  {\bf{(A1)}}.

In the following, we let $\{X_i^N\}_{1\leq i\leq N}$ be a class of random variables with values in  $\mathbb{R}^d$ or $\mathbb{T}^d$, here $\{X_i^N\}_{1\leq i\leq N}$ plays the role of $\{X_i(0)\}_{1\leq i\leq N}$ in {\bf{(A1)}}. More general, the laws of $\{X_i^N\}_{1\leq i\leq N}$ are allowed to be different and these random variables might take values in the whole space $\mR^d$ so that one can compare our result  with the general central limit theorems in the literature. We abuse the notation $\rho_N$ to denote the joint distribution of $\{X_i^N\}_{1\leq i\leq N}$, and let $\bar{\rho}$ denote a probability measure on $\mathbb{R}^d$ (or $\mathbb{T}^d$). For simplicity, we omit  the torus case in the following discussion.

For fixed $\varphi\in \cS(\mR^d)$, define $Y_{\varphi}^N$ by
\begin{align*}
		Y^N_{\varphi}\assign \frac{\sum_{i = 1}^{N}\varphi(X_i^N)-N\left\langle \varphi,\bar{\rho} \right\rangle }{\sqrt{N}}.
\end{align*}
Thus  $\frac{1}{\sqrt{N}}\sum_{i=1}^N(\delta_{X_i^N}-\bar{\rho})$ converges in distribution to the Gaussian variable $\eta_0$ in $\mathcal{S}'(\mathbb{R}^d)$ which satisfies
\begin{equation*}
	\left\langle \eta_0,\varphi\right\rangle \overset{d}{\sim} \mathcal{N}(0,\left\langle \varphi^2,\bar{\rho} \right\rangle-\left\langle \varphi,\bar{\rho}\right\rangle^2  ), \quad \forall\varphi\in \mathcal{S}(\mathbb{R}^d),
\end{equation*}
if and only if the law of $Y_{\varphi}^N$ converges weakly  to $\mathcal{N}(0,\left\langle \varphi^2,\bar{\rho} \right\rangle-\left\langle \varphi,\bar{\rho}\right\rangle^2  )$, which is denoted by  $G_{\varphi}$.

Recall that the bounded Lipschitz distance  $d_{bL}(\mu,\nu)$ between two probability measures $\mu, \nu$ on $\mR^d$ is defined as
\begin{align}\label{df:dis}
d_{bL}(\mu,\nu)=\sup\left\{\int_{\mathbb{R}^d}g(y)\mu(\dif y) -\int_{\mathbb{R}^d}g(y)\nu(\dif y); \quad \|g\|_{L^{\infty}}+\|g\|_{\text{Lip}}\leq 1\right\},
\end{align}
where $\|g\|_{\text{Lip}}$ denotes the Lipschitz constant of $g$. The bounded Lipschitz distance metrizes  the weak convergence (see \cite{villani2008optimalbook}).
 We are going to control the bounded Lipschitz distance between the law of $Y_{\varphi}^N$ and $G_{\varphi}$, denoted by $\mathd_{bL}(\mathcal{L}(Y_{\varphi}^N),G_{\varphi})$,  by the relative entropy $H(\rho_{ N}|\bar{\rho}_N)$ and the bounded Lipschitz distance between $G_{\varphi}$ and $S^N_{\varphi}$, which is defined as
 	\begin{align*}
 S^N_{\varphi}\assign\frac{\sum_{i = 1}^N\varphi(Z_i)-N\left\langle \varphi,\bar{\rho}\right\rangle }{\sqrt{N}},
 \end{align*}
where $\{Z_i\}_{i\in \mathbb{N}}$ is  a class of i.i.d random variables with distribution $\bar{\rho}$.

\begin{proposition}\label{prop:clt}

There is $\lambda>0$ such that for every $\kappa>0$,
	\begin{align*}
		d_{bL}(\mathcal{L}(Y^N_{\varphi}),G_\varphi)\leq \frac{1}{\kappa} H(\rho_{ N}|\bar{\rho}_N) +\frac{\kappa}{2\lambda}+d_{bL}(\mathcal{L}(S^N_{\varphi}),G_\varphi).
	\end{align*}
 By taking $\kappa = \sqrt{2 \lambda H (\rho_N \vert \bar \rho_N )}$, then
	\begin{equation*}
	d_{bL}(\mathcal{L}(Y^N_{\varphi}),G_\varphi)\leq \sqrt{\frac 2 \lambda} \sqrt{H(\rho_N \vert \bar \rho_N )} + d_{bL}(\mathcal{L}(S^N_{\varphi}),G_\varphi).
	\end{equation*}
Assume further
	\begin{equation*}
		H(\rho_N|\bar{\rho}_{ N})\xrightarrow{N\rightarrow \infty} 0,
	\end{equation*}
	then $\mathcal{L}(Y^N_{\varphi})$  converges to $G_\varphi$ weakly, as $N\to \infty$.
\end{proposition}
\begin{proof} From \eqref{df:dis} we know
	\begin{align}
		d_{bL}(\mathcal{L}(Y^N_{\varphi}),G_\varphi)=\sup\left\{\int_{\mathbb{R}^d}g(y)\mathd \mathcal{L}(Y^N_{\varphi}) -\int_{\mathbb{R}^d}g(y)G_{\varphi}(y)\mathd y; \quad \|g\|_{L^{\infty}}+\|g\|_{\text{Lip}}\leq 1\right\}.\label{initial:1}
	\end{align}

	First, we find
	\begin{align*}
		\int_{\mathbb{R}^d}g(y)\mathd \mathcal{L}(Y^N_{\varphi})= &\mathbb{E}g(Y^N_{\varphi})=\mathbb{E}g\left(\frac{\sum_{i = 1}^N\varphi(X_i^N)-N\left\langle \varphi,\bar{\rho} \right\rangle }{\sqrt{N}}\right)
		\\=&\int_{\mathbb{R}^{dN}}g\(\frac{\sum_{i = 1}^N\varphi(x_i)-N\left\langle \varphi,\bar{\rho}\right\rangle }{\sqrt{N}}\)\rho_N\mathd x^N
		\\\assign&\int_{\mathbb{R}^{dN}}g\circ\Phi_N(x^N){\rho}_N\mathd x^N,
	\end{align*}
	with $x^N=(x_1,\dots,x_N)$.
	Then applying the Donsker-Varadhan variational formula \eqref{eq:var} 
	gives
	\begin{align}
		\int_{\mathbb{R}^d}g(y)\mathd \mathcal{L}(Y^N_{\varphi})	\leq \frac{1}{\kappa}\left( H(\rho_{ N}|\bar{\rho}_N)+\log \int_{\mathbb{R}^{dN}} e^{\kappa g\circ \Phi_N(x^N)}\bar{\rho}_N\mathd x^N\right),\label{initial:2}
	\end{align}
	for every $\kappa>0$. Since  $\{Z_i\}_{i\in \mathbb{N}}$ is  a class of i.i.d random variables with distribution $\bar{\rho}$,  we have
	\begin{align}
		\int_{\mathbb{R}^{dN}} e^{\kappa g\circ \Phi_N(x^N)}\bar{\rho}_N\mathd x^N=\mathbb{E}\exp\[\kappa g\(\frac{\sum_{i = 1}^N\varphi(Z_i)-N\left\langle \varphi,\bar{\rho}\right\rangle }{\sqrt{N}}\)\]=\int_{\mathbb{R}}e^{\kappa g(y)}\mathd \mathcal{L}(S^N_{\varphi}).\label{initial:3}
	\end{align}
	Combining \eqref{initial:1}-\eqref{initial:3}, we arrive at
	\begin{align}
		d_{bL}(\mathcal{L}(Y^N_{\varphi}),G_\varphi)\leq& \frac{1}{\kappa} H(\rho_{ N}|\bar{\rho}_N)+d_{bL}(\mathcal{L}(S^N_{\varphi}, G_{\varphi}))\nonumber
		\\&+\sup\left\{\frac{1}{\kappa}\log\int_{\mathbb{R}}e^{\kappa g(y)}\mathd \mathcal{L}(S^N_{\varphi})-\int_{\mathbb{R}}g(y)\mathd \mathcal{L}(S^N_{\varphi}); \quad \|g\|_{L^{\infty}}+\|g\|_{\text{Lip}}\leq 1\right\},\label{initial:4}
	\end{align}
	for every $\kappa>0$.
	
	To handle the last term on the right hand side of \eqref{initial:4}, we need  uniform Gaussian concentration of $\{\mathcal{L}(S^N_{\varphi})\}_{N\in \mathbb{N}}$. That is, there is $a>0$ such that
	\begin{align}\label{eq:a}
		\sup_N\int_{\mathbb{R}}e^{ay^2}\mathd \mathcal{L}(S^N_{\varphi})<\infty.
	\end{align}
	Indeed, recall \eqref{initial:3} and the definition of $\Phi_N$, for \eqref{eq:a} it is sufficient to show that
	\begin{align*}
		&\sup_N	\int_{\mathbb{R}^{dN}} \exp \bigg(aN \bigg| \bigg\langle \varphi,\frac{1}{N}\sum_{i = 1}^N\delta_{x_i}-\bar{\rho} \bigg\rangle\bigg| ^2 \bigg)\bar{\rho}_N\mathd x^N
		\\&=\sup_N	\int_{\mathbb{R}^{dN}} \exp \bigg( aN \bigg| \bigg\langle \varphi-\left\langle \varphi,\bar{\rho}\right\rangle ,\frac{1}{N}\sum_{i = 1}^N\delta_{x_i} \bigg\rangle\bigg| ^2 \bigg)\bar{\rho}_N\mathd x^N<\infty.
	\end{align*}
	This  follows by Lemma \ref{lemma jw} with sufficient small $a$, which depends only on $\|\varphi\|_{L^{\infty}}$. Hence \eqref{eq:a} holds.
	
	Therefore, this  uniform Gaussian concentration \eqref{eq:a} allows us to apply \cite[Theorem 22.10]{villani2008optimalbook}, which   yields that $\{\mathcal{L}(S^N_{\varphi})\}_{N\in \mathbb{N}}$ satisfies the dual formulation of Talagrand-1 inequality uniformly. More precisely, there is $\lambda>0$ such that for any $g\in C_b(\mathbb{R}^d)$, $N\in \mathbb{N}$, and $\kappa\geq 0$,
	\begin{align*}
		\int_{\mathbb{R}}\exp\[{\kappa\inf_{z\in \mathbb{R}}\(g(z)+|y-z|\)}\]\mathd \mathcal{L}(S^N_{\varphi})\leq\exp \(\frac{\kappa^2}{2\lambda}+{\kappa\int_{\mathbb{R}}g(y)\mathd \mathcal{L}(S^N_{\varphi})}\).
	\end{align*}
	Recall \eqref{initial:4}, we are only  interested in functions with Lipschitz constant smaller than 1, for which functions it holds that $\inf_{y\in \mathbb{R}}\(g(z)+|y-z|\)=g(y)$. Thus we have
	\begin{align}
		\frac{1}{\kappa}\log \int_{\mathbb{R}}e^{\kappa g(y)}\mathd \mathcal{L}(S^N_{\varphi})- \int_{\mathbb{R}}g(y )\mathd \mathcal{L}(S^N_{\varphi})\leq \frac{\kappa}{2\lambda}, \quad \forall \kappa>0.\nonumber
	\end{align}
	Taking the above dual formulation of Talagrand-1 inequality into \eqref{initial:4} gives
	\begin{align*}
		d_{bL}(\mathcal{L}(Y^N_{\varphi}),G_\varphi)\leq& \frac{1}{\kappa} H(\rho_{ N}|\bar{\rho}_N)+d_{bL}(\mathcal{L}(S^N_{\varphi}, G_{\varphi}))+\frac{\kappa}{2\lambda},
	\end{align*}
	the result then  follows by the canonical  central limit theorem.
\end{proof}


\br
Similar as the central limit theorems in general, the convergence  still holds if just a fixed number of  elements  are changed. This gives rise to the more reasonable condition
	\begin{align*}
\lim_{m \rightarrow \infty}\lim_{N \rightarrow \infty}		H(\rho_{N-\tau_{N,m}}|\bar{\rho}_{N-m})=0,
	\end{align*}
where $\tau_{N,m}$  is a subset of $\{1,...,N\}$ with $m$ elements and   $\rho_{N-\tau_{N,m}}$ be the joint distribution of $\{X_i^N\}_{i\notin\tau_{N,m}}$. Indeed, define  $Y_{\varphi}^{\tau_{N,m}}$ as
\begin{align*}
	Y^{\tau_{N,m}}_{\varphi}\assign \frac{\sqrt{N}}{\sqrt{N-m}}\( Y^N_{\varphi}- \frac{\sum_{i \in \tau_{N,m}}\varphi(X_i^N)-m\left\langle \varphi,\bar{\rho} \right\rangle }{\sqrt{N}}\).
\end{align*}
It is easy to check that $d_{bL}(\mathcal{L}(Y^N_{\varphi}),\mathcal{L}(Y^{\tau_{N,m}}_{\varphi}))\lesssim m/\sqrt{N}$, which together with  Proposition \ref{prop:clt} applied to $\{X_i^N\}_{i\notin\tau_{N,m}}$ yields the convergence of $\mathcal{L}(Y^N_{\varphi})$  to the Gaussian distribution $G_\varphi$.
\er

\appendix\section{Besov spaces}\label{sec:appa}

In this section we collect useful results related to Besov spaces. Recall that Besov spaces on the torus $B^{\alpha}_{p,q}(\mathbb{T}^d)$ (c.f \cite{triebel2006theory}, \cite{mourrat2017global}), with  $\alpha\in \mathbb{R}$ and $1\leq p,q\leq \infty$,  are defined as the completion of $C^{\infty}$ with respect to the norm
\begin{equation*}
	\|f\|_{B^{\alpha}_{p,q}}\assign \( \sum_{n\geq -1}\left(2^{n\alpha q}\|\mathscr{F}^{-1}(\chi_n\mathscr{F}(f))\|_{L^p(\mathbb{T}^d)}^q \right) \) ^{\frac{1}{q}},
\end{equation*}
where $\mathscr{F}$ represents Fourier transform on $\mathbb{R}^d$ and $\{\chi_n\}_{n\geq -1}:\mathbb{R}^d\rightarrow[0,1]$ are  compact supported smooth functions satisfying
\begin{align*}
	\text{supp}\chi_{-1}\subseteq B(0,\frac{4}{3});\quad \text{supp}\chi_{0}\subseteq B(0,\frac{8}{3})\setminus B(0,\frac{4}{3}),\quad \chi_{n}(\cdot)=\chi_{0}(2^{-n}\cdot) \text{ for } n\geq 0;
	\quad \sum_{n\geq -1} \chi_{n}=1.
\end{align*}
Here $B(0, R)$ denotes the ball of center $0$ and radius $R$.

We collect the following results which are frequently  used  in this article.
\begin{lemma}[{\cite[Proposition4.6]{triebel2006theory}}]\label{lemma embedding}Let $\alpha\in \mathbb{R}$, $\beta\in \mathbb{R}$ and $p_1,p_2,q_1,q_2\in [1,\infty]$. Then the embedding
	\begin{equation*}
		B^{\alpha}_{p_1,q_2}\hookrightarrow B^{\beta}_{p_2,q_2}
	\end{equation*}
is compact if and only if,
\begin{equation*}
	\alpha-\beta>d\(\frac{1}{p_1}-\frac{1}{p_2}\)_{+}.
\end{equation*}
\end{lemma}
\begin{lemma}
	\label{lemma triebel}
	(i) Let $\alpha,
	\beta \in \mathbb{R}$ and $p, p_1, p_2, q \in [1, \infty]$ be such that
	$\frac{1}{p} = \frac{1}{p_1} + \frac{1}{p_2}$. The bilinear map $(u, v)
	\mapsto u \nocomma v$ extends to a continuous map from $B_{p_1, q}^{\alpha}
	\times B_{p_2, q}^{\beta}$ to $B_{p \comma q}^{\alpha \wedge \beta}$ if
	$\alpha + \beta > 0$ (cf. {\cite[Corollary 2]{mourrat2017global}}).
	
	(ii) (Duality.) Let $\alpha\in (0,1)$, $p,q\in[1,\infty]$, $p'$ and $q'$ be their conjugate exponents, respectively. Then the mapping  $(u, v)\mapsto \<u,v\>=\int uv \dif x$  extends to a continuous bilinear form on $B^\alpha_{p,q}\times B^{-\alpha}_{p',q'}$, and one has $|\<u,v\>|\lesssim \|u\|_{B^\alpha_{p,q}}\|v\|_{B^{-\alpha}_{p',q'}}$ (cf.  \cite[Proposition~7]{mourrat2017global}).
\end{lemma}

\begin{lemma}[{{\cite[Corollary
			2.86]{chemin2011fourier}} }]\label{lemma:infity} For any positive real number $\alpha$ and any $p,q\in [1,\infty]$, it holds that
		\begin{align*}
			\|fg\|_{B^{\alpha}_{p,q}}\lesssim\|f\|_{L^{\infty}}\|g\|_{B^{\alpha}_{p,q}}+\|f\|_{B^{\alpha}_{p,q}}\|g\|_{L^{\infty}},
		\end{align*}
		with the proportional constant independent of $f,g$.
\end{lemma}

\begin{lemma}[{\cite[Theorem 2.1 and 2.2]{kuhn2021convolution}}]
	\label{lemma convolution}Let $ \alpha,\beta \in \mathbb{R}, q ,q_1,q_2\in(0, \infty]$ and $p, p_{1}, p_{2} \in[1, \infty]$
	be such that
	$$
	1+\frac{1}{p}=\frac{1}{p_{1}}+\frac{1}{p_{2}},\quad \frac{1}{q}\leq\frac{1}{q_1}+\frac{1}{q_2}.
	$$
	\begin{enumerate}
		\item 	If $f \in B_{p_{1}, q}^{\alpha}$ and $g \in L^{p_{2}},$ then $f * g \in B_{p, q}^{\alpha}$ and
		$$
		\|f * g\|_{B_{p, q}^{\alpha}}\lesssim\| f\|_{B_{p_{1}, q}^{\alpha}} \cdot\left\|g \right\|_{L^{p_{2}}},
		$$
		with the proportional constant independent of $f,g$.
		
		\item	If $f\in B^{\alpha}_{p_1,q_1}$ and $g\in B^{\beta}_{p_2,q_2}$,
		then $f*g\in B^{\alpha+\beta}_{p,q}$ and
		\begin{align*}
			\|f * g\|_{B_{p, q}^{\alpha+\beta}}\lesssim\| f\|_{B_{p_{1}, q_1}^{\alpha}} \cdot\left\|g \right\|_{B^{\beta}_{p_2,q_2}}.
		\end{align*}
	\end{enumerate}

\end{lemma}

Recall the result about smoothing effect of  the heat kernel $\Gamma$.
\begin{lemma}[{\cite[Propositions 3.11, 3.12]{mourrat2017global}}]\label{lemma schau}Let $u\in B^{\alpha}_{p,q}$ for some $\alpha\in \mathbb{R}$, $1\leq p,q\leq \infty$. Then for every $\kappa\geq 0$
	\begin{align*}
		\|\Gamma_t*u\|_{B^{\alpha+2\kappa}_{p,q}}\lesssim t^{-\kappa}\|u\|_{B^{\alpha}_{p,q}},
	\end{align*}
	and
		\begin{align*}
	\|\Gamma_t*u-u\|_{B^{\alpha}_{p,q}}\lesssim t^{\kappa/2}\|u\|_{B^{\alpha+\kappa}_{p,q}}.
	\end{align*}
\end{lemma}

\section{Proof of Lemma \ref{lemma reali}}\label{sec:appb}

In this section we give the proof of Lemma \ref{lemma reali}. First recall the following result from {\cite{flandoli2005stochastic}}.

\begin{lemma}\label{lemma:B2}
	Let $\varphi \rightarrow S (\varphi)$ be a linear continuous mapping from
a separable
Banach space $E$ to $L^0 (\Omega)$ ( random variables with convergence in
probability). Assume that there exists a random variable $C
	(\omega)$ such that for any given $\varphi \in E$ we have
	
	\begin{align}
		| S (\varphi) (\omega) | \leqslant C (\omega) \| \varphi \|_E & \quad
		\tmop{for} \quad\mathbb{P}- \text{a.s.}\quad \omega \in \Omega . \nonumber
	\end{align}
	Then there exists a pathwise realization $\mathcal{S}$ of $S (\varphi)$  from $(\Omega, \mathcal{F}, \mathbb{P})$ to the dual space of $E$  in the sense that
		\begin{align}
		{}[S (\varphi)] (\omega) = & [\mathcal{S} (\omega)] (\varphi), \quad
		\mathbb{P}- a.s, \nonumber
	\end{align}
	for every $\varphi \in E$.
\end{lemma}


\begin{proof}[{Proof of Lemma \ref{lemma reali}}]
	The proof consists of two step. The first step is to find a pathwise
	realization for each $t \in [0, T]$. The second step is justifying the
	pathwise realization forms a progressively measurable process. We denote $\frac{\sqrt{2 \sigma_N}}{\sqrt{N}} $ by $C_N$ below for simplicity.
	
	We first apply {\cite[Lemma 8]{flandoli2005stochastic}} to obtain the
	following equality, for $\varphi
	\in C^{\infty}$,
	
	\begin{align}
		C_N  \sum_{i = 1}^N \int^t_0 \nabla \varphi (X_i) \mathd \nocomma B_s^i =
		& C_N \sum_{i = 1}^N \sum_{k_i \in \mathbb{Z}^d} \langle \nabla \varphi,
		e_{- k_i} \rangle \int^t_0 e_{k_i} (X_i) \mathd B_s^i, \quad \mathbb{P}-a.s.\nonumber
	\end{align} Furthermore, using H\"older's inequality we have
	\begin{align}
		\left| C_N  \sum_{i = 1}^N \int^t_0 \nabla \varphi (X_i) \mathd \nocomma
		B_s^i \right| \leqslant & C_N \sum_{i = 1}^N \left( \sum_{k_i \in
			\mathbb{Z}^{d \nocomma \nocomma}} \langle k_i \rangle^{2 \alpha - 2} |
		\langle \nabla \varphi, e_{- k_i} \rangle |^2 \right)^{\frac{1}{2}}
		 \nonumber\\
		& \times\left( \sum_{k_i \in \mathbb{Z}^{d \nocomma \nocomma}} \langle k_i
		\rangle^{- 2 \alpha + 2} \left| \int^t_0 e_{k_i} (X_i) \mathd B_s^i
		\right|^2 \right)^{\frac{1}{2}} \nonumber\\
		\leqslant & \| \varphi \|_{H^{\alpha}} C_N \sum_{i = 1}^N \left(
		\sum_{k_i \in \mathbb{Z}^{d \nocomma \nocomma}} \langle k_i \rangle^{- 2
			\alpha + 2} \left| \int^t_0 e_{k_i} (X_i) \mathd B_s^i \right|^2
		\right)^{\frac{1}{2}} . \nonumber
	\end{align}
		To apply Lemma \ref{lemma:B2} with $E=H^{\alpha}(\mathbb{T}^d)$ for $\alpha > d / 2 + 1$, it is sufficient to find
		\begin{align}
		& \nocomma \mathbb{E} \left( \sum_{i = 1}^N \left( \sum_{k_i \in
			\mathbb{Z}^{d \nocomma \nocomma}} \langle k_i \rangle^{- 2 \alpha + 2}
		\left| \int^t_0 e_{k_i} (X_i) \mathd B_s^i \right|^2 \right)^{\frac{1}{2}}
		\right)^2 \nonumber\\
		\lesssim_N & \mathbb{E} \left( \sum_{i = 1}^N \sum_{k_i \in \mathbb{Z}^{d
				\nocomma \nocomma}} \langle k_i \rangle^{- 2 \alpha + 2} \left| \int^t_0
		e_{k_i} (X_i) \mathd B_s^i \right|^2 \right) \nonumber\\
		\lesssim_N & \, t \sum_{i = 1}^N \sum_{k_i \in \mathbb{Z}^{d \nocomma
				\nocomma}} \langle k_i \rangle^{- 2 \alpha + 2} < \infty \nocomma .
		\nonumber
	\end{align}
	Therefore, we thus obtain a pathwise realization of $C_N  \sum_{i = 1}^N
	\int^t_0 \nabla \varphi (X_i) \mathd \nocomma B_s^i$ for each $t \in [0,
	T]$, denoted by $\mathcal{M}^N_t$.
	
	Define $\mathcal{M}^N \assign (\mathcal{M}^N_t)_{t \in [0, T]}$. Since the
	stochastic integrals are $t$-continuous, the equality
	
	\begin{align}
		\mathcal{M}^N_t (\varphi) = & C_N  \sum_{i = 1}^N \int^t_0 \nabla \varphi
		(X_i) \mathd \nocomma B_s^i \nonumber
	\end{align}
		holds almost surely for all $t \in [0, T]$ and $\varphi \in C^{\infty}
	(\mathbb{T}^d)$. To justify measurability of $\mathcal{M}^N$. Notice that
	for each $\varphi \in C^{\infty}$, $\langle \mathcal{M}^N, \varphi \rangle=
	\mathcal{M}^N_t(\varphi)$
	is a continuous adapted process. Hence for each $t \in [0, T]$, $\langle
	\mathcal{M}^N, \varphi \rangle : \Omega \times [0, t] \rightarrow
	\mathbb{R}$ 
	is $\mathcal{F}_t \times \mathcal{B} ([0, t])$-measurable. Since
	$C^{\infty}$ is dense in the separable Hilbert space $H^{\alpha}$, using
	Pettis measurability theorem and Lemma \ref{lemma:B2} we thus find $\mathcal{M}^N : \Omega \times
	[0, T] \rightarrow H^{- \alpha}$ is progressively measurable. 
\end{proof}

Data sharing not applicable to this article as no datasets were generated or analysed during the current study.

\bibliographystyle{alpha}
\bibliography{fluctuationvortex}

\end{document}